\documentclass[aps,pra,groupedaddress,amssymb,amsmath,showpacs]{revtex4-1}
\usepackage{dsfont}
\usepackage{bm}
\usepackage{bbm}
\usepackage[mathscr]{euscript}
\usepackage{amsmath}
\usepackage{amssymb}
\usepackage{amsthm}
\usepackage{amsfonts}
\usepackage{graphicx}
\usepackage{color}
\usepackage[colorlinks=true,urlcolor=blue,linkcolor=blue,citecolor=red]{hyperref}

\begin{document}


\title{A Complete Method of Comparative Statics for Optimization Problems \\ \centerline{(Unabbreviated Version)}}


\author{M. Hossein Partovi}
\email[Electronic address: ]{hpartovi@csus.edu}
\affiliation{Department of Physics and Astronomy, California State
University, Sacramento, California 95819-6041}
\author{Michael R. Caputo}
\email[Electronic address: ]{mcaputo@bus.ucf.edu}
\affiliation{Department of Economics University of Central Florida,
P.O. Box 161400, Orlando, Florida 32816-1400 }



\begin{abstract}
A new method of deriving comparative statics information using
generalized compensated derivatives is presented which yields
\textit{constraint-free} semidefiniteness results for any
differentiable, constrained optimization problem.  More generally,
it applies to any differentiable system governed by an extremum
principle, be it a physical system subject to the minimum action
principle, the equilibrium point of a game theoretical problem
expressible as an extremum, or a problem of decision theory with
incomplete information treated by the maximum entropy principle.
The method of generalized compensated derivatives is natural and
powerful, and its underlying structure has a simple and intuitively
appealing geometric interpretation.  Several extensions of the main
theorem such as envelope relations, symmetry properties and
invariance conditions, transformations of decision variables and
parameters, degrees of arbitrariness in the choice of comparative
statics results, and rank relations and inequalities are developed.
These extensions lend considerable power and versatility to the new
method as they broaden the range of available analytical tools.
Moreover, the relationship of the new method to existing
formulations is established, thereby providing a unification of the
main differential comparative statics methods currently in use.  A
second, more general theorem is also established which yields
exhaustive, constraint-free comparative statics results for a
general, constrained optimization problem.  Because of its universal
and maximal nature, this theorem subsumes all other comparative
statics formulations and as such represents a significant
theoretical result.  The method of generalized compensated
derivatives is illustrated with a variety of models, some well
known, such as profit and utility maximization, where several novel
extensions and results are derived, and some new, such as the
principal-agent problem, the efficient portfolio problem, a model of
a consumer with market power, and a cost-constrained profit
maximization model, where a number of interesting new results are
derived and interpreted.  The large arbitrariness in the choice of
generalized compensated derivatives and the associated comparative
statics results is explored and contrasted to the unique eigenvalue
spectrum whence all comparative statics results originate.  The
significance of this freedom in facilitating empirical verification
and hypothesis testing is emphasized and demonstrated by means of
familiar economic models.  Conceptual clarity and intuitive
understanding are emphasized throughout the paper, and particular
attention is paid to the underlying geometry of the problems being
analyzed.

\end{abstract}

\maketitle



\section{INTRODUCTION}

\subsection{R\'{e}sum\'{e} of Current Methods}

Ever since Cournot (1838) proved that a profit maximizing monopoly
would reduce its output and raise its price in response to an
increase in the unit tax on its output, economic theorists have
sought to examine the comparative statics properties of economic
theories. The epitome of such comparative statics occurred in 1886
with Antonelli and in 1915 with Slutsky, when they independently
showed that the comparative statics properties of the archetype
utility maximization model are summarized in a negative semidefinite
matrix, now known as the Antonelli and Slutsky matrices,
respectively. But it was not until the publication of Samuelson's
(1947) dissertation that comparative statics became part of the
economist's tool kit. Samuelson (1947), building upon the pioneering
work of Allen (1938), Hicks (1939), and others, formulated the
following strategy for deriving comparative statics results
associated with optimization problems: (i) assume the second-order
sufficient conditions hold at the optimal solution and apply the
implicit function theorem to the first-order necessary conditions to
characterize the optimal choice functions, (ii) differentiate the
identity form of the first-order necessary conditions with respect
to the parameter of interest, and (iii) solve the resulting linear
system of comparative static equations.  This \textit{primal}
approach is indicative of most comparative statics analyses in
economics to this day, and provides the general framework within
which subsequent generalizations and refinements have been carried
out.  Indeed, despite several extensions of Samuelson's (1947) basic
analysis by subsequent authors, the treatment of comparative statics
in textbooks and monographs (with the exception of Silberberg's
[1974] contributions, discussed below), as well as in most research
work, essentially relies on his original framework and for the most
part ignores the subsequent developments.  This circumstance is
partly a result of the desire for economy of presentation.  But it
is undoubtedly also a consequence of the fact that none of the later
extensions and generalizations of Samuelson's (1947) basic framework
succeeded in arriving at a completely general, constraint-free
semidefiniteness result involving the partial derivatives of the
decision variables with respect to the parameters for the case of a
constrained optimization problem.  In fact the existing literature
does not really contemplate the general existence of such an
unrestricted result [see Silberberg, (1990)].

Our objective in this paper is to show that such a general result in
fact exists, that the formalism for its construction is natural and
grounded in intuition, and that the comparative statics results
derived from it are sufficiently broad and powerful as to render it
a worthwhile tool accessible to most economists.  Indeed only one
extra conceptual ingredient beyond Samuelson's (1947) basic
framework is needed in its construction, namely a clear
understanding of the role of compensated derivatives in formulating
a general comparative statics statement.  This understanding in turn
leads to a recipe for identifying a suitable class of compensated
derivatives in terms of which the unrestricted semidefiniteness
result for a general optimization problem is realized. Before
delving into details, however, it is appropriate to briefly review a
few highlights of the progress since Samuelson's (1947) seminal
work, partly with a view to the concepts and methods presented in
this paper. It is also appropriate at this point to recall that the
standard definition of \textit{definiteness} or
\textit{semidefiniteness} for a matrix, to which we shall adhere,
requires \textit{symmetry} as a prerequisite (see \S IIB for a
definition).

Some twenty seven years after its inception, Samuelson's (1947)
primal method was significantly advanced and enriched by a new
formulation of the problem.  Silberberg (1974), building upon the
work of Samuelson (1965), brought about the change by means of a
clever construction that in retrospect seems quite natural.  By
changing one's viewpoint of the optimization problem, that is, by
simultaneously considering the set of parameters and choice
variables of the original, primal problem as choice variables,
Silberberg (1974) set up an excess optimization problem that
simultaneously contained the primal optimization problem together
with a dual optimization problem, hence the appellation
\textit{primal-dual}.  Optimization of the primal-dual problem with
respect to the original choice variables yields the primal
optimality conditions, while optimization with respect to the
parameters yields the envelope properties of the primal optimization
problem as the first-order necessary conditions, and the fundamental
comparative statics properties of the primal optimization problem as
the second-order necessary conditions. Undoubtedly, the most
important feat achieved by the introduction of the primal-dual
method was the general construction of a semidefinite matrix that
contains the comparative statics properties of the primal
optimization problem.  The most significant shortcoming of this
method, on the other hand, is the fact that its final result is
subject to constraints if the optimization problem is constrained
and the constraint functions depend on the parameters.

A few years later, Hatta (1980) introduced the \textit{gain} method
to deal with a class of optimization problems that are essentially a
nonlinear, multi-constraint generalization of the Slutsky-Hicks
(basic utility maximization) problem.  For unconstrained problems,
the gain method is identical to the primal-dual method of Silberberg
(1974).  For constrained problems, on the other hand, the gain
method yields constraint-free comparative statics results for the
above-mentioned class of problems. This is accomplished by a
procedure that amounts to applying the compensation scheme used in
the standard Slutsky-Hicks problem (cf. \S IIC), albeit in a
generalized form.  Similarly, the method avoids the use of Lagrange
multipliers by a comparative scheme reminiscent of Silberberg's
(1974) primal-dual method, and in effect converts the constrained
problem into an unconstrained one for the gain function.  This
procedure is in essence equivalent to the standard multiplier
method, differing only in the manner in which the auxiliary
functions are introduced.  Otherwise, the flavor of the analysis
with the gain method is similar to that of Silberberg's (1974)
primal-dual method.  Hatta's (1980) work represents a significant
advance since it succeeded in overcoming an important shortcoming in
the primal-dual method by deriving the first complete,
constraint-free, comparative statics results for the above-mentioned
class of problems.  However, it has not spurred further progress in
the subject, nor has it gained wide acceptance by workers in the
field, primarily because its methods are not sufficiently general
and lack a clear conceptual basis.

The methods described above, and indeed most any other comparative
statics analysis, rely upon the differentiability of the primitives
and solution functions of the optimization problem.  The recent work
of Milgrom and Shannon (1994) and Milgrom and Roberts (1994), on the
other hand, breaks with that tradition of dealing with the local
extrema of differentiable functions.  These authors, in contrast,
develop an \textit{ordinal} method for comparative statics, where
ordinality implies invariance with respect to order-preserving
transformations.  Milgrom and Roberts (1994) state the following
properties in defense of the ordinal method: (i) it dispenses with
most of the smoothness assumptions of the traditional methods, (ii)
it is capable of dealing with multiple equilibria and finite
parameter changes, and (iii) it includes a theory of the robustness
of conclusions with respect to assumptions.  It is clear that the
ordinal method is intended to handle complications which lie outside
the scope of differential comparative statics methodology.  As such
the ordinal method is an essentially different, and in certain ways
a complementary, approach to comparative statics.

Before closing this section, it is appropriate to briefly review the
use of
compensated derivatives in differentiable comparative statics analysis.
The interest in compensated comparative statics properties of economic
models has its genesis in the Slutsky matrix of compensated derivatives
of
the Marshallian (or ordinary) demand functions.  Research on compensated
comparative statics properties of general optimization problems,
however,
is of a more recent origin.  The best known contribution of this ilk
is a set of three papers by Kalman and Intriligator (1973), and
Chichilnisky and
Kalman
(1977, 1978), which introduced generalizations that actually predated the
methods described above, albeit within a restricted framework.
In particular, these authors emphasized the significance of compensated
derivatives
in the context of a general class of constrained optimization problems
and
established the existence of a generalized Slutsky-Hicks matrix for such
problems.
However,
they did not succeed in establishing the crucial semidefiniteness
properties of this matrix in general.  Perhaps because their analysis
was
primarily concerned with establishing the existence of solutions using
essentially primal methods, and because their comparative statics
results
were restricted to special forms, their work was largely superseded by
subsequent developments.  Similarly, their construction
and use of compensated derivatives, although a significant advance
toward
a general method of dealing with constraints,
was rendered in the same limited context and was not much pursued by
others.

The last contribution to be mentioned here is the work of Houthakker
(1951-52) which, while perhaps not as well recognized as the other
post-Samuelsonian contributions mentioned above, actually predates
them. In an attempt to quantify the role of \textit{quality} in
consumer demand, Houthakker (1951-52) clearly recognized the
important role played by compensated derivatives, the large number
of ways in which they can be constructed, and how they are related
to a differential characterization of the constraints present in the
problem, albeit in the context of specific examples.  However, he
only succeeded in deriving the desired semidefiniteness condition
for a restricted class of problems, and his contribution did not
lead to any significant development in the subject.

\subsection{Generalized Compensated Derivatives}

As mentioned before, the method presented here is based on the
pivotal role played by compensated derivatives, and in effect
unifies and generalizes the above-described formulations.  As is
often the case, this generalization leads to a conceptually simpler
structure, while yielding a powerful method of deriving
constraint-free comparative statics results. The basic idea
originates in the observation that the natural parameters of a given
optimization problem are not necessarily the most advantageous
variables for formulating comparative statics results. This is
plainly obvious in the Slutsky-Hicks problem where a linear
combination of derivatives in the form of a compensated derivative
must be employed in order to obtain the desired semidefiniteness
properties.  On the other hand, a linear combination of partial
derivatives is, aside from an inessential scale factor, simply a
\textit{directional} derivative pointed in some direction in
parameter space. Since an uncompensated (i.e., a partial) derivative
is also a directional derivative, albeit a very particular one
pointed along one of the coordinate axes in parameter space, it
follows that the distinction between the two is merely a matter of
the choice of coordinates in parameter space and has no intrinsic
standing. Indeed a \textit{local} rotation in parameter space, i.e.,
a rotation whose magnitude and direction may vary from point to
point, can interchange the role of the two.  Such a rotation is of
course equivalent to adopting a new set of parameters for the
optimization problem, suitably constructed as functions of the
original ones. Clearly then, any general formulation of differential
comparative statics must consider the possibility of choosing from
this vastly enlarged class of directional derivatives in parameter
space.

The question that arises at this juncture is whether there exists a
universal criterion for choosing the compensated derivatives so as
to guarantee the desired semidefiniteness properties free of
constraints, and without requiring any restriction on the structure
of the optimization problem.  Remarkably, and in retrospect not
surprisingly, there is a simple and natural answer to this question.
One simply chooses the compensated derivatives in conformity with
the constraint conditions, i.e., along directions that are tangent
to the level surfaces of all the constraint functions at each point
of the parameter space (for a given value of the decision
variables). Equivalently, acting on the constraint functions, the
compensated derivatives are required to return zero at all points of
parameter space (for a given value of the decision variables).  This
requirement is in effect the \textit{ab initio}, differential
implementation of the constraints, as will become clear in the
course of the following analysis.  What is perhaps more remarkable,
however, is that this procedure yields an operationally powerful,
yet intuitively palatable framework that is capable of yielding
varied and novel results, even in the case of unconstrained
problems.

Before embarking upon the following development, it is worth identifying
its main elements.  The cornerstone of the analysis is Lemma 1 in \S IIA,
the main mathematical device is the generalized compensated derivative,
also developed in \S IIA, and the principal result of the analysis is
Theorem 1 in \S IIB.  This is followed by several auxiliary developments
and extensions, most of which are summarized in Theorem 2 and the
Corollaries in \S IIB and IIC.  In \S IID we establish a maximal extension
of Theorem 1, a theoretically important result which implies all other
comparative statics formulations.  The method of generalized compensated
derivatives is applied in \S III to a number of models, some old and
several new, where its power and flavor are illustrated.

\section{DEVELOPMENT OF THE MAIN RESULTS}

This section is devoted to a detailed development of the main tools and
results of the present method according to the ideas described above.
Because the underlying structure of this method is simple and natural, we
find it worthwhile to continue emphasizing the intuitive
aspects of the analysis in the following development.  It is our hope that
a clear conceptual grasp of its basics will encourage its use among a
wider segment of economists.

\subsection{Geometry of Generalized Compensation}

To convey a clear picture of how the ideas described in \S IB are
implemented, we find it helpful to describe the geometrical
background involved in the construction of compensated derivatives
in some detail.  Let us first recall certain basic results from
analysis.  Consider a real-valued, continuously differentiable
function $\phi(\mathbf{x},\mathbf{a})$ defined for $\bf a$ in some
open subset ${\cal P}$ of ${\Re}^{N}$ and $\bf x$ in a subset ${\cal
D}$ of ${\Re}^{M}$.  This function will later be identified with an
objective or constraint function, with the $M$-dimensional vector
$\mathbf{x}= (x_{1}, x_{2},\ldots,x_{M})$ representing the decision
variables in the \textit{decision space}, and the $N$-dimensional
vector $\mathbf{a}= (a_{1}, a_{2},\ldots,a_{N})$ representing a set
of $N$ parameters in the \textit{parameter space} .  For most of
this section, the dependence of $\phi$ on $\bf x$ plays a secondary
role in the discussion, so that it is useful to consider $\bf x$ as
fixed and $\phi$ as a continuously differentiable function defined
on an open subset $\cal P$ of $\Re ^{N}$.  Throughout the paper, we
will symbolize vectors and matrices both in ``vector'' and
``matrix'' notation, as in $\bf v$ and $\sf{M}$, and in
``component'' notation, as in ${v}_{i}$ and ${\sf{M}}_{ij}$.  The
inner product of two vectors $\bf v$ and $\bf u$ is denoted by
$\mathbf{v}\cdot\mathbf{u} \stackrel{\rm{ def}}{=}
{\sum}_{i=1}^{M}{v}_{i}{u}_{i}$.  For a pair of matrices $\sf{M}$
and $\sf{N}$, the product is denoted by ${\sf{M}}{\sf{N}}
\stackrel{\rm{ def}}{=} {\sum}_{j=1}^{M}{\sf{M}}_{ij}{\sf{N}}_{jk}$.
As for inner products between vectors and matrices,
${\sf{M}}\mathbf{v}$ stands for the vector (column matrix)
${\sum}_{j=1}^{M}{\sf{M}}_{ij}{v}_{j}$, $\mathbf{v}^{\dag}{\sf{M}}$
stands for the vector (row matrix)
${\sum}_{j=1}^{M}{v}_{j}{\sf{M}}_{ji}$, so that
$\mathbf{v}^{\dag}{\sf{M}}\mathbf{u}$ stands for the scalar
${\sum}_{i,j=1}^{M}{v}_{i}{\sf{M}}_{ij}{u}_{j}$.

\begin{figure}
\includegraphics[]{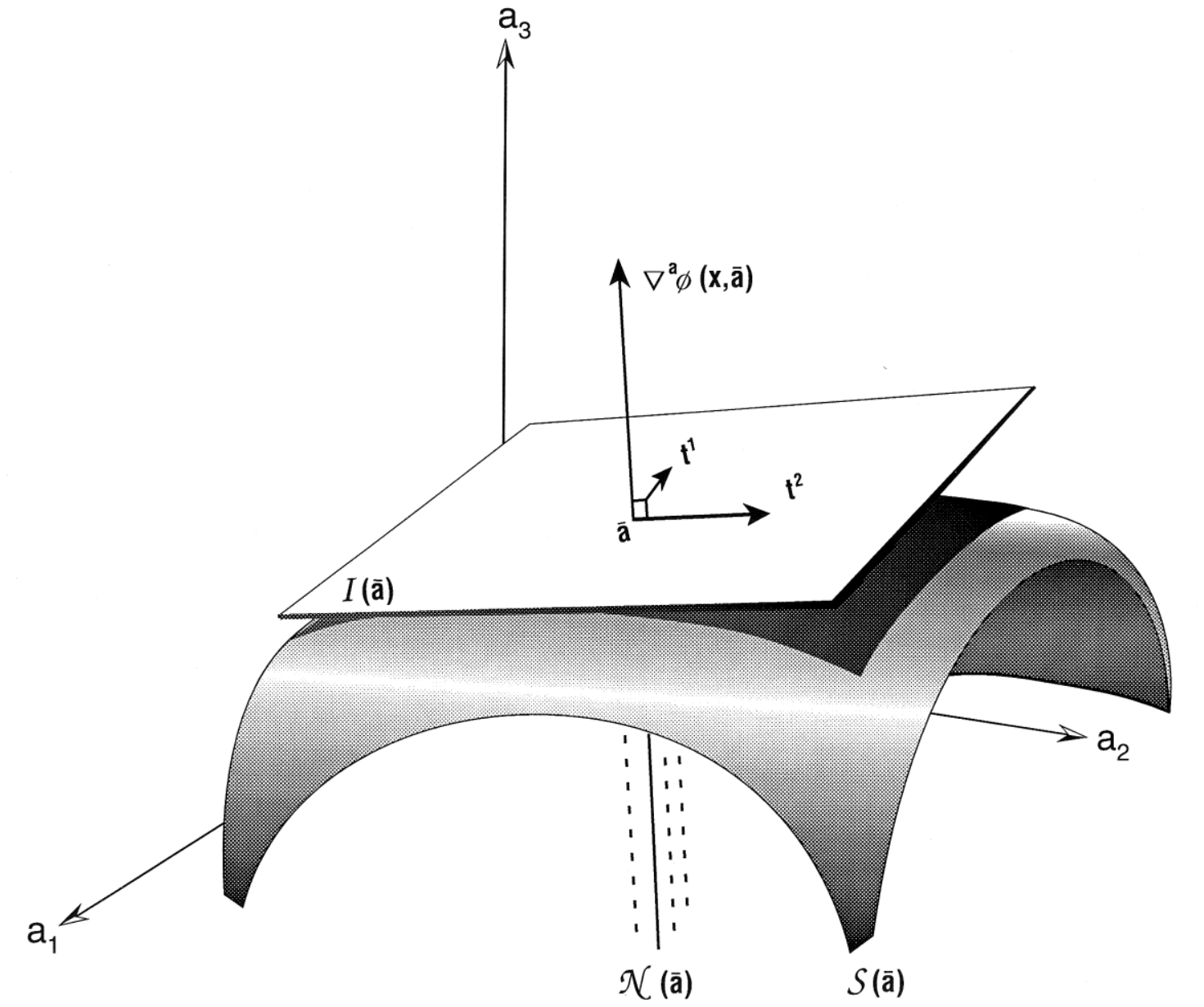}
\caption{An illustration of the tangent plane and normal direction
in a 3-dimensional parameter space.  Isovectors $\mathbf{t}^{1}$ and
$\mathbf{t}^{2}$ are a pair of independent vectors on the tangent
plane.} \label{fig1}
\end{figure}

Given a fixed value of $\mathbf{ x}$, let $\bar{\mathbf{a}}$ be a
point in $\cal P$ and ${\cal S}(\bar{\mathbf{{a}}})$ the
$(N-1)$-dimensional level surface of $\phi$ passing through
$\bar{\mathbf{{a}}}$, that is, the set of points in $\cal P$ at
which $\phi$ assumes the fixed value
$\phi(\mathbf{x},\bar{\mathbf{{a}}})$.  In symbols, ${\cal
S}(\bar{\mathbf{{a}}}) \stackrel{\rm{ {def}}}{=} \{ \mathbf{a} \in
{\cal P}:
\phi(\mathbf{x},\mathbf{a})=\phi(\mathbf{x},\bar{\mathbf{{a}}}) \}
$. Then the gradient of $\phi$ with respect to $\mathbf{a}$ at point
$\bar{\mathbf{a}}$, denoted by $\mathbf{\nabla}^{\mathbf{a}}
\phi(\mathbf{x},\bar{\mathbf{a}})$, is a vector normal to ${\cal
S}(\bar{\mathbf{a}})$ at $\bar{\mathbf{a}}$, and represents the
direction of most rapid change for $\phi$ at that point.  We shall
denote this vector by
$\mathbf{n}(\bar{\mathbf{a}})\stackrel{\rm{{def}}}{=}\mathbf{\nabla}^\mathbf{a}
\phi(\mathbf{x},\bar{\mathbf{a}})$ and its normalized version by
$\hat{\mathbf{n}}(\bar{\mathbf{a}})
\stackrel{\rm{{def}}}{=}{\mathbf{n}}(\bar{\mathbf{a}}) /
\|{\mathbf{n}}(\bar{\mathbf{a}}) \|$, where $\|\mathbf{z}\|$ denotes
the length of the vector $\bf z$. On the other hand, the hyperplane
tangent to the level surface at $\bar{\mathbf{a}}$, denoted by
${\cal I}(\bar{\mathbf{a}})$, is generated by the set of vectors
that are tangent to ${\cal S}(\bar{\mathbf{a}})$ at point
$\bar{\mathbf{a}}$, and represents the directions of no change for
$\phi$ at $\bar{\mathbf{a}}$. Note that
${\mathbf{n}}(\bar{\mathbf{a}})$ is orthogonal to ${\cal I}
(\bar{\mathbf{a}})$.  In plain language, then, the normal vector
${\mathbf{n}}(\bar{\mathbf{a}})$ represents the direction of
\textit{maximum} change, whereas the tangent hyperplane represents
the directions of zero change, or the \textit{null} directions.  We
shall refer to a vector in the tangent hyperplane as an
\textit{isovector}.  Thus an isovector is any vector that points in
a null direction.  Together, the isovectors and the normal vector
span all possible directions, thus providing a convenient local
vector space at point $\bar{\mathbf{a}}$ of the parameter space.
Figure~\ref{fig1} is an illustration of the structure just described
in a three-dimensional parameter space, with $\mathbf{t}^{1}$ and
$\mathbf{t}^{2}$ depicting two isovectors in the tangent hyperplane
${\cal I} (\bar{\mathbf{a}})$.

To define the geometric structure described above more generally, we
first consider the vector space ${\cal N} (\bar{\mathbf{{a}}})$
generated by the normal vector ${\mathbf{n}}(\bar{\mathbf{{a}}})$.
While ${\cal N} (\bar{\mathbf{{a}}})$ is a one-dimensional space in
the present case, it will be generalized to higher dimensions in
more general situations later, and can therefore be referred to as
the \textit{normal hyperplane} in general.  Thus we have associated
with each point $\bar{\mathbf{{a}}}$ of the parameter space a pair
of orthogonal vector spaces ${\cal N} (\bar{\mathbf{{a}}})$ and
${\cal I} (\bar{\mathbf{{a}}})$, the direct sum of which is the
$N$-dimensional local vector space referred to above.  We will
denote the latter by ${\cal T}^{N}(\bar{\mathbf{{a}}})$, and refer
to it as the \textit{tangent space} associated with point
$\bar{\mathbf{{a}}}$. This construction implies, among other things,
that any real, $N$-dimensional vector can be uniquely resolved into
a pair of components, one lying in the normal hyperplane and the
other in the tangent hyperplane.

If we now consider all points of $\cal P$ as endowed with the
tangent space structure just described, there emerges a
configuration of local vector spaces covering all of $\cal P$.  In
terms of these local spaces, the differential structure of $\phi$ is
essentially reduced to that of a function of one variable, namely
the coordinate measured along the normal direction, since at a given
point $\mathbf{a} \in {\cal P}$, the derivative of $\phi$ along any
direction lying in the tangent hyperplane vanishes.  We may
summarize these results by saying that at every point $\mathbf{a}
\in {\cal P}$, and for any direction specified by the (real,
$N$-dimensional) unit vector $\hat{\mathbf{u}}$, (i) the
decomposition $\hat{\mathbf{u}}=[\hat{\mathbf{u}}\cdot
\hat{\mathbf{n}}(\mathbf{a})] \hat{\mathbf{n}} (\mathbf{a}) +
\mathbf{u}^{t}$, where ${\bf{u}}^{t}$ represents the projection of
$\hat{\mathbf{u}}$ onto the tangent hyperplane ${\cal
I}(\mathbf{a})$, is unique, and (ii) the directional derivative of
$\phi$ at point $\mathbf{a}$ in the direction $\hat{\mathbf{u}}$ is
given by $D^{\hat{\mathbf{u}}} \phi
(\mathbf{x},\mathbf{a})\stackrel{\rm{{def}}}{=}\hat{\mathbf{u}}
\cdot \mathbf{\nabla}^\mathbf{a} \phi (\mathbf{x},\mathbf{a}) =
\hat{\mathbf{u}} \cdot {\mathbf{n}}(\mathbf{a}) $. The last property
implies that if $\hat{\mathbf{u}}$ points in a null direction, i.e.,
if $\hat{\mathbf{u}} \in {\cal I}(\mathbf{a})$, then the
corresponding derivative vanishes: $D^{\hat{\mathbf{u}}} \phi
(\mathbf{x},\mathbf{a})=0$.  We shall refer to this condition as the
\textit{null property} of a directional derivative. Note that we
have suppressed the dependence on $\mathbf{{a}}$ in $D^{\mathbf{u}}$
to avoid cluttered notation.  Note also that the length of
$\hat{\mathbf{u}}$ plays no role with respect to the null property,
only the fact that it points in a null direction, or, equivalently,
that it is an isovector.  Thus if $\bf t$ is an isovector of any
length, then $\mathbf{t} \cdot \mathbf{\nabla}^\mathbf{a} \phi
(\mathbf{x},\mathbf{a}) = 0$, so that $\mathbf{t} \cdot
\mathbf{\nabla}^\mathbf{a}$ possesses the null property as well.
Observe that $\mathbf{t} \cdot \mathbf{\nabla}^\mathbf{a}$ is a
linear combination of partial derivatives, equal to $\|\mathbf{t} \|
D^{\hat{\mathbf{t}}}$.  Since our primary focus is on the null
property of this linear combination, we shall continue to (loosely)
refer to $\mathbf{t} \cdot \mathbf{\nabla}^\mathbf{a}$ as a
\textit{directional derivative}, overlooking the inessential scalar
factor $\|\mathbf{t} \|$ in doing so.  We have thus arrived at the
rather obvious conclusion that any directional derivative of a
function in the direction of one of its isovectors has the null
property with respect to that function.

This last property is the basis of our definition of compensated
derivatives: Any linear combination of partial derivatives
possessing the null property with respect to a function $\phi$ is a
\textit{generalized compensated derivative} (hereafter abbreviated
as GCD) with respect to that function.  To avoid trivialities, we
shall exclude identically vanishing linear combinations from this
definition.  Thus any directional derivative in a null direction is
a GCD, with the converse also holding except where
$\mathbf{\nabla}^\mathbf{a} \phi (\mathbf{x},\mathbf{a})$ vanishes.
As an example, consider the case illustrated in Fig.~\ref{fig1}
where $\mathbf{t}^{1}$ and $\mathbf{t}^{2}$ depict a pair of
isovectors at point $\bar{\mathbf{a}}$, which is located on the
level surface ${\cal S}(\bar{\mathbf{a}})$ of
$\phi(\mathbf{x},{\mathbf{a}})$. For this case, then, the
directional derivatives $\mathbf{t}^{1} \cdot
\mathbf{\nabla}^{\mathbf{a}} $ and $\mathbf{t}^{2} \cdot
\mathbf{\nabla}^{\mathbf{a}}$ are a pair of GCD's with respect to
$\phi(\mathbf{x},{\mathbf{a}})$. Their null property, on the other
hand, is seen in Fig.~\ref{fig1}to reflect the simple fact that the
rate of change of a function is zero in directions tangent to its
level surface. Moreover, the reason for making the null property a
defining characteristic is the crucial fact that, when the GCD's
possess this property with respect to the constraint functions, the
resulting semidefiniteness results emerge free of constraints.  This
basic result will be established in this and the following sections.

Although our construction has so far been limited to a
one-dimensional normal hyperplane, we shall consider
higher-dimensional normal hyperplanes in dealing with illustrative
problems, where expressions involving the normal vectors
$\mathbf{\nabla}^\mathbf{a} \phi (\mathbf{x},\mathbf{a})$ and
$\hat{\mathbf{n}}$ will be generalized to include a set of such
vectors. Otherwise, the geometrical features discussed above,
including the decomposition into normal and tangent components,
remain unchanged.

Let us illustrate the construction of GCD's by means of four basic
economic problems.  All but the third assume perfect information and
are thus deterministic, while the third problem involves uncertainty
as well as asymmetry of information between the transacting parties.
In the first problem, the function $\phi$ is given by
$\phi(\mathbf{x},\mathbf{a})\stackrel{\rm{{def}}}{=}m-\mathbf{p
\cdot x}$, where $\bf p$, $\mathbf{x} \in \Re_{+}^{N-1}$, and
$\mathbf{x}$ is a fixed vector here. This is of course the generic
form of the budget constraint that appears in the prototype utility
maximization problem.  Here the variables of interest are the $N$
parameters appearing in $\phi$, identified according to $\mathbf{a}
=
(a_{1},\ldots,a_{N})\stackrel{\rm{{def}}}{=}(p_{1},\ldots,p_{N-1},m)=
(\mathbf{p},m)$.  Then the normal direction is given by ${
\mathbf{\nabla}^\mathbf{a}}\phi(\mathbf{x},\mathbf{a})=(-\mathbf{x},1)$.
Next, we must choose a set of $N-1$ vectors, all orthogonal to the
normal direction just calculated, as a basis for the tangent
hyperplane.  A convenient choice is the set $\mathbf{t}^{\alpha}
\stackrel{\rm{{def}}}{=} (0,\ldots,1,\ldots,0,x_{\alpha})$,
$\alpha=1,\ldots,N-1$, where the component equal to unity occupies
the $\alpha$th position within the parenthesis.  Note that these
basis vectors are neither mutually orthogonal nor of unit length,
but that they are all orthogonal to the gradient vector and
constitute an independent set in the tangent hyperplane.  Clearly,
each $\mathbf{t}^{\alpha}$ is an isovector, so that the set of
directional derivatives formed by taking the inner product
$\mathbf{t}^{\alpha} \cdot {
\mathbf{\nabla}^\mathbf{a}}\stackrel{\rm{{def}}}{=}D_{\alpha}(\mathbf{x},\mathbf{a})$
possesses the desired null property with respect to $\phi$ and is
therefore a set of $N-1$ GCD's.  Note that we have made the
dependence of the GCD's on the points $\bf x$ of the decision space
and $\bf a$ of parameter space explicit in our notation, a practice
which we shall henceforth follow. Remembering that here ${
\mathbf{\nabla}^\mathbf{a}}=({\partial \over \partial
{p}_{1}},\ldots, {\partial \over \partial {p}_{N-1}}, {\partial
\over \partial m})$, one can easily derive the result
$D_{\alpha}(\mathbf{x},\mathbf{a})={\partial \over \partial
{p}_{\alpha}} + {x}_{\alpha} {\partial \over \partial m}$.  The
latter is of course the compensated derivative that appears in the
Slutsky-Hicks problem.  It is worth noting here that (for $N\geq3$)
any linear combination of $D_{\alpha}(\mathbf{x},\mathbf{a})$ will
again yield a GCD. This fact corresponds to the infinite number of
ways one can choose $N-1$ independent basis vectors in ${\cal I}(\bf
{a})$ when $N \geq 3$, and already gives a hint of the generality of
the present method, and as well of the diversity of the results it
can generate. It is also worth noting here that the procedure
described above and applied to this problem is not the only
available method for constructing GCD's, but one that reveals the
underlying geometry most explicitly.  Thus any procedure that yields
a set of directional derivatives with the null property and such
that the corresponding directions in the parameter space span the
tangent hyperplane will suffice. Taking advantage of this freedom,
we will in the following sections introduce simplified procedures
for constructing the set of GCD's for various applications.

While the above treatment illustrates the construction of a GCD for a
constrained problem in the context of a well known
model, our second problem will deal with an unconstrained model and
will give an indication of the calculational novelty of the present
method.

Consider ${\tilde{\phi}} (\mathbf{x}, \mathbf{a})
\stackrel{\rm{{def}}}{=} s[pF(\mathbf{x})-\mathbf{x} \cdot
\mathbf{w}]$, where $\mathbf{x , w} \in {\Re}_{+}^{M}$.  This
$\tilde{\phi}$ is, except for the presence of the \textit{scale
factor} $s > 0$, the familiar function that describes the profit of
a firm producing a single output from $M$ inputs under conditions of
certainty and perfect competition; here $\bf x$ is the vector of
inputs, $F$ is a twice continuously differentiable production
function, $p$ is the output price, and $\bf w$ is the vector of
input prices. Clearly, since the scale factor $s$ is stipulated to
be positive, its magnitude does not affect the optimal values of the
decision variables under profit maximization, nor does it have any
effect on the comparative statics of the problem.  In any case, $s$
will be treated as a parameter whose value will eventually be set
equal to unity.  It should be evident here that the scale factor
serves an auxiliary purpose in this calculation and has no other
purpose, even though there are economically meaningful
interpretations of its role as will be discussed in \S IIB.

Let the parameter set be identified as $\mathbf{a}=
(a_{1},\ldots,a_{M+2})\stackrel{\rm{{def}}}{=}(w_{1},w_{2},\ldots,w_{M},p,s)
= (\mathbf{ w},p,s)$.  Then $\mathbf{\nabla}^\mathbf{a}=({\partial
\over
\partial w_{1}},\ldots, {\partial \over \partial w_{M}},{\partial
\over \partial p},{\partial \over \partial s})$, so that we find ${
\mathbf{\nabla}^\mathbf{a}}{\tilde
\phi}(\mathbf{x},\mathbf{a})=(-s\mathbf{x},sF,\phi)$, where $\phi
(\mathbf{x},\mathbf{a})\stackrel{\rm{{def}}}{=} pF(\bf x)-\mathbf{x}
\cdot\mathbf{w}$ is the standard function describing the firm's
profit. Note that there are now $M+2$ parameters under
consideration, so that $N=M+2$. We can therefore choose the $M+1$
isovectors according to $\mathbf{t}^{\alpha}
\stackrel{\rm{{def}}}{=} (0,0,\ldots,1,\ldots,0,0,{s{x}_{\alpha}
\over \phi})$, $\alpha=1,2,\ldots,M$, where the component equal to
unity occupies the $\alpha$th position within the parenthesis, and
$\mathbf{t}^{M+1} \stackrel{\rm{{def}}}{=} (0,\ldots,0,1,-s{F \over
\phi})$.  The corresponding set of compensated derivatives is easily
found to be $D_{\alpha}(\mathbf{x},\mathbf{a})={\partial \over
\partial w_{\alpha}} + ({s{x}_{\alpha} \over \phi}) {\partial \over
\partial s}$, $\alpha=1,2,\ldots,M$, and $D_{M+1}(\mathbf{x},\mathbf{a})={\partial \over
\partial p} -s{F \over \phi} {\partial \over
\partial s}$.  Note that
each $D_{\alpha}(\mathbf{x},\mathbf{a})$ has the null property with
respect to $\tilde \phi$.  Therefore these $M+1$ directional
derivatives constitute a set of GCD's for this problem and will be
used in \S IIIA to derive certain novel comparative statics results
for the standard profit maximization model.

It is appropriate at this juncture to emphasize an important feature
of GCD's briefly mentioned earlier.  While the customary meaning of
\textit{compensation} refers to a correction term that accounts for
the effect of a \textit{constraint} in the problem, the archetypical
example being the correction for the income effect in the utility
maximization problem, no necessary connection with constraints is
implied in the case of \textit{generalized compensation}, as is
clearly illustrated in the profit maximization case treated above.
Indeed any problem, constrained or not, will admit the use of
generalized compensation if its parameter space is larger than
one-dimensional. Furthermore, as the introduction of the scale
parameter $s$ for the problem treated above shows, the parameter
space can always be enlarged, so that the restriction to more than
one dimension is really no restriction at all.

The next problem to be considered here involves multiple
constraints.  It is therefore appropriate to first generalize the
notion of GCD's to the case of a set of functions
$\phi^{k}(\mathbf{x},\mathbf{a})$, $k=1,2,\ldots,K$, $K < N$, each
of which possesses the properties ascribed to
$\phi(\mathbf{x},\mathbf{a})$ in \S IIA.  To that end, consider a
fixed value of $\bf x$, and let $\bf a$ be a point in $\cal P$. Then
the normal hyperplane ${\cal N} (\mathbf{a})$ is defined to be the
vector space generated by the set of vectors
${\mathbf{n}}^{k}(\mathbf{a})
\stackrel{\rm{{def}}}{=}\mathbf{\nabla}^\mathbf{a}
\phi^{k}(\mathbf{x},\mathbf{a})$, $k=1,2,\ldots,K$.  Note that in
case the normal vectors constitute an \textit{independent} set,
i.e., if the only solution to $\sum_{k=1}^{K} c^{k}
{\mathbf{n}}^{k}(\mathbf{a})=0$, where every $c^{k}$ is a real
number, is $c^{k}=0$ for all $k$, then the normal hyperplane will be
$K$-dimensional. The subset of the tangent space ${\cal
T}^{N}(\mathbf{a})$ orthogonal to ${\cal N} (\mathbf{a})$ is then
the tangent hyperplane ${\cal I} (\mathbf{a})$ as before. Note that
when ${\cal N} (\mathbf{a})$ has dimension $K$, i.e., when
$\dim[{\cal N} (\mathbf{a})]=K$, then $\dim[{\cal I}
(\mathbf{a})]=N- K$. The case of interest for economic problems
usually corresponds to the condition $\dim[{\cal N} (\mathbf{a})]=K$
(known in decision space as a \textit{constraint qualification}
condition) with $K<N$. In any case, it should be emphasized here
that what is needed for the construction of GCD's is an independent
set of isovectors which can be characterized as a set of
$N-\dim[{\cal N} (\mathbf{a})]$ independent, $N$-dimensional vectors
each of which is orthogonal to ${\mathbf{n}}^{k}(\mathbf{a})$ for
every value of $k$.  Clearly, there is no need for an explicit
construction of ${\cal N} (\mathbf{a})$ here.

Having established the generalization to the multiconstraint case,
we now turn to the third problem illustrating the construction of
GCD's.  This is the principal-agent problem with hidden actions
(also known as \textit{moral hazard}) where a firm, the
\textit{principal}, intends to hire an individual, the
\textit{agent}, to work on a certain venture on a contractual basis
(see, e.g., Mas-Collel et al. 1995).  We shall return to this
problem in \S IIIB, where we define and formulate the problem, then
show that it can be transformed to the following convenient form:
\[
{\min \;}_\mathbf{x} {\;}{\sum}_{i=1}^{M}{x}_{i}{\cal P}_{i}^{I}
\;\;s.t. {\:} {\:} {B}^{k}-{\sum}_{i=1}^{M}v({x}_{i}){\cal
P}_{i}^{k} =0, \;\;{s}^{k}-{\sum}_{i=1}^{M}{\cal P}_{i}^{k}=0,
\;k=I,II,
\]
Here ${\cal P}_{i}^{k}$, where $0 < {\cal P}_{i}^{k} < 1$,
$i=1,\ldots, M$, $k=I,II$, is the probability that the $i$th profit
level, ${\pi}_{i}$, is realized for the firm given that the agent
performs at effort level $k$, with $I$ and $II$ corresponding to
\textit{high} and \textit{low effort} respectively.  The decision
variable ${x}_{i}$, on the other hand, is the agent's compensation
in case the $i$th profit level is realized. Furthermore,
${B}^{k}\stackrel{\rm{{def}}}{=}{c}^{k}+\bar{ u}$, where ${c}^{k}$
is the agent's disutility of working at effort level $k$, $\bar{ u}$
is the market price of the agent's services, and $v(x)-{c}^{k}$ is
the agent's utility function.  Finally, the parameters ${s}^{k}$ are
a pair of positive auxiliary variables which are introduced for
calculational convenience and will eventually be set equal to unity.
The parameter set is thus identified as $({\cal P}_{i}^{k},{B}^{k},
{s}^{k})$, with $i=1,2,\ldots ,M$ and $k=I,II$, for a total of
$2(M+2)$, to be eventually reduced to $2(M+1)$ upon setting
${s}^{k}=1$.

Rather than following the procedure used in the other illustrative
problems of this section, here we shall develop an intuitive
generalization of the GCD's already constructed for the Slutsky-Hicks
problem.  We accomplish this in two steps.  In the first step, we note
that for each value of $k$, the constraint ${B}^{k}-
{\sum}_{i=1}^{M}v({x}_{i}){\cal P}_{i}^{k} =0$ is analogous to a budget
constraint with ${B}^{k}$ standing for income, $v({x}_{i})$ for the
quantity of the $i$th good, and the probability ${\cal P}_{i}^{k}$ for the
price of that good.  This analogy immediately gives us the two-term
structure $\partial / \partial {\cal P}_{i}^{k} + v({x}_{i}) \partial /
\partial {B}^{k}$, $k=I,II$ which is clearly compensated (i.e., possess
the null property) with respect to the constraints involving
${B}^{k}$. The second step is to amend this structure so as to
extend the null property to the constraints involving ${s}^{k}$. One
can determine by analogy and verify by inspection that the addition
of $\partial / \partial {s}^{k}$ to the above structure provides the
desired extension of the null property to the constraints involving
${s}^{k}$ without disturbing the same property with respect to the
constraints involving ${B}^{k}$.  In short, the resulting structure,
${d}_{i}^{k}\stackrel{\rm{{def}}}{=}\partial /
\partial {\cal P}_{i}^{k} + v({x}_{i}) \partial / \partial {B}^{k} +
\partial / \partial {s}^{k}$, has the required null property with respect
to all constraints.  We can thus define
\[
D_{\alpha}(\mathbf{x},\mathbf{a}) \stackrel{\rm{{def}}}{=}
{d}_{\alpha}^{I}, \; \alpha =1,2,\ldots,M, \;
D_{\alpha}(\mathbf{x},\mathbf{a}) \stackrel{\rm{{def}}}{=}
{d}_{\alpha - M}^{II}, \; \alpha =M+1,M+2,\ldots,2M,
\]
as a set of GCD's for our problem.  We will use this set in \S IIIB to
derive comparative statics information for the principal-agent problem.

The fourth problem to be considered here is that of the
Pareto-optimal allocation of a given bundle of goods to a number of
individuals with given utility functions.  We will treat this
problem in full generality, primarily to illustrate the method for a
case with multiple, nontrivial constraints.  We shall assume here an
appropriate set of regularity conditions (e.g., twice continuously
differentiable, quasi-concave, strongly monotonic utility functions)
such that the desired optimality condition can be formulated as a
maximization problem subject to multiple constraints.  Returning to
the problem of allocating a resource bundle of $G$ goods,
${\omega}_{i} > 0$, $i=1,2,\ldots,G$, to $H$ individuals in a
Pareto-optimal manner, we consider maximizing
$u^{1}(\mathbf{x}^{1},\mathbf{b})$ by an appropriate choice of the
variables $\mathbf{x}^{h}\in {\Re}_{+}^{G}$, $h=1,2,\ldots,H$, or
equivalently $\mathbf{x} \in {\Re}_{+}^{M}$, $M=G \times H$, subject
to a set of $K=H-1+G$ constraints. These constraints are
$\phi^{k}(\mathbf{x},\mathbf{a})\stackrel{\rm{{def}}}{=}\bar{{u}}^{k}-
u^{k}(\mathbf{x}^{k},\mathbf{b})=0$ for $k=2,3,\ldots,H$, and
$\phi^{k}(\mathbf{x},\mathbf{a})\stackrel{\rm{{def}}}{=}{\omega}_{k-H}
- {\sum}_{h=1}^{H}{x}_{k-H}^ {h} =0$ for $k=H+1,H+2,\ldots,H+G$.
Here $u^{h}(\mathbf{x}^{h},\mathbf{b})$ is the utility function of
the $h$th agent, $\bf b$ is a set of $L$ parameters, and
$\bar{{u}}^{h}$ is the fixed utility level of the $h$th agent.
Moreover, the set of $N=L+H-1+G$ parameters comprising the vector
$\bf a$ are identified as
$(b_{1},\ldots,b_{L},\bar{{u}}^{2},\ldots,\bar{{u}}^{H},{\omega}_{1},
\ldots, {\omega}_{G})$.  We shall assume as given a continuously
differentiable, interior solution to this maximization problem, and
concern ourselves with the construction of the GCD's for the
problem.

As a first step in the construction of the GCD's, we consider the
$H-1+G$ normal vectors with components $\partial
\phi^{k}(\mathbf{x},\mathbf{a}) / \partial {a}_{\alpha}$,
$\alpha=1,2, \ldots,L+H-1+G$.  For $k=2,3,\ldots,H$, these
components are given by $- \partial u^{k}(\mathbf{x}^{k},\mathbf{b})
/
\partial {b}_{\alpha} $ for $\alpha=1,2,\ldots,L$, by
${\delta}_{\alpha , L+k-1}$ for $\alpha=L+1, L+2,\ldots,L+H-1$, and
they vanish for $\alpha=L+H,L+H+1,\ldots,L+H-1+G$. For
$k=H+1,H+2,\ldots,H+G$, the components of the normal vectors vanish
for $\alpha =1,2, \ldots, L+H-1$, and are equal to
${\delta}_{k-1,\alpha-L}$ for $\alpha=L+H, L+H+1,\ldots,L+H-1+G$. We
can compactly rewrite these expressions as
$\mathbf{\nabla}^\mathbf{a} \phi^{k}(\mathbf{x},\mathbf{a}) =[ -
\mathbf{\nabla}^\mathbf{b}u^{k}(\mathbf{x}^{k},\mathbf{b}),
\mathbf{e}_{H-1}^{k-1},\mathbf{0} ]$ for $k=2,3,\ldots,H$, and
$\mathbf{\nabla}^\mathbf{a} \phi^{k}(\mathbf{x},\mathbf{a})
=[\mathbf{0},\mathbf{0}, \mathbf{e}_{G}^{k-H}]$ for
$k=H+1,H+2,\dots,H+G$ in a shorthand notation where
$\mathbf{e}_{J}^{j}$ denotes a $J$-dimensional unit vector pointing
in the $j$th direction.  The next step is the construction of a set
of isovectors $\mathbf{t}^{\alpha}$.  There will be $L$ of these if
the normal vectors are independent, which is the typical situation
and will be assumed to be the case here.  As in the case of the
previous models, we can construct these by inspection.  For example,
we can choose, for $\alpha=1,\ldots,L$, $\mathbf{t}^{\alpha}
\stackrel{\rm{{def}}}{=} [\mathbf{e}_{L}^{\alpha}, {\partial \over
\partial {b}_{\alpha}} u^{2}(\mathbf{x}^{2},\mathbf{b}),{\partial \over
\partial {b}_{\alpha}} u^{3}(\mathbf{x}^{3},\mathbf{b}),\ldots, {\partial
\over \partial {b}_{\alpha}}
u^{H}(\mathbf{x}^{H},\mathbf{b}),\mathbf{0} ]$.

The set of $L$ compensated derivatives corresponding to the
isovector set just constructed can be expressed as
$D_{\alpha}(\mathbf{x},\mathbf{a})=\mathbf{t}^{\alpha} \cdot
\mathbf{\nabla}^\mathbf{a}={\partial \over \partial b_{\alpha}} +
{\sum}_{h=2}^{H} [{\partial \over \partial {b}_{\alpha}}
u^{h}(\mathbf{x}^{h} ,\mathbf{b})] {\partial \over \partial
\bar{{u}}^{h}}$, $\alpha=1,2,\ldots,L$.  As usually formulated, the
general Pareto-optimal allocation problem yields multiple solutions
(the Pareto set). Consequently, additional conditions must be
imposed in order to render the solution unique and comparative
statics questions meaningful.  Once this is done, the GCD set
constructed above can be used to derive constraint-free comparative
statics results for the Pareto-optimal allocation problem according
to Eq.~(\ref{7}) of \S IIB in a straightforward manner.  Note that
this GCD set does not involve the parameters ${\omega}_{i}$, and
that we would have arrived at the same GCD's had we elected not to
include the $\omega$'s in the parameter set in the first place. This
is so since in any event the partial derivatives with respect to the
remaining parameters $\bf b$ and $\bar{{u}}^{k}$ annihilate the
resource constraint equations, i.e., they are already
``compensated'' with respect to the latter.  This will always occur
if (a) certain parameters are the only ones appearing in certain
constraint equations, and (b) those parameters do not appear in any
other constraint equations.  Note also that the GCD's just
constructed were not compensated with respect to the objective
function, i.e., $D_{\alpha}(\mathbf{x},\mathbf{a})
u^{1}(\mathbf{x}^{1},\mathbf{b})$ does not vanish by construction.
As we will show below, this property is not needed for obtaining
constraint-free comparative statics results.  The fact that this
property was incorporated into the profit maximization problem
treated above, on the other hand, was simply motivated by the
calculational power that generalized compensation makes available
for deriving various comparative statics results, as will be
demonstrated in \S IIC and in the applications.

Having demonstrated the construction of GCD's for four basic
economic models, we now resume the main development and proceed to
establish the most important property of a GCD, namely its
constraint conformance property mentioned earlier.  To that end, we
start by reversing the roles of $\bf x$ and $\bf a$, consider the
latter as fixed, and proceed to apply the construction developed
above in parameter space to decision space.  In particular, an
isovector $\bf s$ in decision space is an $M$-dimensional vector
which is orthogonal to all the normal directions associated with the
set $\phi^{k}(\mathbf{x},\mathbf{a})$, i.e., $\mathbf{s}\cdot
{\nabla}^\mathbf{x}{\phi}^{k} (\mathbf{x},\mathbf{a})=0$,
$k=1,\ldots,K$. If, moreover, some of the functions in the set, say
those corresponding to $k=1,\ldots,C$, with $C \leq K$, serve as
constraint functions by virtue of the requirement that they must
equal zero, then the corresponding isovectors in decision space are
precisely those that \textit{conform} to the constraints by virtue
of pointing in directions that are tangent to the level surfaces
defined by the constraints in decision space.  The isovectors in
decision space are therefore \textit{conforming} vectors in this
sense, and for that reason play a crucial role in the construction
of the matrix that conveys constraint-free comparative statics
results.

To make these assertions more transparent, let us now consider a
restriction of the above constructions in the decision and parameter
spaces to the case where the two vector arguments of the functions
$\phi^{k}(\mathbf{x},\mathbf{a})$ are functionally related.
Specifically, let $\mathbf{x}$ be a continuously differentiable
function of $ \bf a$, i.e., $\mathbf{x} = \mathbf{x}(\bf a)$, and
consider the restricted set of functions
$\phi^{k}(\mathbf{x}(\mathbf{a}),\mathbf{a})$.  In applications, the
vector-valued function $\mathbf{x}(\cdot)$, which we shall refer to
as a \textit{decision function}, is derived from some optimality
condition, subject to constraints if present. Such a constraint
condition therefore implies that
$\phi^{k}(\mathbf{x}(\mathbf{a}),\mathbf{a})=0$ for $k=1,2,\ldots,C$
and $\bf a$ in some subset of $\cal P$. In other words, when
restricted to $\mathbf{x} = \mathbf{x}(\bf a)$, each function in the
constraint set identically vanishes in the variable $\bf a$.  To
exploit this property, let us apply the now restricted parameter
space GCD's $D_{\alpha}(\mathbf{x}(\mathbf{a}),\mathbf{a})$,
henceforth abbreviated by $D_{\alpha}(\mathbf{a})$, to the
constraint identities. The result is that for all $\mathbf{a} \in
\cal P$ and $k=1,\ldots,C$,
\[0=D_{\alpha}(\mathbf{a})\phi^{k}(\mathbf{x}(\mathbf{a}),
\mathbf{a})={\sum}_{i=1}^{M}
{\phi}_{,i}^{k}(\mathbf{x}(\mathbf{a}),\mathbf{a}){x}_{i;\alpha}(\mathbf{a})
+{\phi}_{;\alpha}^{k}(\mathbf{x}(\mathbf{a}),\mathbf{a}), \] where
${\phi}_{,i}^{k}(\mathbf{x},\mathbf{a})\stackrel{\rm{{def}}}{=}{\partial
\over
\partial {x}_{i}}
{\phi}^{k}(\mathbf{x},\mathbf{a})$,
${x}_{i;\alpha}(\mathbf{a})\stackrel{\rm{{def}}}{=}D_{\alpha}
(\mathbf{a}){x}_{i}(\mathbf{a})=\mathbf{t}^{\alpha}
\cdot { \mathbf{\nabla}^\mathbf{a}}{x}_{i}(\mathbf{a})$, and
${\phi}_{;\alpha}^{k}(\mathbf{x},\mathbf{a})\stackrel{\rm{{def}}}{=}
D_{\alpha}(\mathbf{a}){\phi}^{k}(\mathbf{x},\mathbf{a})$.
However, owing to the null property of the GCD's, the second term on
the right hand side of the above equations vanishes, leaving behind
the term involving the compensated derivatives of the decision
functions.   Note the appearance of decision space normal vectors
${\phi}_{,i}^{k}(\mathbf{x},\mathbf{a})$ in the surviving term of
the above equation.  Also note the introduction of a notational
convention here whereby a subscript occurring to the right of a
comma signifies partial differentiation, whereas a subscript
occurring to the right of a semicolon signifies directional
differentiation corresponding to a GCD.  Moreover, Latin subscripts
are used to denote differentiation with respect to decision
variables, while Greek indices are used for differentiation in
parameter space.

The surviving term in the last-stated equation above is in fact the
inner product of ${x}_{i;\alpha}(\bf a)$ with decision space normal
vectors ${\phi}_{,i}^{k}(\mathbf{x},\mathbf{a})$, so that its
vanishing for every value of $k$ and $\alpha$ implies the
orthogonality of every generalized compensated derivative of
$\mathbf{x}(\mathbf{a})$ to every normal vector associated with the
constraint surfaces in decision space.  Equivalently, the stated
orthogonality implies the tangency of ${x}_{i;\alpha}(\bf a)$ to the
decision space surface ${\phi}^{k}(\mathbf{x},\mathbf{a})=0$ (with
$\bf a$ fixed) for every value of $\alpha$ and $k$.  We have thus
established the fact that the application of parameter space GCD's
to the decision functions produces isovectors in decision space,
i.e., vectors that conform to the constraints.  It is useful to
examine more closely how this crucial property of GCD's emerges.
Because the restricted constraint functions
${\phi}^{k}(\mathbf{x}(\bf a),\mathbf{a})$ vanish identically in
$\bf a$, and since a ``small'' change $\delta \bf
a\stackrel{\rm{{def}}}{=}\epsilon {\sf t}^{\alpha}$ which points in
a null direction $\mathbf{t}^{\alpha}$ induces no variation in
${\phi}^{k}(\mathbf{x}(\bf a),\mathbf{a})$ arising from its
dependence on the second argument, there cannot be any variation
arising from the change in the first argument, $\epsilon
\mathbf{t}^{\alpha} \cdot {\nabla}^\mathbf{a}
\mathbf{x}(\mathbf{a})$, either.  But the last statement
characterizes $\mathbf{t}^{\alpha} \cdot {\nabla}^\mathbf{a}
\mathbf{x}(\mathbf{a})$ as an isovector in decision space. Recalling
that ${\sf t}^{\alpha} \cdot
{\nabla}^\mathbf{a}\stackrel{\rm{{def}}}{=}D_{\alpha}(\mathbf{a})$,
we arrive at the desired result.  This is the most important
property of a GCD and, under the conditions stipulated above, may be
summarized as

\textbf{Lemma 1.}  \textit{Every generalized compensated derivative
of a decision function conforms to the constraints in decision
space.}

We pause here to emphasize the profoundly dual nature of this
result, as is already evident in the arguments preceding the
statement of the lemma: a directional derivative which annihilates
the constraint functions in \textit{parameter space}, when applied
to decision functions, will produce a vector which conforms to the
constraints in \textit{decision space}.

\subsection{The Main Theorems}

Having described the construction and properties of GCD's in the previous
section in some detail, we are now in a position to establish the main
theorem of this paper.  Consider the optimization problem
\begin{equation}
{\max \;}_\mathbf{x} {\;}f(\mathbf{x},\mathbf{a}) {\:} {\:} s.t.
{\:} {\:} {g}^{k}(\mathbf{x},\mathbf{a})=0,\label{1}
\end{equation}
where $f , {g}^{k}$ are twice continuously differentiable with
$k=1,2,\ldots,K$ and $M , N > K$.  As stipulated above, $\mathbf{x}
$ and $\mathbf{a}$ are points in decision and parameter spaces,
respectively. To avoid trivialities, we shall furthermore require
the set of constraint functions to be independent at the optimum
point.  This requirement implies that the set of parameter space
normal vectors is linearly independent and the parameter space
normal hyperplane is of dimension $K$ at the optimum point.  It is a
simple matter to show that this requirement also implies a parallel
condition in decision space, the well known constraint qualification
condition.

The task before us now is the establishment of the comparative
statics corresponding to a given solution of the above problem.  To
that end, let us suppose there exists a unique, continuously
differentiable, interior solution to the above problem specified by
the (vector-valued) decision function $\mathbf{x}(\mathbf{a})$.  In
other words, we suppose the existence of an open set ${\cal P}$ in
parameter space and a set ${\cal D}$ in decision space together with
a continuously differentiable function $\mathbf{x}(\cdot)$ from
${\cal P}$ to ${\cal D}$ such that for each value of $\mathbf{a} \in
{\cal P}$, the point $\mathbf{x}(\mathbf{a})$ is in the interior of
$\cal D$ and possesses the constrained maximum property stated in
Eq.~(\ref{1}). Note that we have not included any inequality
constraints, such as those arising from nonnegativity conditions
imposed on decision variables, in Eq.~(\ref{1}) since the interior
nature of the solution in effect obviates any such conditions.

In order to characterize the constrained maximum property of
$\mathbf{x}(\mathbf{a})$, we construct a set of GCD's with respect
to the constraint functions given in Eq.~(\ref{1}), possibly
including the objective function as well,  according to the
procedure explained in \S IIA.  It is also useful to introduce the
notation
$\mathbf{h}(\mathbf{a})\stackrel{\rm{{def}}}{=}{\sum}_{\alpha=1}^{A}{\eta}_{\alpha}
\mathbf{x}_ {;\alpha}(\mathbf{a})$, where ${\eta}_{\alpha}$ is an
arbitrary vector of real numbers, $A$ is the number of GCD's, and
$\mathbf{x}_{;\alpha}(\mathbf{a})=\mathbf{t}^{\alpha} \cdot
{\nabla}^\mathbf{a} \mathbf{x}(\mathbf{a})=D_{\alpha}(\mathbf{a})
\mathbf{x}(\mathbf{a})$ using the notation established in \S IIA.
According to Lemma 1 established above,
$\mathbf{x}_{;\alpha}(\mathbf{a})$ conforms to the constraints in
decision space for every value of $\alpha$.  But then the same is
implied for $\mathbf{h}(\mathbf{a})$, since it is a linear
combination of conforming vectors (recall that conforming vectors
are elements of a linear space, the tangent hyperplane in decision
space).  By construction, then, the vector $\mathbf{h}(\mathbf{a})$
conforms to the constraints.

Using the above construction, we can use the maximum condition stated in
(\ref{1}) to assert that for any $\epsilon$ of sufficiently small
magnitude, we have the condition
\begin{equation}
f(\mathbf{x}(\mathbf{a})+\epsilon \mathbf{h}(\mathbf{a}),\mathbf{a}
) \leq f(\mathbf{x}(\mathbf{a}),\mathbf{a}),\label{2}
\end{equation}
a statement that holds \textit{free of constraints}.  But then the
arbitrary sign of $\epsilon$ in Eq.~(\ref{2}) can be used in
conjunction with the differentiability property of $f$ to deduce
that the directional derivative of $f$ in the direction of
$\mathbf{h}(\mathbf{a})$ must be nonpositive as well as nonnegative,
hence the result that it must vanish;
\begin{equation}
{\sum}_{i=1}^{M}{h}_{i}(\mathbf{a}){f}_{,i}(\mathbf{x}(\mathbf{a}),\mathbf{a})
=0,\label{3}
\end{equation}
Eq.~(\ref{3}) is the constraint-free, necessary, first-order
condition implied by the constrained maximum property stated in
Eq.~(\ref{1}).

We next explore the consequences of the first-order condition just
derived.  This condition simply characterizes
${f}_{,i}(\mathbf{x}(\mathbf{a}),\mathbf{a})$ as being orthogonal to
the tangent hyperplane in decision space, hence as belonging to the
normal hyperplane in decision space.  But the last statement implies
that, for each value of $\bf a$, the gradient vector
${f}_{,i}(\mathbf{x}(\mathbf{a}),\mathbf{a})$ may be expressed as a
linear combination of the normal vectors to the constraint surfaces
in decision space. Since the latter are simply given by the
gradients of the constraint functions in decision space, we find
\begin{equation}
{f}_{,i}(\mathbf{x}(\mathbf{a}),\mathbf{a})
\stackrel{\rm{{def}}}{=}-{\sum}_{k=1}^{K} {\lambda}_{k}(\mathbf{a})
{g}_{,i}^{k}(\mathbf{x}(\mathbf{a}),\mathbf{a}). \label{4}
\end{equation}
The scalars ${\lambda}_{k}$ that enter the linear combination on the
right hand side of Eq.~(\ref{4}) are of course the familiar
multipliers of Lagrange's method of constrained maximization (with
the negative sign inserted for convenience).

Thus far we have essentially established the validity of Lagrange's
method for dealing with constrained optimization problems using our
construction.  Having established that connection, we now recall the
second-order, necessary conditions associated with the maximization
problem of Eq.~(\ref{1}) (see, e.g., Mas-Colell et al. 1995;
Takayama 1985).  That condition states that for any conforming
vector $\mathbf{l}$, i.e., any vector $\mathbf{l}$ which satisfies
the condition ${\sum}_{i=1}^{M}
{l}_{i}{g}_{,i}^{k}(\mathbf{x}(\mathbf{a}),\mathbf{a}) =0$ for every
value of $k$, we have the inequality
\begin{equation}
{\sum}_{i,j=1}^{M}{l}_{i}{l}_{j}[{f}_{,ij}(\mathbf{x}(\mathbf{a}),\mathbf{a})
+{\sum}_{k=1}^{K}
{\lambda}_{k}(\mathbf{a}){g}_{,ij}^{k}(\mathbf{x}(\mathbf{a}),\mathbf{a})]
\leq 0. \label{5}
\end{equation}
Since in the present instance the vector $\mathbf{h}(\mathbf{a})$
conforms to the constraints by construction, Eq.~(\ref{5}) would
hold \textit{free of constraints} if $\mathbf{h}(\mathbf{a})$ is
substituted for $\mathbf{l}$. Recalling that
$\mathbf{h}(\mathbf{a})$ is equal to ${\sum}_{\alpha=1}^{A}
{\eta}_{\alpha} \mathbf{x}_ {;\alpha}(\mathbf{a})$, and remembering
that the vector of real numbers ${\eta}_{\alpha}$ appearing therein
is arbitrary, we recognize in Eq.~(\ref{5}) the statement that the
matrix
\[{\sf{\Omega}}_{\alpha \beta}(\mathbf{a})
\stackrel{\rm{{def}}}{=}-
{\sum}_{i,j=1}^{M}{x}_{i;\alpha}(\mathbf{a}){x}_{j
;\beta}(\mathbf{a})[{f}_{,ij}(\mathbf{x}(\mathbf{a}),\mathbf{a})
+{\sum}_{k=1}^{K}
{\lambda}_{k}(\mathbf{a}){g}_{,ij}^{k}(\mathbf{x}(\mathbf{a}),\mathbf{a})]
\] is positive semidefinite.  The symmetry property of
${\sf{\Omega}}$ is a consequence of the symmetry of
${f}_{,ij}(\mathbf{x}(\mathbf{a}),\mathbf{a})$ and
${g}_{,ij}^{k}(\mathbf{x}(\mathbf{a}),\mathbf{a})$, the latter
following from the symmetry of the matrix of second-order partial
derivatives for twice continuously differentiable functions. Let us
recall here that a real matrix $\sf A$ is by definition
\textit{positive definite} or \textit{semidefinite} if (a) it is
symmetric and (b) for every real vector $\bf v \neq 0$, the
quadratic form $\mathbf{v}^{\dag} {\sf A} \mathbf{v}$ is
\textit{positive definite} or \textit{semidefinite} respectively.

The above expression for ${\sf{\Omega}}$ can be rewritten by first
applying a GCD, say ${D}_{\alpha}(\mathbf{a})$, to Eq.~(\ref{4}),
obtaining
\begin{eqnarray}
&&{\sum}_{j=1}^{M}[{f}_{,ij}(\mathbf{x}(\mathbf{a}),\mathbf{a})
+{\sum}_{k=1}^{K}
{\lambda}_{k}(\mathbf{a}){g}_{,ij}^{k}(\mathbf{x}(\mathbf{a}),\mathbf{a})]{x}_{j
;\alpha}(\mathbf{a})=\nonumber\\
&&-[{f}_{,i; \alpha}(\mathbf{x}(\mathbf{a}),\mathbf{a})
+{\sum}_{k=1}^{K} {\lambda}_{k}(\mathbf{a}){g}_{,i;
\alpha}^{k}(\mathbf{x}(\mathbf{a}),\mathbf{a})] -{\sum}_{k=1}^{K}
{\lambda}_{k;
\alpha}(\mathbf{a}){g}_{,i}^{k}(\mathbf{x}(\mathbf{a}),\mathbf{a}).
\label{6}
\end{eqnarray}
If Eq.~(\ref{6}) is now multiplied by ${x}_{i; \beta}$ and summed
over $i$, the second term on the right hand side vanishes by the
conformance property of ${x}_{i; \beta}$ that follows from Lemma 1
of \S IIA, leaving a simplified expression for $\sf{\Omega}$:
\begin{equation}
{\sf{\Omega}}_{\alpha \beta}(\mathbf{a}) = {\sum}_{i=1}^{M}{x}_{i;
\beta}(\mathbf{a})[{f}_{,i;
\alpha}(\mathbf{x}(\mathbf{a}),\mathbf{a}) +{\sum}_{k=1}^{K}
{\lambda}_{k}(\mathbf{a}){g}_{,i;
\alpha}^{k}(\mathbf{x}(\mathbf{a}),\mathbf{a})]. \label{7}
\end{equation}
It is appropriate to recall here that the semidefiniteness property
of $\sf \Omega$ implies that it is symmetric, i.e.,
${\sf{\Omega}}_{\alpha \beta}(\mathbf{a})={\sf{\Omega}}_{\beta
\alpha}(\mathbf{a})$, a fact that we have already used in writing
Eq.~(\ref{7}) and will continue to utilize throughout this work. We
have thus arrived at the result that the matrix $\sf{\Omega}$, which
is a linear combination of the partial derivatives of the decision
functions with respect to the parameters, is positive semidefinite,
free of constraints.  We shall refer to a matrix possessing these
properties as a \textit{comparative statics matrix} (or CSM in
abbreviated form) for the optimization problem.  The unrestricted
existence of a comparative statics matrix for a general, constrained
optimization problem is the central result of our analysis:

\textbf{Theorem 1.} \textit{The constrained optimization problem
defined by Eq.~(\ref{1}) et seq. admits of a constraint-free
comparative statics matrix $\sf{\Omega}$ given in Eq.~(\ref{7}).}

It is worth reemphasizing here that there is a large freedom of
choice in the construction of CSM's, a feature that will be explored
in the following and summarized in Theorem 2.  While this freedom
may be exploited to generate different forms of comparative statics
for a given optimization problem, it is well to remember that all
such matrices convey no more information than is contained in the
second-order, necessary conditions expressed in Eq.~(\ref{5}). These
conditions in turn originate in the local concavity of the
underlying constrained maximization problem defined in
Eq.~(\ref{1}).  A more intuitive discussion of these and related
matters will be given at the end of \S IIC.

Having established the main result of this paper, we now proceed to
discuss how such features as the envelope and homogeneity properties are
realized in the present framework.  We will then consider the
nonuniqueness features of our construction and develop a characterization
of the associated arbitrariness in the resulting CSM's.

As one might surmise, the null property of our GCD's ensures that
the envelope property holds for the general constrained optimization
problem defined by Eq.~(\ref{1}) as a property of the optimized
objective function and without the intrusion of constraint
functions.  To see this explicitly, let us start by introducing the
\textit{value function}
$V(\mathbf{a})\stackrel{\rm{{def}}}{=}{f}(\mathbf{x}(\mathbf{a}),\mathbf{a})$,
and observe that it is also equal to
${f}(\mathbf{x}(\mathbf{a}),\mathbf{a}) +{\sum}_{k=1}^{K}
{\lambda}_{k}(\mathbf{a})
{g}^{k}(\mathbf{x}(\mathbf{a}),\mathbf{a})$ by the identical
vanishing of the restricted constraint functions for all $\mathbf{a}
\in \cal P $. Next consider the result of applying a GCD, say
$D_{\alpha}(\mathbf{a})$, to the value function as given by the
second expression above. Then, using Eq.~(\ref{4}) and the null
property ${g}_{;\alpha}^{k}(\mathbf{x}(\mathbf{a}),\mathbf{a})=0$,
we find the desired result that
$D_{\alpha}(\mathbf{a})V(\mathbf{a})={V}_{;\alpha}(\mathbf{a})={f}_{;
\alpha}(\mathbf{x}(\mathbf{a}), \mathbf{a})$, which is a statement
of the envelope property in terms of GCD's.  Note that in case the
GCD's have the (optional) null property with respect to the
objective function, i.e., if ${f}_{; \alpha}(\mathbf{x}(\mathbf{a}),
\mathbf{a})=0$, then the envelope property reduces to the statement
that ${V}_{;\alpha}(\mathbf{a})=0$, i.e., that the value function
also possesses the null property.  In other words, if the the null
property of the GCD's is extended to to the objective function, then
the value function has vanishing derivatives in the null directions
$\mathbf{t}^{\alpha}$.  Note that this condition is equivalent to a
set of $A$ first-order partial differential equations satisfied by
the value function.  Needless to say, each of these is an invariance
property of the value function.

Next we explore the consequences of any useful symmetry, or
invariance property that the objective and constraint functions
might possess. Consider the optimization problem defined by
Eq.~(\ref{1}), together with the assumed solution
$\mathbf{x}(\mathbf{a})$, with $\mathbf{a} \in \cal P $.  The task
before us is to find out how a given invariance possessed by the
objective and constraint functions of this problem, such as
homogeneity of a given degree, gives rise to a corresponding
property of the decision functions.  To that end, let us start by
characterizing the class of symmetries we would like to consider,
limiting our attention to those invariances that are likely to play
a significant role in economic applications.

Suppose there exist a pair of real, continuously differentiable
vector functions $\mathbf{X}(\mathbf{x})$ and
$\mathbf{A}(\mathbf{a})$, from $\cal D $ and $\cal P $ to
${\Re}^{M}$ and ${\Re}^{N}$ respectively, and consider the
differential operator $J$ defined by
\[ J(\mathbf{x},\mathbf{a})\stackrel{\rm{{def}}}{=}{\sum}_{i=1}^{M} {X}_{i}(\mathbf{x}){\partial \over
{\partial {x}_{i}}}+{\sum}_{\mu=1}^{N}
{A}_{\mu}(\mathbf{a}){\partial \over {\partial {a}_{\mu}}}. \]
Suppose further that the action of $J$ on the objective and
constraint functions can be described by suitable functions, that
is,
\[ J(\mathbf{x},\mathbf{a})f(\mathbf{x},\mathbf{a})=\rm{ F}(f(\mathbf{x},\mathbf{a})), \]
\begin{equation}
J(\mathbf{x},\mathbf{a}){g}^{k}(\mathbf{x},\mathbf{a})=\rm{
G}^{k}({g}^{k} (\mathbf{x},\mathbf{a})), k=1,\ldots,K, \label{8}
\end{equation}
for every $\mathbf{x} \in \cal D $ and $\mathbf{a} \in \cal P $.
Here $\rm F$ and $\rm{ G}^{k}$ are continuously differentiable
functions from ${\Re}^{1}$ to ${\Re}^{1}$, with $\rm{ G}^{k}(0)=0$,
$k=1,\ldots,K$.

These conditions, obscure as they may appear, actually have a
straightforward interpretation as invariance conditions.  They
essentially state that if the objective and constraint functions are
evaluated at the ``slightly'' displaced values of their arguments
$\mathbf{x}+\epsilon J(\mathbf{x},\mathbf{a})
\mathbf{x}=\mathbf{x}+\epsilon \mathbf{X}(\mathbf{x})$ and
$\mathbf{a}+\epsilon J(\mathbf{x},\mathbf{a})
\mathbf{a}=\mathbf{a}+\epsilon \mathbf{A}(\mathbf{a})$, instead of
$\bf x$ and $\bf a$ respectively, where $\epsilon$ is a ``small''
real number, then the underlying optimization problem remains
unchanged to first order in $\epsilon$. Given a suitable uniqueness
requirement on the solution function $\mathbf{x}(\mathbf{a})$, this
first-order invariance condition would imply that the modified
objective and constraint functions define a solution that differs
from the solution of the original problem only by quantities of
second order in $\epsilon$. In other words, they imply that the
decision vector $\mathbf{x}(\mathbf{a})+\epsilon
\mathbf{X}(\mathbf{x})$ differs from the solution
$\mathbf{x}(\mathbf{a}+\epsilon \mathbf{A}(\mathbf{a}))$ by
second-order quantities only (Euclidean norm implied).  The last
statement can then be converted into the following invariance
property for $\mathbf{x}(\mathbf{a})$ in the limit of vanishing
$\epsilon$:
${X}_{i}(\mathbf{x}(\mathbf{a}))-{\sum}_{\mu=1}^{N}{A}_{\mu}(\mathbf{a}){x}_{i,\mu}(\mathbf{a})=0$.

The heuristic argument given in the foregoing paragraph can be formalized
in a straightforward manner.  To do so, we start with the assumed
invariance conditions obeyed by the objective and constraint functions,
Eqs.\ (\ref{8}), and proceed to differentiate these with respect to the
decision variables.  The result for the objective function is
\[ {\sum}_{i=1}^{M} [{X}_{i,j}(\mathbf{x}){f}_{,i} (\mathbf{x},\mathbf{a})+{X}_{i}(\mathbf{x}){f}_{,ij} (\mathbf{x}, \mathbf{a})] +{\sum}_{\mu=1}^{N}
{A}_{\mu}(\mathbf{a}) {f}_{,\mu j}(\mathbf{x},\mathbf{a})=\rm{ F}'
(f (\mathbf{x}, \mathbf{a}) ){f}_{,j} (\mathbf{x}, \mathbf{a}), \]
where a prime signifies differentiation with respect to the
argument. A similar equation results for each constraint function
${g}^{k}$. Next we multiply the resulting equation for ${g}^{k}$ by
${\lambda}_{k}$, sum this over $k$, and add the result to the above
equation for $f$ in order to obtain an equation involving the
partial derivatives of the combination
$f(\mathbf{x})+{\sum}_{k=1}^{K} {\lambda}_{k}
 {g}^{k}(\mathbf{x}, \mathbf{a})\stackrel{\rm{{def}}}{=}
{L}(\mathbf{x}, \mathbf{a})$, the Lagrange function associated with
the optimization problem.  Here, as elsewhere in this paper, the
dependence of $L$ on the multipliers has been suppressed to avoid
clutter in the notation.  Since only partial derivatives of $L$ with
respect to $\bf x$ and $\bf a$ will appear in our notation, the
suppression of $\lambda$ should not cause any confusion. Restricting
the resulting equation for $L$ to the solution of the optimization
problem by substituting $\mathbf{x}(\mathbf{a})$ for $\bf x$ and
using the first-order condition given in Eq.~(\ref{4}), we find
\begin{eqnarray}
&&{\sum}_{i=1}^{M} {X}_{i}(\mathbf{x}(\mathbf{a})){L}_{,ij}
(\mathbf{x}(\mathbf{a}), \mathbf{a}) +{\sum}_{\mu=1}^{N}
{A}_{\mu}(\mathbf{a}) {L}_{,\mu j}(\mathbf{x}(\mathbf{a}),\mathbf{a})= \nonumber \\
&&-{\sum}_{k=1}^{K} {\lambda}_{k}
{g}_{,j}^{k}(\mathbf{x}(\mathbf{a}), \mathbf{a})[\rm{ F}' (f
(\mathbf{x}(\mathbf{a}), \mathbf{a}) )- \rm{ G}^{k
\prime}({g}^{k}(\mathbf{x}(\mathbf{a}), \mathbf{a})] \nonumber
\end{eqnarray}
Similarly, restricting the assumed invariance condition for
${g}^{k}$ to the solution $\mathbf{x}=\mathbf{x}(\mathbf{a})$, we
find ${\sum}_{i=1}^{M} {X}_{i}
(\mathbf{x}(\mathbf{a})){g}_{,i}^{k}(\mathbf{x}(\mathbf{a}),
\mathbf{a})+{\sum}_{\mu=1}^{N} {A}_{\mu}(\mathbf{a})
{g}_{,\mu}^{k}(\mathbf{x}(\mathbf{a}), \mathbf{a})=0$.  On the other
hand, a differentiation of the $k$th constraint equation with
respect to ${a}_{\mu}$ leads to
${g}_{,\mu}^{k}(\mathbf{x}(\mathbf{a}), \mathbf{a})
+{\sum}_{i=1}^{M}{g}_{,i}^{k}({x}(\mathbf{a}), \mathbf{a})
{x}_{i,\mu}(\mathbf{a})=0$.  Eliminating
${g}_{,\mu}^{k}(\mathbf{x}(\mathbf{a}), \mathbf{a})$ from the last
two equations, we get
\[ {\sum}_{i=1}^{M}{Z}_{i}(\mathbf{a}){g}_{,i}^{k}(\mathbf{x}(\mathbf{a}), \mathbf{a})=0, k=1,\ldots,K, \]
where $\mathbf{Z}(\mathbf{a})\stackrel{\rm{{def}}}{=}
\mathbf{X}(\mathbf{x}(\mathbf{a}))-{\sum}_{\mu=1}^{N}{A}_{\mu}(\mathbf{a})\mathbf{x}_{,\mu}(\mathbf{a})$.
Note that the relation just derived is an orthogonality condition in
decision space which identifies $\mathbf{Z}$ as a conforming vector.
If the equation derived above for $L$ is multiplied by ${Z}_{j}$ and
summed over $j$, we find, taking account of the conformance property
of $\mathbf{Z}$, the result
\[ {\sum}_{i,j=1}^{M} {X}_{i}(\mathbf{x}(\mathbf{a})){L}_{,ij} (\mathbf{x}(\mathbf{a}), \mathbf{a}){Z}_{j}(\mathbf{a}) +{\sum}_{j=1}^{M} {\sum}_{\mu=1}^{N}
{A}_{\mu}(\mathbf{a}) {L}_{,\mu
j}(\mathbf{x}(\mathbf{a}),\mathbf{a}) {Z}_{j}(\mathbf{a}) =0. \]

To
develop this equation further, we differentiate the first-order
condition given in Eq.~(\ref{4}) with respect to ${a}_{\mu}$ and use
the conformance property of $\mathbf{Z}$ to rewrite the term
${\sum}_{j=1}^{M} {\sum}_{\mu=1}^{N} {A}_{\mu}(\mathbf{a}) {L}_{,\mu
j}(\mathbf{x}(\mathbf{a}),\mathbf{a}) {Z}_{j}(\mathbf{a})$ as
$-{\sum}_{i,j=1}^{M} {\sum}_{\mu=1}^{N} {A}_{\mu}(\mathbf{a})
{L}_{,ij}(\mathbf{x}(\mathbf{a}),\mathbf{a})
{x}_{i,\mu}(\mathbf{a}){Z}_{j}(\mathbf{a})$.  Finally, using this
last relation in the equation for ${L}$, we arrive at the result
\begin{equation}
{\sum}_{i,j=1}^{M} {Z}_{i}(\mathbf{a}){L}_{,ij}
(\mathbf{x}(\mathbf{a}), \mathbf{a}){Z}_{j}(\mathbf{a})=0. \label{9}
\end{equation}
It is at this point that a local uniqueness condition must be
imposed on the solution $\mathbf{x}(\mathbf{a})$ to insure the
desired invariance property. To avoid inessential technical
complications, we will assume here that the matrix ${L}_{,ij}$ is
negative definite, a condition which is sufficient but not
necessary.  While it is common practice in textbook expositions to
make this assumption throughout, we have only relied on the
semidefiniteness of ${L}_{,ij}$ elsewhere in this paper.  Given the
stronger assumption then, we immediately conclude from Eq.~(\ref{9})
that $\mathbf{Z}$ must vanish, a conclusion which is equivalent to
the invariance property we arrived at heuristically before [and
incorporated in Theorem 2 as Eq.~(\ref{16})].

To demonstrate the effectiveness of this invariance result, let us
consider how it applies to the case of homogeneity, or scale
invariance properties.  Suppose the objective and constraint
functions of a given optimization problem satisfy the invariance
conditions given in Eqs.\ (\ref{8}) with
$\mathbf{X}(\mathbf{x})=\eta \mathbf{x}$ and
$\mathbf{A}(\mathbf{a})= \mathbf{a}$.  Then, according to
Eq.~(\ref{16}), the solution to this problem will satisfy the
invariance condition ${\sum}_{\mu=1}^{N}
{a}_{\mu}{x}_{i,\mu}(\mathbf{a})=\eta {x}_{i}(\mathbf{a})$.  This
last condition characterizes $\mathbf{x}(\mathbf{a})$ as a
homogeneous function of degree $\eta$ (via Euler's Theorem).  As an
example of this, consider the problem of utility maximization
treated in \S IIA, where the underlying symmetry resides in the fact
that a uniform inflation or deflation in income and all the prices
has no effect on the consumption bundle demanded. Formally, this
symmetry corresponds to the fact that the invariance conditions
given in Eqs.\ (\ref{8}) are satisfied by the objective (utility)
and constraint (budget constraint) functions of the problem with
$\mathbf{X}(\mathbf{x})=0$ and $\mathbf{A}(\mathbf{a})= \mathbf{a}$,
hence $\eta =0$. This characterizes the decision function as
homogeneous of degree zero, in full harmony with the underlying
symmetry; in symbols, ${\sum}_{\mu=1}^{N-1} {p}_{\mu} {\partial
\over {\partial {p}_{\mu}}}\mathbf{x}(\mathbf{p},m)+m {\partial
\over
\partial m}\mathbf{x}(\mathbf{p},m)=0$.  An even simpler case of scale
invariance is one which we have already exploited, namely the
introduction and use of the scale parameter $s$ in the basic profit
maximization model in \S IIA.  This invariance can be interpreted in
several economically meaningful ways.  For example, one can consider
the fact that the amount of revenues minus expenditures being
maximized in that model may be expressed in various multiples of a
given currency, or even in different currencies, and this should
have no effect on the decision functions.  Alternatively, one might
interpret $1-s >0$ as a flat tax rate on profits and the objective
function as the net profit of the firm, which again has no effect on
the decision functions.  Mathematically, on the other hand, one
chooses $\mathbf{X}(\mathbf{x})=\mathbf{0}$ and
$\mathbf{A}(\mathbf{a})= (0,0,\ldots,0,s)$, so that
$J(\mathbf{x},\mathbf{a})=s \partial /
\partial s$ and $\eta =0$.  Then the invariance conditions in Eqs.\
(\ref{8}) are satisfied, and (with $s\neq 0$) we have the result
that $
\partial \mathbf{x}(\mathbf{a})/ \partial s=0$, the obvious result.  Other
examples of the use of Eq.~(\ref{16}) will be considered in the
applications.

Having dealt with the envelope and symmetry properties of the
optimization problem defined by Eq.~(\ref{1}), we now proceed to
study several structural features of $\sf{\Omega}$, the general CSM
constructed for that problem.  The first question concerns the
definiteness of $\sf{\Omega}$, or more specifically, whether its
rank is lower than its order, and if so, whether there exists an
upper bound on this rank.  Let us recall that we are dealing with
$M$ decision variables, $N$ parameters, $K$ (independent)
constraints, and a CSM of order $A \leq N$ defined by [the equation
preceding Eq.~(\ref{6})]
\[{\sf{\Omega}}_{\alpha \beta}(\mathbf{a}) \stackrel{\rm{{def}}}{=}-
{\sum}_{i,j=1}^{M}{x}_{i;\alpha}(\mathbf{a}){x}_{j
;\beta}(\mathbf{a})[{f}_{,ij}(\mathbf{x}(\mathbf{a}),\mathbf{a})
+{\sum}_{k=1}^{K}
{\lambda}_{k}(\mathbf{a}){g}_{,ij}^{k}(\mathbf{x}(\mathbf{a}),\mathbf{a})].
\] The task before us is to establish an upper bound for the rank of
this matrix.  Let us first emphasize that an upper bound is all that
can in general be established, since the $M \times A$ matrix
${x}_{i;\alpha}(\mathbf{a})$ that appears in the above equation can
have an arbitrarily small rank, including zero, implying the same
for $\sf{\Omega}$.  For example, for the optimization problem ${\max
\;}_\mathbf{x} {\;}[G(\mathbf{x})+H(\mathbf{a})]$, the decision
functions do not depend on the parameters at all, causing ${x}_{i;
\alpha}(\mathbf{a})$ to vanish identically.  This implies the
vanishing of $\sf{\Omega}$ and its rank.

Let the rank of $\sf{\Omega}$ be denoted by ${\rho}^{\sf{\Omega}}$.
Then an upper bound to ${\rho}^{\sf{\Omega}}$ can be readily derived
from the theorems that (a) the rank of an ${N}^{R} \times {N}^{C}$
matrix cannot exceed $ \min ({N}^{R},{N}^{C})$, and (b) the rank of
a product cannot exceed that of any of its factors.  Now ${x}_{i;
\alpha}(\mathbf{a})$ is a factor for $\sf{\Omega}$, and its rank
cannot exceed $M-K$ on account of the constraints, as an argument
below will confirm.  Therefore, the rank of $\sf{\Omega}$ cannot
exceed the smaller of $M-K$ and $A$, i.e., ${\rho}^{\sf{\Omega}}
\leq \min(M-K,A)$.  To establish the rank property of the matrix
${x}_{i; \alpha}(\mathbf{a})$ just used, recall Lemma 1 of \S IIA
which states that each column of this matrix must conform to the
constraints, i.e., for every value of $\alpha$,
${\sum}_{i=1}^{M}{g}_{,i}^{k}(\mathbf{x}(\mathbf{a}),\mathbf{a})
{x}_{i;\alpha}(\mathbf{a}) =0$, $k=1,2,\ldots,K$.  In other words,
each of the $A$ derivative vectors ${x}_{i;\alpha}(\mathbf{a})$ is
orthogonal to the $K$ normal vectors
${g}_{,i}^{k}(\mathbf{x}(\mathbf{a}),\mathbf{a})$ in decision space.
Now if the constraints are independent, as we have assumed
throughout, then the $K$ normal vectors constitute an independent
set of dimension $K$, forcing the $A$ derivative vectors to be
linearly dependent and at most of dimension $M-K$, since the sum of
these two dimensions cannot exceed $M$, the dimension of the
decision space.  We have thus established that when viewed as a
matrix, ${x}_{i;\alpha}(\mathbf{a})$ can at most have $M-K$ linearly
independent rows.  But this implies that the rank of
${x}_{i;\alpha}(\mathbf{a})$ is no larger than $M-K$, which is the
desired result.

As a first example of the rank inequality formula established above,
consider the profit maximization problem of \S IIA where $K=0$ and
$A=M+1$.  Here, we find ${\rho}^{\sf{\Omega}} \leq \min(M,M+1)=M$,
again implying that the full, $(M+1) \times (M+1)$ CSM will be
singular. However, if following standard practice one uses the $M
\times M$ submatrix of the latter corresponding to the partial
derivatives of the input factors with respect to input prices, there
will no longer be a necessary rank reduction.  As a second example,
consider the Slutsky-Hicks problem, also considered in \S IIA.  Then
with $K=1$, we have ${\rho}^{\sf{\Omega}} \leq \min(M-1,M)=M-1$,
implying that the $M \times M$ Slutsky matrix is necessarily
singular since its order exceeds its rank at least by one, a well
known result.  As a third example, we consider the allocation
problem considered in \S IIA, where $M=G \times H$, $N=L+G+H- 1$,
and $K=G+H-1$.  The resulting $L \times L$ CSM will have have a rank
no larger than $\min((G-1) \times (H-1),L)$.  As a typical
situation, let us take $H \gg L=G \ge 1$.  The corresponding CSM
will then be of order $L$ and rank $L$ or less, and will not convey
detailed comparative statics information about the $G \times H$
optimum allocation levels ${x}_{i}^{h}(\mathbf{a})$ since there are
many more of these than there are rows in the CSM.  Generally
speaking, this state of affairs prevails when the dimensions of the
decision and parameter spaces are widely different. Incidentally,
with $N \gg M$, one ends up with a highly redundant CSM with
${\rho}^{\sf{\Omega}} \ll A$.  It is appropriate at this point to
emphasize that there may very well be further rank reductions of
$\sf{\Omega}$ in specific cases resulting from the special
properties of the objective function, and that the above result
represents the rank reduction that is imposed by the underlying
geometry of the GCD's, independently of the specific properties of
the objective function.

Our next task in this section is a characterization of the
arbitrariness in the result of Theorem 1.  Specifically, we must
consider the problem defined by Eq.~(\ref{1}) et seq., and seek to
classify and characterize all possible CSM's associated with this
problem.  Note that we are taking the set of decision variables as
well as the parameter set as given and fixed, thereby excluding from
the present discussion the arbitrariness associated with these
choices.  It must be emphasized here that although ordinarily there
is a ``natural,'' or ``sensible,'' choice of decision and parameter
sets associated with a given problem, there does exist in principle
the possibility of considering other sets constructed from the given
ones, or even considering smaller or larger sets by, e.g., ignoring
certain parameters as uninteresting or irrelevant and conversely
augmenting the parameter set by introducing auxiliary parameters, or
in the case of decision variables in a constrained problem,
discarding a number of constraint equations by solving for a subset
of the decision variables and conversely.  Moreover, these
alternative choices are not always mere mathematical curiosities
devoid of meaning or use.  Indeed we will exploit these extra
degrees of freedom in our treatment of the models considered in \S\S
IIIA and B, where the utility of the alternative choices of both the
parameter and decision sets will be evident.  However, as already
stipulated, the following discussion will only consider the
arbitrariness resulting from the choice of isovectors and GCD's,
relegating those corresponding to enlarging or contracting the
parameter and decision sets to \S III.

Let us recall, in connection with the problem defined by
Eq.~(\ref{1}) et seq., that we defined and employed a complete set
of GCD's according to $D_{\alpha}(\mathbf{a})
\stackrel{\rm{{def}}}{=}{\sum}_{\mu=1}^{N}{t}^{\alpha}_{\mu}
{\partial / \partial {a}_{\mu}}$.  Recall also that the set of
isovectors $\mathbf{t}^{\alpha}$, $\alpha=1,\ldots,A$, is supposed
to be linearly independent and span the tangent hyperplane (which is
of dimension $A$).  Consider now a different choice for the set of
isovectors, ${\tilde{\mathbf{t}}}^{\alpha}$, $\alpha=1,\ldots,A$,
with the same properties as $\mathbf{t}^{\alpha}$, which must
therefore have the same dimension $A$.  Since, by supposition,
either set is linearly independent and spans the tangent hyperplane
(viewed as a vector space), there must exist a nonsingular matrix
${\sf C}$ of order $A$ that expresses the new set as a linear
combination of the old, and conversely under the inverse matrix
${{\sf C}}^{-1}$; to wit, ${\tilde{t}}^{\alpha}_{\mu}
\stackrel{\rm{{def}}}{=}{\sum}_{\beta=1}^{A} {{\sf C}}_{\alpha
\beta} {t}^{\beta}_{\mu}$.  The set of GCD's corresponding to the
new isovectors is related to the old set according to ${\tilde
{D}}_{\alpha}(\mathbf{a})
\stackrel{\rm{{def}}}{=}{\sum}_{\beta=1}^{A}{{\sf C}}_{\alpha
\beta}D_{\beta}(\mathbf{a})$, with the inverse transformation
effected by means of ${\sf C}^{-1}$.  The transformation rule
between the corresponding CSM's follows straightforwardly from that
of the GCD's.  To wit,
\begin{equation}
{\tilde{\sf{\Omega}}}(\mathbf{a})={\sf
C}{\sf{\Omega}}(\mathbf{a}){{\sf C}}^{\dag}, \label{10}
\end{equation}
where ${\tilde{\sf{\Omega}}}$ is the CSM constructed from the new
isovector set ${\tilde{\mathbf{t}}}^{\alpha}$ and ${\sf C}^{\dag}$
stands for the transpose of $\sf C$.  Again, the inverse
transformation from ${\tilde{\sf{\Omega}}}$ to $\sf{\Omega}$ exists
and is implemented by ${{\sf C}}^{-1}$.  Note that according to
Eq.~(\ref{10}), the two matrices $\sf{\Omega}$ and
${\tilde{\sf{\Omega}}}$ are \textit{congruent}.  As a pair of CSM's,
on the other hand, $\sf{\Omega}$ and ${\tilde{\sf{\Omega}}}$ are
essentially equivalent in the sense that (a) the semidefiniteness of
one implies that of the other, and (b) the two are of equal rank.
Property (a) follows from the fact that the symmetry of
$\sf{\Omega}$ implies that of ${\tilde{\sf{\Omega}}}$, a fact which
readily follows from Eq.~(\ref{10}), together with the observation
that for every pair of real vectors $\mathbf{ v}$ and
${\tilde{\mathbf{v}}}$ related by $\mathbf{v}={{\sf C}}^{\dag}
{\tilde{\mathbf{v}}}$, the quadratic forms $\mathbf{v}^{\dag}
\sf{\Omega} \mathbf{v}$ and ${\tilde{\mathbf{v}}}^{\dag}{\tilde{
\sf{\Omega}}} {\tilde{\mathbf{v}}}$ are equal. Property (b) is a
direct consequence of the nonsingular nature of the transformation
matrix ${\sf C}$.  Intuitively, properties (a) and (b) follow from
the observation that each of the two sets of isovectors from which
$\sf{\Omega}$ and ${\tilde{\sf{\Omega}}}$ are constructed forms a
basis for the tangent hyperplane in the parameter space of the
optimization problem, and as such must provide a description fully
equivalent to the other.  Following standard matrix nomenclature, we
shall refer to a pair of CSM's related according to Eq.~(\ref{10})
as a \textit{congruent} pair. It should be pointed out here,
however, that congruency does not imply similarity of properties as
the two matrices can be quite different with respect to such matters
as observability and empirical verification.  An example of a
congruent but rather dissimilar pair of CSM's for the basic profit
maximization model will be discussed in \S IIIA.

Thus far we have only considered complete sets of GCD's, i.e., those
constructed from a set of isovectors which constitute a basis
for the tangent hyperplane.  At this point we shall relax the
assumption of completeness and explore the consequences of a reduced or
dependent set of isovectors.  Specifically, we shall
characterize the CSM's that are constructed from such incomplete sets,
thereby gaining further insight as to how several CSM's, possibly of
different order and rank, can provide comparative statics information
about the same optimization problem.

Let us first deal with dependent sets of isovectors.  Now a set of
vectors is dependent if at least one member of the set is
expressible as a linear combination of the rest.  It is clear that
GCD's constructed from such sets will be redundant in the sense that
at least one of them will be equal to a linear combination of the
rest.  Mathematically, this redundancy translates into the existence
of a zero eigenvalue and a concomitant reduction of rank for the CSM
(assuming full rank to start with).  Clearly then, no comparative
statics content is lost by discarding the redundancies until the
isovector set becomes independent and none of the GCD's is equal to
a combination of the others.  To restate this argument in
quantitative terms, let us recall that a set of isovectors
${\tilde{\mathbf{t}}} ^{\alpha}$, $\alpha =1,\ldots,A$, is dependent
if there exists a set of scalars ${\eta}_{\alpha}$, not all zero,
such that the vector represented by the linear combination
${\sum}_{\alpha=1}^{A} {\eta}_{\alpha}{\tilde{\mathbf{t}}}^{\alpha}$
vanishes.  But this directly implies the vanishing of
${\sum}_{\alpha=1}^{A}
{\eta}_{\alpha}{\tilde{D}}_{\alpha}(\mathbf{a})$, hence also of
${\sum}_{\alpha=1}^{A} {\eta}_{\alpha} {\tilde{\sf{\Omega}}}_{\alpha
\beta}(\mathbf{a})$, where ${\tilde{D}}_{\alpha}$ and
${\tilde{\sf{\Omega}}}$ respectively represent the set of GCD's and
the CSM constructed from the isovector set ${\tilde{\mathbf{t}}}
^{\alpha}$.  The vanishing of the last-stated sum verifies the claim
that ${\tilde{\sf{\Omega}}}$ has a vanishing eigenvalue, hence a
rank no larger than $A-1$.  Thus, carrying through the
transformation process between the two sets of isovectors, GCD's,
and CSM's as was done above, we arrive at the result that there
exists a \textit{singular} matrix ${\sf E}$, with the property
${\tilde{t}}^{\alpha}_{\mu} ={\sum}_{\beta=1}^{A} {{\sf E}}_{\alpha
\beta} {t}^{\beta}_{\mu}$, such that
\begin{equation}
{\tilde{\sf{\Omega}}}(\mathbf{a})={\sf
E}{\sf{\Omega}}(\mathbf{a}){{\sf E}}^{\dag}. \label{11}
\end{equation}
Note that it follows directly from Eq.~(\ref{11}) that
${\tilde{\sf{\Omega}}}$ is a semidefinite matrix if ${\sf{\Omega}}$
is.

At this point one can discard a dependent isovector, i.e., one
that can be expressed in terms of the others, and proceed to examine the
remaining set for further linear dependences.  It should be clear at this
point that only redundant information is conveyed by means of dependent
isovectors and GCD's.  This is essentially the reason for the
stipulation that our isovector sets be independent.  On the other hand,
it must be emphasized that CSM's resulting from dependent isovector
sets and redundant GCD's are nevertheless valid comparative statics
statements, even though they are not optimally concise and may in fact
fail to capture all the comparative statics information pertinent to the
model in question.

The last remark returns us to the issue of completeness.  Recalling
that we are still dealing with the problem defined by Eq.~(\ref{1})
et seq., we proceed to contemplate a new, smaller set of isovectors
${\tilde{t}}^{\alpha}_{\mu}$, with $\mu =1,\ldots,N$, $\alpha
=1,\ldots,B$, and $B<A$; both sets are assumed to be linearly
independent.  Under these circumstances, we shall refer to
${\tilde{t}}^{\alpha}_{\mu}$ as a \textit{contracted} set and the
transformation as a \textit{contraction}.  This transformation is
implemented by a rectangular matrix ${\sf R}$ of size $B \times A$:
${\tilde{t}}^{\beta}_{\mu} ={\sum}_{\alpha=1}^{A} {{\sf R}}_{\beta
\alpha} {t}^{\alpha}_{\mu}$, $\beta=1,\ldots,B$.  Continuing with
the process of transformation as was carried out above, we arrive at
the result that the $B \times B$ matrix
\begin{equation}
{\tilde{\sf{\Omega}}}(\mathbf{a})={\sf
R}{\sf{\Omega}}(\mathbf{a}){{\sf R}}^{\dag} \label{12}
\end{equation}
is the CSM corresponding to the contracted set.  Again one can show
directly from Eq.~(\ref{12}) that ${\tilde{\sf{\Omega}}}$ is a
semidefinite matrix.  It is clear that the information conveyed by
the contracted CSM is in general reduced relative to the original
CSM.  For example, if ${{\sf{\Omega}}}$ has full rank, then
${\tilde{\sf{\Omega}}}$ will be lower in rank by $A-B$, hence lower
in information content as well.

The last subject to be considered in this section is the effect of
generalized transformations on the CSM.  Here we are contemplating the
possibility of defining new decision variables and parameters in terms of
the old, essentially what amounts to a change of coordinates in the
decision and parameter spaces.  We shall require the transformation
functions to be continuously differentiable with nonvanishing Jacobian
determinants (in appropriate neighborhoods of the respective spaces) in
order to insure a one-to-one mapping between the two sets of coordinates.
While not as directly useful in applications as the specific ones
considered above, these general transformations nevertheless provide a
framework for contemplating a wide class of CSM's associated with an
optimization problem.

Consider the problem defined in Eq.~(\ref{1}) et seq. and the
associated CSM in Eq.~(\ref{7}).  Let us rewrite the latter in the
generic form
\[ {\sf{\Omega}}_{\alpha \beta}(\mathbf{a}) \stackrel{\rm{{def}}}{=}
{\sum}_{i=1}^{M}{\sum}_{\mu=1}^{N}{C}_{i \mu}^{\alpha
\beta}(\mathbf{a}){x}_{i, \mu}(\mathbf{a}). \] This equation in
effect defines the functions ${C}_{i \mu}^{\alpha \beta}(\cdot)$.
Let us now consider a transformation to the generalized decision
variables and parameters $\tilde{\mathbf{x}}$ and
$\tilde{\mathbf{a}}$, which are related to the original ones
according to
\begin{equation}
{x}_{i}\stackrel{\rm{def}}{=}{{\xi}}_{i}({\tilde{\mathbf{x}}}), \:\:
{a}_{\mu}\stackrel{\rm{def}}{=}{{\alpha}}_{\mu}({\tilde{\mathbf{a}}}),
\label{13}
\end{equation}
where $i=1,\ldots,M$ and $\mu=1,\ldots,N$.  The transformation functions
${\xi}_{i}(\cdot)$ and ${\alpha}_{\mu}(\cdot)$ are assumed to obey the
regularity conditions stated above, so that the inverse transformations
${\tilde{\xi}}_{i}(\cdot)$ and ${\tilde{\alpha}}_{\mu}(\cdot)$ exist and
are continuously differentiable functions on ${\Re}^{M}$ and ${\Re}^{N}$
respectively.  An application of the chain rule then leads to
\begin{equation}
{\tilde{\sf{\Omega}}}_{\alpha \beta}({\tilde{\mathbf{a}}})
\stackrel{\rm{def}}{=}
{\sum}_{j=1}^{M}{\sum}_{\nu=1}^{N}{\tilde{C}}^{\alpha \beta}_{j
\nu}({\tilde{\mathbf{a}}}){\tilde{x}}_{j, \nu}({\tilde{\mathbf{a}}})
\label{14}
\end{equation}
as the transformed CSM.  Here ${\tilde{C}}_{j \nu}^{\alpha \beta}
={\sum}_{i=1}^{M}{\sum}_{\mu=1}^{N} [\partial {\xi}_{i} / \partial
{\tilde{x}}_{j}]{C}_{i\mu}^{\alpha \beta}[ \partial {\tilde{\alpha}}_{\nu}
/ \partial {a}_{\mu}]$.  These results clearly demonstrate that the
freedom in deciding the complexion of the comparative statics information
for a given optimization problem is essentially limitless.  Whether any
given complexion is useful or even meaningful for the problem at hand is
an altogether separate question.  We shall illustrate the workings of
these general transformations in \S IIIA and IIIB.

We close this section by summarizing the envelope, invariance, rank, and
transformation properties established above.

\textbf{Theorem 2.} \textit{The following properties hold for the
constrained optimization problem defined by Eq.~(\ref{1}) et seq.:}

\rm{ (i)} \textit{The value function satisfies the envelope
property}
\begin{equation}
{V}_{;\alpha}(\mathbf{a})={f}_{; \alpha}(\mathbf{x}(\mathbf{a}),
\mathbf{a}).\label{15}
\end{equation}
\textit{Furthermore, if the GCD's possess the null property with
respect to the objective function as well, the value function
satisfies the stronger condition ${V}_{;\alpha}(\mathbf{a})=0$.}

\textit{\rm {(ii)} If the optimization problem satisfies the
invariance conditions stated in Eq.~(\ref{8}), then the decision
functions possess the invariance property}
\begin{equation}
{X}_{i}(\mathbf{x}(\mathbf{a}))-{\sum}_{\mu=1}^{N}{A}_{\mu}
(\mathbf{a}){x}_{i,\mu}(\mathbf{a})=0.
\label{16}
\end{equation}

\rm {(iii)} \textit{The rank of $\sf{\Omega}$ is no larger than the
smaller of $M- K$ and $A$, i.e., ${\rho}^{\sf{\Omega}} \leq
\min(M-K,A)$.}

\rm {(iv)} \textit{The matrix ${\tilde{\sf{\Omega}}}$ obtained from
$\sf{\Omega}$ by means of the transformation ${\sf
T}{\sf{\Omega}}(\mathbf{a}){{\sf T}}^{\dag}$, as in Eqs.\
(\ref{10})-(\ref{12}), is semidefinite and constitutes a CSM.
Furthermore, if $\sf T$ is nonsingular, as in Eq.~(\ref{10}), then
${\tilde{\sf{\Omega}}}$ is equivalent to ${{\sf{\Omega}}}$ as a CSM,
whereas if $\sf T$ is singular, as in Eq.~(\ref{11}), or if $T$ is
rectangular representing a contraction, as in Eq.~(\ref{12}), then
${\tilde{\sf{\Omega}}}$ has a rank lower than $A$ and may not be
equivalent to $\sf{\Omega}$ as a CSM.}

rm {{(v)} \textit{If the generalized decision variables and
parameters defined in Eqs.\ (\ref{13}) et seq. are employed in place
of the original ones, then the corresponding CSM is given by
Eq.~(\ref{14}).}

\subsection{Further Developments}

As mentioned in the Introduction, Samuelson (1947) established the
foundations of comparative statics methodologies and Silberberg (1974)
generalized and advanced that work to the point of constructing a
semidefinite matrix conveying the comparative statics
properties of a general optimization problem.  However, Silberberg's
(1974)
construction has a serious shortcoming in dealing with constrained
optimization
problems, namely, the subjection of the said matrix to the constraints.
The construction described above and summarized in Theorem 1 removes this
limitation in a general way.  Naturally, this raises the question of just
how this is accomplished, and the relation, if any, between the present
method and that of Silberberg (1974).  We
will answer this question by deriving the result of Theorem 1 from that of
Silberberg (1974).

To that end, let us recall Silberberg's [1974, Eq.~(\ref{10})]
result as applied to the constrained optimization problem defined by
Eq.~(\ref{1}) et seq.  Stated in our notation, the result is the
statement that the $N \times N$ matrix
\begin{eqnarray}
{{\sf S}}_{\mu \nu}(\mathbf{a}) && \stackrel{\rm{{def}}}{=}
{\sum}_{i=1}^{M} {x}_{i, \nu} (\mathbf{a})[{f}_{,\mu i}
(\mathbf{x}(\mathbf{a}),\mathbf{a}) +{\sum}_{k=1}^{K}
{\lambda}_{k}(\mathbf{a}){g}_{,\mu i}^{k}(\mathbf{x}(\mathbf{a}),\mathbf{a})] \nonumber \\
&& +{\sum}_{k=1}^{K} {\lambda}_{k,\nu}(\mathbf{a}){g}_{,
\mu}^{k}(\mathbf{x}(\mathbf{a}),\mathbf{a})  \label{17}
\end{eqnarray}
is positive semidefinite, \textit{subject to constraints}. Stated
more explicitly, the latter qualification implies that for every
real, $N$-dimensional vector $\mathbf{q}$ such that ${\sum}_{\mu
=1}^{N} {q}_{\mu} {g}_{,
\mu}^{k}(\mathbf{x}(\mathbf{a}),\mathbf{a})=0$, $k=1,\ldots,K$, the
quantity ${\sum}_{\mu, \nu=1}^{N} {q}_{\mu} {{\sf S}}_{\mu
\nu}(\mathbf{a}) {q}_{\nu}$, or $\mathbf{q}^{\dag} {\sf S}
(\mathbf{a})\mathbf{q}$ for short, is nonnegative. Restated in
geometrical terms, the last statement implies that the quadratic
form $\mathbf{q}^{\dag}{\sf S} (\mathbf{a})\mathbf{q}$ is
nonnegative provided that the vector $\bf q$ lies in the parameter
space tangent hyperplane defined by the constraint functions.  On
the other hand, the isovectors $\mathbf{t}^{\alpha}$,
$\alpha=1,\ldots,A$, possess this property by construction (given in
\S IIA), as does any linear combination of them.  Let us consider
then, for any arbitrary real vector ${\eta}_{\alpha}$,
$\alpha=1,\ldots,A$, the linear combination
${\sum}_{\alpha=1}^{A}{\eta}_{\alpha} \mathbf{t}^{\alpha}$, and use
this in place of $\bf q$ in the quadratic form above.  The resulting
expression can in turn be viewed as another quadratic form composed
of the $A \times A$ symmetric matrix $\mathbf{t}^{\dag \alpha}{\sf
S} (\mathbf{a})\mathbf{t}^{\beta}$ and the arbitrary vector
${\eta}_{\alpha}$.  But then the arbitrary nature of
${\eta}_{\alpha}$ permits one to conclude that $\mathbf{t}^{\dag
\alpha} {\sf S}(\mathbf{a})\mathbf{t}^{\beta}$ is positive
semidefinite. If the expression for ${\sf S} (\mathbf{a})$ given in
Eq.~(\ref{17}) is now substituted in $\mathbf{t}^{\dag \alpha} {\sf
S}(\mathbf{a})\mathbf{t}^{\beta}$ and the symmetry of the
second-order partial derivatives of $f$ and ${g}^{k}$ are used, the
parameter space partial derivatives in that expression turn into
GCD's by virtue of the identity $\mathbf{t}^{\alpha} \cdot
{\nabla}^\mathbf{a}=D_{\alpha}(\mathbf{a})$.  This causes the last
term in the expression for ${\sf S} (\mathbf{a})$ to drop out of the
resulting form by virtue of Lemma 1, leaving behind an expression
that is found to be identical with the CSM matrix $\sf{\Omega}$
given in Eq.~(\ref{7}).  This completes the deduction of Theorem 1
from Silberberg's (1974) result.

Next we consider the work of Hatta (1980) referred to in the
Introduction, and show that its main comparative statics result is
in fact a special case of Theorem 1.  The optimization problem
treated by Hatta (1980), given in his Eq.~(\ref{10}), is a special
case of our Eq.~(\ref{1}), and appears as
\[ {\max \;}_\mathbf{x} {\;}f(\mathbf{x},\mathbf{p}) {\:} {\:} s.t. {\:} {\:}
{g}^{l}(\mathbf{x},\mathbf{p},{\sf{\kappa}})\stackrel{\rm{{def}}}{=}{\kappa}_{l}-
{k}^{l}(\mathbf{x},\mathbf{p})=0, {\:}{\:}l=1,2,\ldots,K \] in our
notation, where $\bf x$ and $\bf p$ are $M$-dimensional vectors. The
crucial property of this special form is the occurrence of the
parameters ${\sf{\kappa}}$ in a separable, linear manner in the
constraint equations, and their absence from the objective function.
This special structure makes it possible to construct a set of GCD's
patterned after those customarily used for the Slutsky-Hicks problem
(cf., \S IIA or \S IIIB), namely
$D_{\alpha}(\mathbf{x},\mathbf{a})\stackrel{\rm{{def}}}{=}{\partial
/
\partial {p}_{\alpha}} +{\sum}_{l=1}^{K} [\partial {k}^{l}(\mathbf{x},
\mathbf{p}) / \partial {p}_{\alpha} ] {\partial / \partial
{\kappa}_{l}}$.  Note that the parameter set here is identified as
$\mathbf{a}\stackrel{\rm{{def}}}{=} (\mathbf{p}, \sf{\kappa})$.
Using these GCD's in Eq.~(\ref{7}), we find the result that the
matrix
\begin{eqnarray}
{\sum}_{i=1}^{M}[{f}_{,i
\alpha}(\mathbf{x}(\mathbf{p},&&{\sf{\kappa}}),\mathbf{p})
-{\sum}_{l=1}^{K} {\lambda}_{l}(\mathbf{p},{\sf{\kappa}}){k}^{l}_{,i
\alpha}(\mathbf{x}(\mathbf{p},{\sf{\kappa}}),\mathbf{p})] \{
\partial {x}_{i}(\mathbf{p},{\sf{\kappa}}) /
\partial {p}_{\beta}  \nonumber \\
&&+{\sum}_{l=1}^{K} [\partial {k}^{l}(\mathbf{x}, \mathbf{p}) /
\partial {p}_{\beta} ] \partial {x}_{i}(\mathbf{p},{\sf{\kappa}})/
\partial {\kappa}_{l} \}, \nonumber
\end{eqnarray}
is positive semidefinite.  This statement is the same as Hatta's
(1980) Theorems 6 and 7, his main comparative statics results.  Note
that because of the special structure of the problem, compensation
terms appear only in the partial derivatives of the decision
functions, $\mathbf{x}_{;\alpha}(\mathbf{p},{\sf{\kappa}})
=D_{\alpha}(\mathbf{x},\mathbf{a})\mathbf{x}(\mathbf{p},{\sf{\kappa}})$.
These compensated derivatives are denoted by ${s}_{p}(\mathbf{p},
\mathbf{x}^{\star}(\mathbf{p}, {\sf{\kappa}}))$ and termed ``the
\textit{Slutskian} substitution matrix'' by Hatta (1980), while the
Lagrange multipliers ${\lambda}_{l}$ are represented by
${\phi}_{\kappa}$ in his notation.  An examination of the manner in
which the quantities ${s}_{p}(\mathbf{p},
\mathbf{x}^{\star}(\mathbf{p}, {\sf{\kappa}}))$ are derived by Hatta
(1980), on the other hand, reveals that they are constructed in
conformity to the constraints, i.e., precisely according to the
definition of our GCD's, although this property is obscured by the
presentation. Moreover, the method of their construction
specifically relies on the special role played by the parameters
$\sf{\kappa}$ and is therefore limited to the assumed form of the
problem.  This completes the discussion of how our method relates to
the existing ones for dealing with constrained problems.

Next we turn to developing the results of Theorem 1 in more detail
for certain generic forms that naturally arise within the context of
economic problems.  This will also allow us to introduce simplified
variants of the construction methods given in \S\S IIA and B.  The
first such form arises in the case of unconstrained optimization, to
wit, ${\max \;}_\mathbf{x} {\;}f(\mathbf{x},\mathbf{a})$.  As
discussed in \S IIA, it is expedient to deal instead with the
modified problem ${\max \;}_\mathbf{x}
{\;}\tilde{f}(\mathbf{x},\mathbf{b})$, where
$\tilde{f}(\mathbf{x},\mathbf{b})\stackrel{\rm{{def}}}{=}sf(\mathbf{x},\mathbf{a})$,
and where $ s > 0$ is viewed as the ($N+1$)st parameter in the
modified problem.  Since the decision functions for these two
problems must be the same, we have the statement that ${\partial
\over \partial s}\tilde{\mathbf{x}}(\mathbf{b})=0$ for the modified
problem, so that we can use $\mathbf{x}(\mathbf{a})$ to denote the
decision functions for both problems without any fear of confusion.
Obviously, other quantities of interest in the original problem can
be recovered from those of the modified version by setting the
auxiliary parameter $s$ equal to unity.

An important advantage of introducing the auxiliary parameter is
that it allows a natural choice of GCD's, as well as a simple method
of constructing them.  The idea is to start with the set of partial
derivatives ${\partial \over \partial {a}_{\alpha}}$ and linearly
combine each of these with ${\partial \over \partial s}$ in such a
way that the application of the resulting combination to the
(modified) objective function $\tilde{f}(\mathbf{x},\mathbf{b})$
returns zero \textit{identically}.  Applying this prescription, one
can readily construct the set of GCD's given by
${D}_{\alpha}(\mathbf{x},\mathbf{b})=\tilde{f}'(\mathbf{x},\mathbf{b})
{\partial \over
\partial {a}_{\alpha}}-\tilde{f}_{,\alpha}(\mathbf{x},\mathbf{b}){\partial \over
\partial s}$, $\alpha =1,\ldots, N$, where $\tilde{f}'(\mathbf{x},\mathbf{b})\stackrel{\rm{{def}}}{=}{\partial \over \partial s}\tilde{f}(\mathbf{x},\mathbf{b})$.  This brute force method of constructing GCD's will be
referred to as the \textit{one-term compensation} method.  It is
simple and direct, and will suffice when there is only one target
function, such as in the present case, or the case of constrained
optimization with only one constraint and no stipulation that the
GCD's have the null property with respect to the objective function
of the problem.  Needless to say, this already encompasses a number
of important cases of interest in economic problems.

The next step is the construction of the CSM, using the GCD's
defined above by the method of one-term compensation.  First, let us
note that $\tilde{f}'(\mathbf{x},\mathbf{b}) =
f(\mathbf{x},\mathbf{a})$, and
${x}_{i;\alpha}(\mathbf{a})=f(\mathbf{x}(\mathbf{a}),\mathbf{a}){x}_{i,\alpha}(\mathbf{a})$
since the decision functions are independent of the auxiliary
parameter.  On the other hand,
$\tilde{f}_{,i;\beta}(\mathbf{x},\mathbf{b})={f}(\mathbf{x},\mathbf{a})sf_{,i
\beta}(\mathbf{x},\mathbf{a})-s{f}_{,\beta}(\mathbf{x},\mathbf{a})
{f}_{,i}(\mathbf{x},\mathbf{a})$, $\beta =1,\ldots, N$.  Restricting
all of these expressions to the solution
$\mathbf{x}=\mathbf{x}(\mathbf{a})$ and substituting them in
Eq.~(\ref{7}) while setting the auxiliary parameter $s$ equal to
unity, we find, upon some rearrangement, the matrix
\begin{equation}
{\sf{{\Omega}^{ A1}}}_{\mu \nu}(\mathbf{a}) \stackrel{\rm{{def}}}{=}
f(\mathbf{x}(\mathbf{a}),\mathbf{a}) {\sum}_{i=1}^{M}{x}_{i,
\nu}(\mathbf{a}) {[\log f(\mathbf{x}(\mathbf{a}),\mathbf{a})]}_{,i
\mu} \label{18}
\end{equation}
as a CSM for the original problem.

A useful variant of ${\sf{{\Omega}^ { A1}}}$ results upon carrying
out the implied differentiations in Eq.~(\ref{18}) and using the
first-order condition.  However, it is more instructive to derive
this variant directly from Theorem 1.  To do so, let us recall that
Theorem 1 requires the null property of the GCD's only with respect
to the constraint functions, leaving the same property with respect
to the objective function as an option.  Therefore, Theorem 1 can be
applied to \textit{unconstrained} problems provided that GCD's are
replaced with \textit{ordinary} partial derivatives.  This is
because, for unconstrained problems, the null property required of
the GCD's with respect to the constraint functions is vacuously
satisfied by ordinary partial derivatives.  Therefore, we have a
corollary of Theorem 1 stating that the matrix
\begin{equation}
{\sf{{\Omega}^{ A2}}}_{\mu \nu}(\mathbf{a}) \stackrel{\rm{{def}}}{=}
{\sum}_{i=1}^{M}{x}_{i, \nu}(\mathbf{a}){f}_{,i
\mu}(\mathbf{x}(\mathbf{a}),\mathbf{a}) \label{19}
\end{equation}
is a CSM.  This is the desired variant mentioned above, a form
that also follows directly from Silberberg's (1974) general theorem
discussed
earlier.  Although Eqs.\ (\ref{18}) and \ (\ref{19}) are trivially
related, it is useful to record them both since they do lead to different
forms of CSM's.

\textbf{Corollary A.} \textit{If $K=0$ in Theorem 1, i.e., in the
absence of constraints, ${\sf{{\Omega}^{ A1}}}$ and ${\sf{{\Omega}^{
A2}}}$ given by Eqs.\ (\ref{18}) and \ (\ref{19}) respectively, are
comparative statics matrices.}

Returning now to the general problem posed in Eq.~(\ref{1}), let us
again introduce the auxiliary parameter $s$ and treat the modified
problem ${\max \;}_\mathbf{x} {\;}sf(\mathbf{x},\mathbf{a}){\:} s.t.
{\:} {g}^{k}(\mathbf{x},\mathbf{a})=0$, $k=1,\ldots,K$.  We start by
using the general construction scheme of \S IIA to develop a set of
GCD's, denoted by $D_{\alpha}(\mathbf{x},\mathbf{a})$ as before,
with respect to the constraint functions
${g}^{k}(\mathbf{x},\mathbf{a})$; the associated null property is
$D_{\alpha}(\mathbf{x},\mathbf{a}){g}^{k}(\mathbf{x},\mathbf{a})=0$,
$\alpha=1,\ldots,A$ and $k=1,\ldots,K$. Next, we use the auxiliary
variable $s$ in conjunction with the one-term compensation method to
extend the null property just stated to the objective function.  In
this manner we arrive at the set of modified GCD's
$\tilde{D}_{\alpha}(\mathbf{x},\mathbf{b})\stackrel{\rm{{def}}}{=}f(\mathbf{x},\mathbf{a})D_{\alpha}(\mathbf{x},\mathbf{a})-s{f}_{;
\alpha}(\mathbf{x},\mathbf{a}){\partial \over {\partial s}}$, where
$\mathbf{b}=(\mathbf{a}, s)$ and a semicolon represents compensated
differentiation with respect to the unmodified GCD's, as before.
These modified GCD's have the null property with respect to the
modified objective function as well as the constraint functions:
$\tilde{D}_{\alpha}(\mathbf{x},\mathbf{b})\tilde{f}(\mathbf{x},\mathbf{b})=\tilde{D}_{\alpha}(\mathbf{x},\mathbf{b}){g}^{k}(\mathbf{x},\mathbf{a})=0$,
$\alpha=1,\ldots,A$ and $k=1,\ldots,K$, where
$\tilde{f}(\mathbf{x},\mathbf{b})\stackrel{\rm{{def}}}{=}sf(\mathbf{x},\mathbf{a})$
as before.

Our next step is the construction of the CSM, remembering that the
decision functions are independent of the auxiliary parameter, and that
results appropriate to the original problem are recovered from those of
the modified version upon setting $s$ equal to unity.  Carrying out the
necessary steps, we find, after some rearrangement,
\begin{equation}
{\sf{{\Omega}^{ B}}}_{\alpha
\beta}(\mathbf{a})\stackrel{\rm{{def}}}{=}f(\mathbf{x}(\mathbf{a}),\mathbf{a})
{\sum}_{i=1}^{M}{x}_{i; \beta}(\mathbf{a}) \{ {[\log
f(\mathbf{x}(\mathbf{a}),\mathbf{a})]}_{,i;\alpha} +{\sum}_{k=1}^{K}
{\lambda}_{k}(\mathbf{a}){g}_{,i;
\alpha}^{k}(\mathbf{x}(\mathbf{a}),\mathbf{a}) \} . \label{20}
\end{equation}
Again, as in the case of Eq.~(\ref{18}), a variant of this equation
can be derived by applying the first-order condition in combination
with Lemma 1 established in \S IIA; this variant turns out to be
$\sf{\Omega}$ of Theorem 1.  Both variants are useful, however, and
can lead to rather different forms of CSM's.  Therefore, we will
record the result in Eq.~(\ref{20}) as a corollary to Theorem 1.

\textbf{Corollary B.} \textit{Under the conditions of Theorem 1,
${\sf{{\Omega}^{ B}}}$ given in Eq.~(\ref{20}) is a comparative
statics matrix.}

Let us mention in passing that for $K=1$, the case of one constraint, the
method of one-term compensation can be usefully employed in the
construction of $\sf{\Omega}$ and ${\sf{{\Omega}^{ B}}}$.

In \S IIB we discussed the large arbitrariness in the choice of the CSM's
associated with a given optimization problem, and have otherwise provided
examples of this arbitrariness in the construction of isovector
sets and GCD's.  We will close this subsection by exploring and
highlighting
the invariant characteristics that lie at the root of all these CSM's.
Not surprisingly, it is again the underlying geometric structure that is
most effective in providing an intuitive picture of what is going on.
Since the characteristics we seek are essentially the same for constrained
and unconstrained problems, we shall discuss them within the latter
context for simplicity.

Consider an unconstrained optimization problem with the objective
function $f(x, \mathbf{a})$ and $M=1$, i.e., one decision variable,
endowed with the regularity properties set out in \S IIA.  Suppose
further that there exists a local maximum of $f$ as a function of
$x$ at some interior point $x(\mathbf{a})$.  Clearly, for a fixed
value of $\bf a$, $f(x,\mathbf{a})$ is concave at $x=x(\mathbf{a})$,
and ${f}_{,11}(x(\mathbf{a}),\mathbf{a})$ is a negative semidefinite
number equal to the \textit{curvature} of the graph of $f$ at the
maximum point.  Next, consider the analogous situation in $M$
dimensions, with $\mathbf{x} \in {\Re}^{M}$ and
${f}_{,i}(\mathbf{x}(\mathbf{a}),\mathbf{a})=0$, where
$\mathbf{x}(\mathbf{a})$ is the maximizing point. The concavity
condition is now equivalent to the negative semidefiniteness of the
matrix of second-order partial derivatives,
${f}_{,ij}(\mathbf{x}(\mathbf{a}),\mathbf{a})$.  Since this matrix
is real and symmetric, an equivalent condition is the negative
semidefiniteness of its eigenvalues, the set of which is known as
the \textit{spectrum} of the matrix. Geometrically, these
eigenvalues can be thought of as the \textit{principal curvatures}
of the surface ${x}_{M+1}-f(\mathbf{x},\mathbf{a})=0$ in ${\Re
}^{M+1}$ at the maximum point
$(f(\mathbf{x}(\mathbf{a}),\mathbf{a}),\mathbf{x}(\mathbf{a}))$, as
may be surmised by analogy with the one-dimensional case. Note that
this surface is the graph of $f$ as a function of $\mathbf{x}$, as
suggested by analogy to the one-dimensional case.  To see this
analogy more clearly, let us imagine, for $M=2$, the surface of a
smooth hilltop as representing the graph, with elevation
representing the value of the objective function and the apex as the
maximum point.  Each vertical plane containing the apex will
intersect the surface of the hilltop in a curve which has a maximum
at the apex, and a curvature associated with the maximum.  The set
of curvatures so defined has a maximum, which is the largest of the
principal curvatures, and the associated vertical plane is a
\textit{principal plane}.  Next, we restrict the set of vertical
planes to those orthogonal to the principal plane(s) already found,
and repeat the procedure.  In two dimensions, of course, there is
only one such plane, the curvature associated with which is the
minimum principal curvature. For $M>2$, this process continues
through $M$ steps, culminating in a nonincreasing sequence of $M$
principal curvatures. These are the eigenvalues constituting the
spectrum of ${f}_{,ij}(\mathbf{x}(\mathbf{a}),\mathbf{a})$.  Note
that for a spherical hilltop ($M=2$), the two curvatures are equal,
corresponding to a spectrum composed of two equal eigenvalues.  The
most general case for $M=2$, on the other hand, is that of an
ellipsoidal hilltop, corresponding to an unequal pair of nonpositive
eigenvalues.  The singular cases here include one vanishing
eigenvalue, corresponding to a cylindrical hilltop, and two
vanishing eigenvalues, corresponding to a flat hilltop.  To avoid
misunderstanding, it should be stated that the curvatures considered
here are the so-called \textit{extrinsic} curvatures, not to be
confused with the \textit{intrinsic}, or \textit{Gaussian}
curvatures.

Using this geometrical construction, one can reformulate the
second-order conditions by stating that the set of principal
curvatures, equivalently the spectrum of the objective function at
the maximum point, must be nonpositive.  It is this semidefiniteness
of the spectrum of ${f}_{,ij}(\mathbf{x}(\mathbf{a}),\mathbf{a})$
that underlies all comparative statics results. This is seen most
directly by noting that, in the absence of constraints, the matrix
$\sf{\Omega}$ of Eq.~(\ref{7}) is related to
${f}_{,ij}(\mathbf{x}(\mathbf{a}),\mathbf{a})$ according to
\begin{equation}
{\sf{\Omega}}_{\mu \nu}(\mathbf{a}) = -{\sum}_{i,j=1}^{M}{x}_{i,
\mu}(\mathbf{a}){f}_{,ij}(\mathbf{x}(\mathbf{a}),\mathbf{a}) {x}_{j,
\nu}(\mathbf{a}). \label{21}
\end{equation}
First, observe that this relation shows ${\sf{\Omega}}_{\mu
\nu}(\mathbf{a})$ to be symmetric as a direct consequence of the
symmetry of ${f}_{,ij}(\mathbf{x}(\mathbf{a}),\mathbf{a})$.  This
property, on the other hand, guarantees that ${\sf{\Omega}}_{\mu
\nu}(\mathbf{a})$ has a spectrum consisting of $A$ real eigenvalues,
together with a set of corresponding eigenvectors.  We will denote
these by ${\mu}^{(\gamma)}$ and $\mathbf{z}^{(\gamma)}$,
respectively. In other words, we have ${\sum}_{\nu=1}^{A}
{\sf{\Omega}}_{\mu \nu}(\mathbf{a}){z}^{(\gamma)}_{\nu}
={\mu}^{(\gamma)}{z}^{(\gamma)}_{\mu}$ for all values of $\gamma$
and $\mu$. Without loss of generality, we can assume that the
eigenvectors have unit length, i.e., that ${\sum}_{\mu=1}^{A}
{z}^{(\gamma)}_{\mu}{z}^{(\gamma)}_{\mu}=1$. Now multiply
Eq.~(\ref{21}) on both sides by the eigenvector
$\mathbf{z}^{(\gamma)}$; the result is ${\mu}^{(\gamma)} =-
{\sum}_{i,j=1}^{M}{q}_{i}^{(\gamma)}
{f}_{,ij}(\mathbf{x}(\mathbf{a}),\mathbf{a}) {q}_{j}^{(\gamma)}$,
where ${q}_{i}^{(\gamma)}\stackrel{\rm{{def}}}{=}{\sum}_{\mu=1}^{A}
{z}^{(\gamma)}_{\mu} {x}_{i, \mu}(\mathbf{a})$.  Finally, if we
denote the $M$ eigenvalues (or principal curvatures) and
(normalized) eigenvectors of
${f}_{,ij}(\mathbf{x}(\mathbf{a}),\mathbf{a})$ by ${m}^{(I)}$ and
$\mathbf{b}^{(I)}$, we find from the last equation
\begin{equation}
{\mu}^{(\gamma)} = -{\sum}_{I=1}^{M}{(\mathbf{q}^{(\gamma)}\cdot
\mathbf{b}^{(I)})}^{2} {m}^{(I)}, \label{22}
\end{equation}
where we have used the spectral decomposition
${f}_{,ij}={\sum}_{I=1}^{M} {m}^{(I)}{b}_{i}^{(I)}{b}_{j}^{(I)}$ to
replace the Hessian matrix in favor of its invariant
characteristics.  Equation (\ref{22}) clearly shows that the
eigenvalues of $\sf{\Omega}$ are linear combinations of those of
${f}_{,ij}(\mathbf{x}(\mathbf{a}),\mathbf{a})$ with
\textit{nonpositive} coefficients.  Inasmuch as the spectrum of
${f}_{,ij}(\mathbf{x}(\mathbf{a}),\mathbf{a})$ is known to be
nonpositive, one is assured of the nonnegativity of the spectrum of
$\sf{\Omega}$.

Recall now that ${f}_{,ij}(\mathbf{x}(\mathbf{a}),\mathbf{a})$ and
its spectrum $ \{ {m}^{(I)} \}$ are basically unique for a given
problem (and a given local maximum, if there is more than one).  The
corresponding CSM, on the other hand, can be constructed in an
infinite variety of ways.  For example, one may recall from Theorem
2-(iv) that any matrix of the form ${\sf T}{\sf{\Omega}}{\sf
T}^{\dag}$ is a CSM for a given model if $\sf{\Omega}$ is.  All of
these matrices, and many more, qualify as CSM's by virtue of the
fact that the spectrum $ \{ {m}^{(I)} \}$ is nonpositive, as is
evident in Eq.~(\ref{22}).  In short, the concavity of the objective
function, or the nonpositivity of the spectrum of
${f}_{,ij}(\mathbf{x}(\mathbf{a}),\mathbf{a})$, is the invariant
property that guarantees the existence of a great variety of CSM's.

\subsection{A Universal Comparative Statics Matrix}

The preceding discussion has underlined the fact that all
comparative statics information stemming from the underlying
optimization hypothesis originates from the semidefinite spectrum of
principal curvatures at the optimum point.  Given a constrained
optimization problem, then, there arises the question of whether
there exists a universal construction that yields the corresponding
comparative statics information exhaustively (i.e., from which any
other CSM can be derived)?  The answer is affirmative, as we will
now show.  Specifically, we shall construct a comparative statics
matrix for the problem defined in Eq.~(\ref{1}) without any recourse
to the geometry of the parameter space or generalized compensated
derivatives.  This construction is almost entirely anchored in
decision space, and achieves its results by relying on intrinsic,
projective techniques.  The result is a constraint-free, maximal
CSM, i.e., a semidefinite matrix with the highest possible rank for
the given parameter set.  We shall refer to this matrix as the
\textit{universal} CSM. The maximal property of the universal CSM
implies that any other CSM which uses the same parameter set or a
subset thereof, or indeed any comparative statics statement
derivable from the underlying optimization problem which is
expressible in terms of the said parameter set is implied by the
universal CSM.  Not surprisingly, the universal CSM is not as
convenient to construct or apply as those constructed by the method
of generalized compensated derivatives utilized throughout this
work.  From a theoretical point of view, on the other hand, the
existence theorem for the universal matrix established below is a
more fundamental result than that of Theorem 1 (which it implies),
insofar as it provides a constructive proof for the existence of a
constraint-free, \textit{maximal} comparative statics matrix for a
general optimization problem.

Let us recall the general optimization problem
$$ {\max \;}_\mathbf{x} {\;}f(\mathbf{x},\mathbf{a}) {\:} {\:} s.t. {\:} {\:}
{g}^{k}(\mathbf{x},\mathbf{a})=0, {\:}k=1,2,\ldots,K, $$ subject to
the regularity and existence conditions stipulated in \S IIB. As
elsewhere in this work, we shall assume that the constraints are
independent at the solution point, i.e., that the set of $K$ normal
vectors in decision space, ${g}_{,i}^{k}(\mathbf{x}(\mathbf{a}),
\mathbf{a})$, $k=1,2,\ldots,K$, is linearly independent (which
assumption, it may be recalled, is a constraint qualification
condition). Equivalently, the symmetric, $K \times K$ matrix ${\cal
G}_{kk'} (\mathbf{x}(\mathbf{a}), \mathbf{a})
\stackrel{\rm{{def}}}{=}
{\sum}_{i=1}^{M}{g}_{,i}^{k}(\mathbf{x}(\mathbf{a}),
\mathbf{a}){g}_{,i}^{k'}(\mathbf{x}(\mathbf{a}), \mathbf{a})$
possesses full rank and is therefore invertible.  Thus the normal
hyperplane (generated by the set of normal vectors) in decision
space is of dimension $K$, and its orthogonal complement, the
tangent hyperplane, has dimension $M-K$.  Consequently, any decision
space vector $\bf w$ possesses a unique orthogonal decomposition
$\mathbf{w}=\mathbf{w}^{n}+\mathbf{w}^{t}$ corresponding to the
above decomposition.  According to standard results of linear
algebra, there exists a symmetric, idempotent, $M \times M$
\textit{projection} matrix ${\sf Q}$ such that $\mathbf{w}^{n}={\sf
Q}\mathbf{w}$ and $\mathbf{w}^{t}=({\sf I}-{\sf Q})\mathbf{w}$,
where $\sf I$ is the identity matrix.  The construction of ${\sf Q}$
may be accomplished by first establishing an orthonormal basis
$\mathbf{e}^{k}$, $k=1,2,\ldots,K$, for the normal hyperplane (e.g.,
by orthonormalizing the set of normal vectors ${g}_{,i}^{k}$) and
then using the formula ${\sf
Q}_{ij}={\sum}_{k=1}^{K}\mathbf{e}_{i}^{k} \mathbf{e}_{j}^{k} $. The
end result of that procedure is the compact formula ${\sf
Q}_{ij}(\mathbf{x}(\mathbf{a}), \mathbf{a})=
{\sum}_{k,k'=1}^{K}{g}_{,i}^{k} (\mathbf{x}(\mathbf{a}),
\mathbf{a}){\cal G}^{-1}_{kk'} (\mathbf{x}(\mathbf{a}),
\mathbf{a}){g}_{,j}^{k'} (\mathbf{x}(\mathbf{a}), \mathbf{a})$.
Indeed one can verify by direct calculation that $\sf Q$ as given
satisfies the stated properties, and therefore that $\sf Q$ and
${\sf I}-{\sf Q}$ are the projection matrices onto the normal and
tangent hyperplanes, respectively.

At this point we recall from \S IIB that the matrix
$${\sf L}_{ij}(\mathbf{x}(\mathbf{a}),\mathbf{a}) \stackrel{\rm{{def}}}{=}
{L}_{,ij}(\mathbf{x}(\mathbf{a}),\mathbf{a})={f}_{,ij}(\mathbf{x}(\mathbf{a}),\mathbf{a})
+{\sum}_{k=1}^{K} {\lambda}_{k}(\mathbf{a})
{g}_{,ij}^{k}(\mathbf{x}(\mathbf{a}), \mathbf{a})$$ is negative
semidefinite subject to the constraints. In other words,
$\mathbf{v}^{\dag}{\sf L}\mathbf{v} \leq 0$ provided $\bf v$ is an
isovector in decision space, i.e., a vector in the (decision space)
tangent hyperplane. Since for any $M$-vector $\bf w$, $({\sf I}-{\sf
Q})\mathbf{w}$ is an isovector, we have the result that $({\sf
I}-{\sf Q}){\sf L}({\sf I}-{\sf Q})$ is a negative semidefinite
matrix.

Next we differentiate the first order conditions
${L}_{,i}={f}_{,i}+{\sum}_{k=1}^{K} {\lambda}_{k} {g}_{,i}^{k}=0$
restricted to the optimum point $\mathbf{x}=\mathbf{x}(\mathbf{a})$
with respect to ${a}_{\mu}$, and left-multiply the result by the
projection matrix ${\sf I}-{\sf Q}$, also restricted to the optimum
point.  The result is
$${\sum}_{i=1}^{M} \{ {({\sf I}-{\sf Q})}_{li}[{L}_{,i \mu}
+{\sum}_{j=1}^{M}{L}_{,ij}{x}_{j, \mu}] \} =0, $$ where we have left
out the arguments of various functions for brevity. The next step is
to left-multiply the last equation by ${\mathbf{x}_{,\nu}}^{\dag}$,
the partial derivative of the decision vector with respect to
${a}_{\nu}$.  This multiplication leads to
$${\sum}_{i,l=1}^{M} {x}_{l, \nu} \{ {({\sf I}-{\sf Q})}_{li}[{L}_{,i \mu}
+{\sum}_{j=1}^{M}{L}_{,ij}{\sum}_{m=1}^{M}{({\sf Q}+{\sf I}-{\sf
Q})}_{jm}{x}_{m, \mu}] \}=0, $$ where we have also inserted the
identity matrix in the form ${({\sf Q}+{\sf I}-{\sf Q})}_{jm}$
between $\sf L$ and ${\mathbf{x}_{,\mu}}$ in the above equation.
With this insertion, we have arranged for the appearance of the
matrix ${\mathbf{x}_{,\nu}}^{\dag}({\sf I}-{\sf Q}){\sf L}({\sf I}-
{\sf Q})\mathbf{x}_{,\mu}$ in the equation, a matrix which is
negative semidefinite by virtue of the fact, established in the
preceding paragraph, that ${({\sf I}-{\sf Q})}{\sf L}({\sf I}-{\sf
Q})$ is negative semidefinite.  The remaining terms in the above
equation must therefore constitute a positive semidefinite matrix.
Thus the matrix
$${\sf U}_{\nu \mu}\stackrel{\rm{{def}}}{=}{\sum}_{i,l=1}^{M} {x}_{l, \nu}
\{ {({\sf I}-{\sf Q})}_{li}[{L}_{,i \mu}
+{\sum}_{j=1}^{M}{L}_{,ij}{\sum}_{m=1}^{M}{\sf Q}_{jm} {x}_{m, \mu}] \}$$
is positive semidefinite (which, it may be recalled, implies that it is
symmetric as well).

To convert the above form for ${\sf U}$ to a CSM, we must now eliminate
the decision functions from the term ${\sum}_{m=1}^{M}{\sf Q}_{jm} {x}_{m,
\mu}$.  To that end, we use the definition of $\sf Q$ to write
${\sum}_{m=1}^{M}{\sf Q}_{jm} {x}_{m,
\mu}={\sum}_{m=1}^{M}{\sum}_{k,k'=1}^{K}{g}_{,j}^{k} {\cal G}^{-1}_{kk'}
{g}_{,m}^{k'}{x}_{m, \mu}$.  Then, using the identity ${g}^{k}_{, \mu} +
{\sum}_{m=1}^{M}{g}^{k}_{,m}{x}_{m,\mu} =0$ which results from a
differentiation of the constraint functions restricted to the optimum
point with respect to ${a}_{\mu}$, we arrive at

\textbf{Theorem 3.} \textit{The constrained optimization problem
defined by Eq.~(\ref{1}) et seq. admits of the following
constraint-free, positive semidefinite, comparative statics matrix:}
\begin{eqnarray}
{\sf U}_{\mu
\nu}(\mathbf{a})=&&{\sum}_{i,j=1}^{M}{x}_{i,\mu}(\mathbf{a}){[{\sf
I}-{\sf Q}(\mathbf{x}(\mathbf{a}),\mathbf{a})]}_{ij}
\{ {L}_{,j \nu}(\mathbf{x}(\mathbf{a}),\mathbf{a}) - \nonumber \\
&&
{\sum}_{l=1}^{M}{L}_{,jl}(\mathbf{x}(\mathbf{a}),\mathbf{a}){\sum}_{k,
k' =1}^{K}{g}^{k}_{,l}(\mathbf{x}(\mathbf{a}),\mathbf{a}){\cal
G}^{-1}_{kk'}(\mathbf{x}(\mathbf{a}),\mathbf{a}){g}^{k'}_{,\nu}(\mathbf{x}(\mathbf{a}),\mathbf{a})
\}. \label{023}
\end{eqnarray}
\textit{The rank of ${\sf U}$ is no larger than the smaller of $M-K$
and $N$, i.e., ${\rho}^{\sf U} \leq \min (M-K,N)$.}

The rank property of $\sf U$ stated above is readily established by
an argument patterned after the one leading to Theorem 2-(iii).  The
maximal nature of $\sf U$, on the other hand, is easily ascertained
by noting that the rank of any CSM associated with the underlying
optimization problem is limited by (i) the dimension of the
\textit{constrained} decision space, namely $M-K$, and (ii) the size
of the parameter set, namely $N$.  Thus no CSM can be of higher rank
than the smaller of these two integers, which shows that $\sf U$ has
the highest possible rank.  It is worth recalling here that zeros in
the spectrum of $\sf L$ or other ``exceptional'' circumstances may,
in specific cases, reduce the rank of a CSM below the values deduced
by general arguments.

We shall conclude our discussion of the universal CSM by
establishing its relation to $\sf \Omega$ of Theorem 1.  To that
end, let ${t}^{\alpha}_{\mu}$, $\alpha =1,2, \ldots, A$, $\mu
=1,2,\ldots,N$, be the set of isovectors used in the construction of
$\sf \Omega$.  Then the desired relation is given by ${\sf
\Omega}_{\alpha \beta} {=}{\sum}_{\mu,
\nu=1}^{N}{t}^{\alpha}_{\mu}{\sf U}_{\mu \nu}{t}^{\beta}_{\nu}$.  To
establish this result, first note that
${\sum}_{\mu=1}^{N}{t}^{\alpha}_{\mu}{x}_{i,\mu}$ by definition
equals ${x}_{i;\alpha}$, a generalized compensated derivative of a
decision function, and that according to Lemma 1 any such derivative
conforms to the constraints.  But then we must have
${\sum}_{i=1}^{M}{x}_{i;\alpha}{\sf Q}_{ij}=0$, $\alpha
=1,2,\ldots,A$, $j=1,2,\ldots,M$, since $\sf Q$ is the projection
matrix onto the normal hyperplane and annihilates any conforming
vector upon contraction.  This is a crucial property in the
derivation, and one that clearly shows how the application of
generalized compensated derivatives in effect obviates the use of
projection matrices.  To complete the derivation, we substitute the
defining expression for $\sf U$ given in Theorem 3 in the expression
${\sum}_{\mu, \nu=1}^{N}{t}^{\alpha}_{\mu}{\sf U}_{\mu
\nu}{t}^{\beta}_{\nu}$, and observe that the contributions involving
the projection operator $\sf Q$, including the second term within
the square brackets of Eq.~(\ref{023}), vanish by virtue of the
annihilation property just remarked upon.  Taking account of the
symmetry of either CSM, we find that the remaining contributions are
exactly equal to ${\sf \Omega}_{\alpha \beta}$, as claimed.

The relation between ${\sf \Omega}$ and ${\sf U}$ stated above
provides further insight into their rank properties.  Recalling the
theorems pertaining to the rank of a matrix quoted in \S IIB, we
conclude that the rank of ${\sf \Omega}$ is no larger than the
smaller of the ranks of ${\sf U}$ and ${t}^{\alpha}_{\mu}$, with the
latter viewed as an $A \times N$ matrix.  Since the latter rank is
never larger than $A$, we arrive at the result of Theorem 2-(iii)
that ${\rho}^{\sf \Omega} \leq \min (M- K,N,A)=\min (M-K,A)$, since
$A$ is never larger than $N$.  But this derivation also demonstrates
that in those cases where $N \leq M-K$ and $A < N$, the rank of $\sf
U$ will in general be larger than that of $\sf \Omega$, simply
because ${t}^{\alpha}_{\mu}$ is at most of rank $A$.  The latter is
in turn caused by the fact that the process of compensation
prohibits differentiation with respect to the \textit{normal}
directions in parameter space, while the present construction does
not.  As we have already remarked, and will be amply demonstrated in
the applications, one can always compensate for this shortfall by
introducing auxiliary parameters into the objective or constraint
functions.

As remarked earlier, the construction and application of the universal CSM
is not as straightforward as in the case of $\sf \Omega$.  For example, a
comparison of the two for the Slutsky-Hicks problem clearly shows that
$\sf U$ is less intuitive and more laborious to construct and, moreover,
that it emerges in a rather unfamiliar form which requires further
manipulation to appear as the standard result.

\section{APPLICATIONS}

The primary purpose of this section is to illustrate the workings of
the method of generalized compensated derivatives.  A second,
parallel objective is to present and discuss certain novel results
that naturally emerge from the applications of this method.  We have
already introduced a number of model problems in connection with the
construction of GCD's, and will be referring to these in the
following sections. Furthermore, we will not repeat the regularity
conditions stipulated in \S IIB for the objective and constraint
functions, assuming instead that they are appropriately fulfilled in
each application.  Throughout this section, we shall summarize the
main results obtained for the analyzed models in italicized and
numbered statements, each referred to as a \textit{property}.
Furthermore, we shall distinguish results which we believe to be new
by indenting them as new paragraphs.

\subsection{Profit Maximization}

Consider the basic profit maximization model introduced in \S IIA,
where one seeks to maximize the objective function $f(\mathbf{x},
\mathbf{a}) \stackrel{\rm{{def}}}{=} pF(\mathbf{x})-\mathbf{x} \cdot
\mathbf{w}$ by an appropriate choice of input factors $\bf x$. Since
this is an unconstrained problem, we can directly apply
Eq.~(\ref{19}) of Corollary A to it.  Taking the set of parameters
$\mathbf{a}$ to equal the single parameter $p$, we find (i)
$\partial F({x}(\mathbf{w},p))/
\partial p \geq 0$, while choosing $\mathbf{w}$ as the parameter set,
we obtain (ii) the matrix ${\sf W}_{\mu \nu}
(\mathbf{w},p)\stackrel{\rm{{def}}}{=}\partial
{x}_{\mu}(\mathbf{w},p) /
\partial {w}_{\nu}$ is negative semidefinite.  Since
semidefiniteness entails symmetry, the familiar reciprocity
relations among the demand functions are already implicit in
Property (ii).  To establish the homogeneity property of the demand
functions, on the other hand, we appeal to Theorem 2.  Here we
choose $\mathbf{X}(\mathbf{x})\stackrel{\rm{{def}}}{=}0$ and
$\mathbf{A}(\mathbf{a})\stackrel{\rm{{def}}}{=}\mathbf{a}$, where
the set of parameters $\mathbf{a}$ is now enlarged to
$(\mathbf{w},p)$. It can then be readily verified that the objective
function $f$ satisfies the prerequisites of that theorem (assuming
definiteness for the Hessian matrix associated with $F$).  The
conclusion, given in Eq.~(\ref{11}), is (iii) ${\sum}_{\mu=1}^{M}
{w}_{\mu} \partial \mathbf{x}(\mathbf{w},p) /
\partial {w}_{\mu} + p \partial \mathbf{x}(\mathbf{w},p) / \partial p =0$,
which, by way of Euler's theorem, asserts that the demand functions
(hence also the supply function) are homogeneous of degree zero in
the prices.  These are of course the standard properties of the
basic profit maximization model.

To provide further illustrations of the results of \S II, we now
proceed to analyze the above model by explicitly employing two
different sets of GCD's.  For the first, we take the set of
parameters to be $\mathbf{a}\stackrel{\rm{{def}}}{=}(\mathbf{w},p)$,
$M+1$ in number.  Then, starting with the $M$ partial derivatives
$\partial / \partial {w}_{\alpha}$, we construct a set of GCD's
equipped with the null property with respect to $f$ by the method of
one-term compensation, using $p$ as the compensating parameter.  The
result is $D_{\alpha}(\mathbf{x},\mathbf{a})=F(\mathbf{x}) \partial
/
\partial {w}_{\alpha}+{x}_{\alpha} \partial / \partial p$,
$\alpha=1,\ldots,M$. Substituting this set in Theorem 1, we find the
result that the matrix
\[ {\sum}_{\sigma=1}^{M}[{\delta}_{\mu \sigma}-{x}_{\mu}(\mathbf{w},p)
{w}_{\sigma}/pF(\mathbf{x}(\mathbf{w},p))] \{
{x}_{\sigma,\nu}(\mathbf{w},p) + [{x}_{\nu}(\mathbf{w},p)
{x}_{\sigma,M+1}(\mathbf{w},p) ]/ F(\mathbf{x}(\mathbf{w},p)) \} \]
is negative semidefinite.  Upon denoting the combination
${x}_{\mu}(\mathbf{w},p)/F(\mathbf{x}(\mathbf{w},p))$ by
${\zeta}_{\mu}(\mathbf{w},p)$, we can rearrange this result as

\textbf{Property (iv)} \textit{The matrix
\begin{equation}
{\sf Z}_{\mu \nu}(\mathbf{w},p)\stackrel{\rm{{def}}}{=}\partial
{\zeta}_{\mu}(\mathbf{w},p) / \partial {w}_{\nu}+
{\zeta}_{\nu}(\mathbf{w},p) \partial {\zeta}_{\mu}(\mathbf{w},p) /
\partial p \label{23}
\end{equation}
is negative semidefinite.}

This last result is formally identical to one derived by Paris,
Caputo, and Holloway (1993, Eq.~11) in their treatment of the
long-run competitive equilibrium model of the firm under conditions
of free entry and exit.  These authors set up an optimization
problem modeling the behavior of a representative firm under
equilibrium conditions, including a restriction to zero profits.  It
turns out that the resulting model is equivalent to profit
maximization restricted to zero profits, hence the formal identity
of their result to Property (iv).  Moreover, their derivation
continues to be valid at all (positive) profit levels, so that the
two results are essentially identical, albeit with differing
interpretations.  These authors then relate their result to that of
the minimum average cost model considered in earlier literature,
thereby pointing to a second method of deriving Eq.~(\ref{23}) at
zero profits. At any rate, the new insight gained here is the
realization that the comparative statics results conveyed by
Eq.~(\ref{23}) are those of any profit maximizing firm at any profit
level, and are therefore not specific characteristics of firms under
long-run equilibrium conditions.  We will return to $\sf Z$ below
and point out further properties of this matrix in connection with
Theorem 2.

At this point we turn to the second set of GCD's, slightly different
from the first set we just employed to derive Property (iv).  This
is the set we developed for the profit maximization problem in \S
IIA, where we enlarged the parameter set to include an auxiliary
parameter.  In effect, we will now seek to maximize the modified
objective function ${\tilde f}(\mathbf{x}, \mathbf{a})
\stackrel{\rm{{def}}}{=} s[pF(\mathbf{x})-\mathbf{x} \cdot
\mathbf{w}]$ with the understanding that the auxiliary parameter
$s>0$ must be set equal to unity to recover the original problem.
There are $M$ decision variables as before, but the parameter set
has been enlarged to
$\mathbf{a}\stackrel{\rm{{def}}}{=}(\mathbf{w},p,s)$, $M+2$ in
number. We will employ the set of $M+1$ GCD's developed for this
problem in \S IIA, and use these in Eq.~(\ref{7}), with $K=0$
corresponding to the absence of constraints. As a preliminary step,
let us carry out the differentiations required in Eq.~(\ref{7}),
then set $s$ equal to unity.  This gives us ${x}_{i;
\beta}(\mathbf{a})=[{x}_{i, \nu}(\mathbf{a}), {x}_{i,
M+1}(\mathbf{a})]$, and ${f}_{,i; \alpha}
(\mathbf{x}(\mathbf{a}),\mathbf{a})=[-{\delta}_{i
\mu},{F}_{,i}(\mathbf{x}(\mathbf{a}))]$, where $\alpha,
\beta=1,\ldots,M+1$, $i, \mu, \nu =1,\ldots,M$, and $\alpha=\mu,
\beta=\nu$ for $\mu, \nu=1,\ldots,M$. The CSM can now be assembled
according to Eq.~(\ref{7}).  It is actually more convenient at this
point to state the semidefiniteness properties of the resulting CSM
in two parts.  The first part conveys the symmetry of the CSM, and
it consists of the familiar reciprocity relations
${x}_{\mu,\nu}(\mathbf{w},p)={x}_{\nu,\mu}(\mathbf{w},p)$, as well
as (v) ${x}_{\mu,M+1} (\mathbf{w},p)=-\partial
F(\mathbf{x}(\mathbf{w},p)) /
\partial {w}_{\mu}$.  The second part conveys the specific
semidefiniteness property of the CSM, and it is most conveniently
stated as the nonnegativity of the quadratic form ${\sum}_{\alpha,
\beta=1}^{M+1} {b}_{\alpha}{\sf{\Omega}}_{\alpha \beta}(\mathbf{a})
{b}_{\beta}$.  Here $\mathbf{b}$ stands for the $(M+1)$-vector
$(\mathbf{q},q)$, where $\bf q$ and $q$ respectively represent an
arbitrary $M$-vector and a scalar.  Upon some rearrangement using
the reciprocity relations of $\sf W$, we can recast this result as
the inequality
\begin{equation}
{\sum}_{\mu,\nu=1}^{M}({q}_{\mu}-q{w}_{\mu}/p){x}_{\mu,\nu}(\mathbf{w},p)({q}_{\nu}-q{w}_{\nu}/p)
\leq 0. \label{24}
\end{equation}

We now proceed to consider certain consequences of Property (v) and
inequality (\ref{24}).  Let us start transforming (v) by writing
$\partial F(\mathbf{x}(\mathbf{w},p)) / \partial {w}_{\mu} =
{\sum}_{i=1}^{M}{F}_{,i}(\mathbf{x}(\mathbf{w},p)){x}_{i,
\mu}(\mathbf{w},p)$ and then eliminating
${F}_{,i}(\mathbf{x}(\mathbf{w},p))$ in favor of ${w}_{i}/p$ by
means of the first-order conditions.  An application of the
reciprocity relations to the emergent expression then yields the
homogeneity result found under (iii) by way of Euler's theorem. Here
we have derived the homogeneity property as a consequence of the
semidefiniteness of an enlarged CSM.

Next we turn to inequality (\ref{24}), and let
$q={\sum}_{\sigma=1}^{M}{q}_{\sigma}{l}_{\sigma}$, where $\bf l$ is any
real vector of dimension $M$.  This leads to the result that the matrix
\begin{equation}
 {\sum}_{\sigma, \rho =1}^{M}({\delta}_{\mu \sigma}-
{l}_{\mu}{w}_{\sigma}/p)(\partial {x}_{\sigma}(\mathbf{w},p) /
\partial {w}_{\rho}) ({\delta}_{\rho \nu}-{w}_{\rho }{l}_{\nu}/p)
\label{25}
\end{equation}
is negative semidefinite for any choice of $\bf l$.  This is a
powerful result; for example, it can be used to derive Properties
(ii) and (iv) by simply choosing $\bf l$ respectively equal to zero
and $\mathbf{x}(\mathbf{w},p)/F(\mathbf{x}(\mathbf{w},p))$ in
(\ref{25}). Furthermore, it can also be used to derive a sharpened
version of (ii) by choosing $\bf l$ equal to $[\partial
\mathbf{x}(\mathbf{w},p) /
\partial p] /[\partial F(\mathbf{x}(\mathbf{w},p)) / \partial p]$.  Upon
some rearrangement using the reciprocity relations and Euler's
theorem, this choice transforms the matrix in (\ref{25}) into a
useful form as summarized below.

\textbf{Property (vi)} \textit{the matrix
\begin{equation}
{\sf W}^{\star}_{\mu \nu}(\mathbf{w},p)\stackrel{\rm{{def}}}{=}{\sf
W}_{\mu \nu}(\mathbf{w},p) +{(\partial {x}_{\mu}(\mathbf{w},p) /
\partial p)( \partial {x}_{\nu}(\mathbf{w},p) / \partial p) \over
\partial F(\mathbf{x}(\mathbf{w},p)) / \partial p } \label{26}
\end{equation}
is negative semidefinite.}

This result, although not present in the standard treatments of the
profit maximization model in the literature, can nevertheless be
easily derived by standard methods, e.g., by exploiting the
convexity of the profit function with respect to all prices and
applying the envelope relations.  What turns Eq.~(\ref{26}) into a
powerful result is the realization that the second matrix on the
right-hand side is positive semidefinite since its denominator
$\partial F(\mathbf{x}(\mathbf{w},p)) / \partial p$ is nonnegative
by Property (i) and its numerator is of the generically positive
semidefinite complexion ${A}_{\mu}{A}_{\nu}$.  This fact implies
that the semidefiniteness of ${\sf W}^{\star}$ is a stronger
statement than that of ${\sf W}$.  In particular, the nonpositivity
of the diagonal elements of ${\sf W}^{\star}$ yields a sharpened
version of the standard result on the slopes of the demand functions
with respect to their own prices.  Since this is a rather
significant result and has obvious implications for issues related
to empirical testing, we will consider it in some detail.

Let us start by restating the condition ${\sf W}_{\mu \mu}^{\star}\leq 0$
in a different but equivalent form as

\textbf{Property (vii)} \textit{The slope of the demand functions
with respect to their own prices are bounded above according to}
\begin{equation}
\partial {x}_{\mu}(\mathbf{w},p) / \partial {w}_{\mu} \leq -{(\partial
{x}_{\mu}(\mathbf{w},p) / \partial p)}^{2}/[\partial
F(\mathbf{x}(\mathbf{w},p)) /
\partial p].  \label{27}
\end{equation}
Clearly, this inequality is in general stronger than its standard
counterpart, $\partial {x}_{\mu}(\mathbf{w},p) / \partial {w}_{\mu}
\leq 0 $, and more significantly, it sets limits on the magnitude of
the slopes in question.  To see a numerical example of this
comparison, let us consider the Cobb-Douglas production function,
given by ${F}^{CD}(\mathbf{x})\stackrel{\rm{{def}}}{=}{F}_{0} \,
{\prod}_{i=1}^{M} {x}_{i}^{{\gamma}_{i}}$, where
$0<{{\gamma}_{i}}<1$,
${\sum}_{i=1}^{M}{\gamma}_{i}\stackrel{\rm{{def}}}{=}\gamma<1$,
$i=1,\ldots,M$, and ${F}_{0}$ is a fixed positive number.  A
straightforward, if tedious, calculation gives the results
\begin{eqnarray}
{\partial {x}^{CD}_{\mu}(\mathbf{w},p) \over \partial {w}_{\nu}}&&=
-{{x}^{CD}_{\mu}(\mathbf{w},p) \over {w}_{\nu}}({\delta}_{\mu
\nu}+{{\gamma}_{\nu} \over 1-\gamma}), \nonumber \\
{\partial {x}^{CD}_{\mu}(\mathbf{w},p) \over \partial
p}&&={{x}^{CD}_{\mu}(\mathbf{w},p) \over p}{1 \over 1-\gamma}, \nonumber \\
{\partial {F}^{CD}(\mathbf{w},p) \over \partial
p}&&={{F}^{CD}(\mathbf{w},p) \over p}{\gamma \over 1-\gamma}
\nonumber
\end{eqnarray}
for the Cobb-Douglas production function.  It is worth noting that
the remarkable simplicity of these results is essentially a
consequence of the special scaling properties of the model.  These
scaling properties are most clearly evident in the set of relations
$\partial \log F^{CD}(\mathbf{x}) / \partial \log {x}_{i} =
{\gamma}_{i}$.  Using these results, we can calculate the bound
given on the right-hand side of (\ref{27}) for the Cobb-Douglas
case, and compare it to zero (which is the bound given by the
standard result) as well as to the exact result (i.e., $\partial
{x}^{CD}_{\mu}(\mathbf{w},p) / \partial {w}_{\nu}$) calculated for
the model. This comparison is most usefully stated in terms of the
elasticity ${\varepsilon}_{\mu \mu}\stackrel{\rm{{def}}}{=}\partial
\log[ {x}_{\mu}(\mathbf{w},p)] / \partial \log ({w}_{\mu})$ rather
than the dimensional slope given in (\ref{27}).  A calculation using
the above results shows that the demand elasticities for the
Cobb-Douglas case, ${\varepsilon}^{CD}_{\mu \mu}$, are bounded above
by $-{\gamma}_{\mu} / \gamma(1-\gamma)$ according to (\ref{27}).
This should be compared to the (weaker) bound of zero given by the
standard result, and to the exact value ${\varepsilon}^{CD}_{\mu
\mu} =-[1+{\gamma}_{\mu}/(1-\gamma)]$.  For a two-factor model with
${\gamma}_{1}={\gamma}_{2}=1/3$, for example, these three numbers
are $-3/2$, $0$, and $-2$, respectively.  Clearly, the restrictive
conditions given in Eqs.\ (\ref{26}) and (\ref{27}) convey more
stringent consequences of the profit maximization hypothesis than
those given by the standard result, thus providing a correspondingly
more rigorous means of testing that hypothesis.

Not surprisingly, there exists a sharpened bound on the elasticity
of the supply function as well.  The standard result, stated above
as Property (i), guarantees that $\partial F({x}(\mathbf{w},p))/
\partial p \geq 0$.  On the other hand, a procedure similar to the
one leading to (\ref{27}) yields

\textbf{Property (viii)} \textit{The slope of the supply function is
bounded below according to the following inequality: }
\begin{equation}
\partial F(\mathbf{x}(\mathbf{w},p)) / \partial p \geq -{\sum}_{\mu,\nu=1}^{M}
({w}_{\mu}/p) ({w}_{\nu}/p) \partial {x}_{\mu}(\mathbf{w},p) /
\partial {w}_{\nu} .  \label{28}
\end{equation}
Since the right-hand side of (\ref{28}) is positive semidefinite, we
have in this inequality a sharpened version of the standard result.
Using the Cobb-Douglas case as an example again, we find that the
elasticity of the supply function, $\sigma \stackrel{\rm{{def}}}{=}
\partial \log [F(\mathbf{x}(\mathbf{w},p))] / \partial \log(p)$, is
bounded below by $\gamma / (1- \gamma)$.  This bound turns out to be
the same as the exact value ${\sigma}^{CD}$, and greater than the
bound of zero given by the standard result.

For reference purposes, we will outline here generalized versions of
(\ref{27}) and (\ref{28}) appropriate to the multi-output profit
maximization model.  This model is defined by $f(\mathbf{x},
\mathbf{a}) \stackrel{\rm{{def}}}{=} \mathbf{p}\cdot
\mathbf{F}(\mathbf{x})-\mathbf{x} \cdot \mathbf{w}$, where the
production function and the output price have been promoted to the
status of $G$-dimensional vectors.  The basic CSM, constructed for
the set of parameters $\mathbf{a}
\stackrel{\rm{{def}}}{=}(\mathbf{w}, \mathbf{p})$, is most
conveniently described as a $2 \times 2$ block matrix
\begin{equation}
 \sf{T} \stackrel{\rm{{def}}}{=}
 \begin{bmatrix} -\sf{W} & -\sf{M} \\
 \sf{Q} & \sf{P}
 \end{bmatrix},
\end{equation}
where $\sf W$ is the matrix of the partial derivatives of the input
factors with respect to the input prices introduced above, $\sf P$
is its output counterpart ${\sf
P}_{rs}(\mathbf{w},\mathbf{p})\stackrel{\rm{{def}}}{=}\partial
{F}_{r}(\mathbf{x}(\mathbf{w},\mathbf{p})) / \partial {p}_{s}$,
$r,s=1,\ldots,G$, ${\sf M }_{\mu r}(\mathbf{w},\mathbf{p})
\stackrel{\rm{{def}}}{=}\partial {x}_{\mu}(\mathbf{w},\mathbf{p}) /
\partial {p}_{r}$, and ${\sf Q}_{r \mu }(\mathbf{w},\mathbf{p})
\stackrel{\rm{{def}}}{=}\partial
{F}_{r}(\mathbf{x}(\mathbf{w},\mathbf{p})) /
\partial {w}_{\mu}$.  The last two are typically rectangular matrices,
respectively $M \times G$ and $G \times M$ in size.  The comparative
statics properties of this model are summarized in the positive
semidefiniteness of the block matrix $\sf{T}$.  This statement in turn
implies the positive semidefiniteness of the blocks $-\sf W$ and $\sf P$,
as well as the equality $-{\sf M}={\sf Q}^{\dag}$ required by the
symmetry of $\sf{T}$.  These are the standard comparative statics
results of the model, and in effect generalize the basic results of the
single-output model to the present case.

One can improve upon these standard results by treating $\sf{T}$ as
a whole and seeking to optimize the bounds on the slopes of the
demand and supply functions as was done above.  Remarkably, this can
be accomplished in the following elegant manner.  Consider a real
$(M+G)$-vector $\mathbf{z}=(\mathbf{u},\mathbf{v})$, $\bf u$ and
$\bf v$ ``block vectors'' of dimension $M$ and $G$ respectively, and
form the quadratic form $\mathbf{z}^{\dag} \sf{T} \mathbf{z}$.  This
quantity is guaranteed to be nonnegative by the semidefiniteness of
$\sf{T}$, a property which implies the inequality $\mathbf{u}^{\dag}
\sf{W} \mathbf{u} \leq -\mathbf{u}^{\dag} \sf{M} \mathbf{v}
-\mathbf{v}^{\dag} {\sf{M}}^{\dag}\mathbf{u} + \mathbf{v}^{\dag}
\sf{P} \mathbf{v}$. At this point, we seek to minimize the
right-hand side of this inequality with $\bf v$ as choice variables
and $\bf u$ as parameters.  This requirement amounts to looking for
the least upper bound to the eigenvalues of $\sf{W}$.  The solution
to this optimization problem is found to be
$\mathbf{v}^{\star}={\sf{P}}^{- 1}{\sf{M}}^{\dag} \mathbf{u}$,
which, upon substitution in the inequality found above, leads to the
result
\[ \mathbf{u}^{\dag} \sf{W} \mathbf{u} \leq -\mathbf{u}^{\dag} \sf{M} {\sf{P}}^{-
1}{\sf{M}}^{\dag} \mathbf{u}. \] Recalling at this point that $\bf
u$ is an arbitrary vector, we can restate this inequality as the
negative semidefiniteness of the matrix ${\sf W}+{\sf M}{\sf
P}^{-1}{\sf M}^{\dag}$.  An entirely analogous procedure yields a
corresponding result for the output block $\sf{P}$.  These two
results, which are the generalizations of Properties (vii) and
(viii) of the single-output case, are summarized as

\textbf{Property (ix)} \textit{The pair of matrices
\begin{equation}
{\sf W}^{\star}\stackrel{\rm{{def}}}{=}{\sf W}+{\sf M}{\sf
P}^{-1}{\sf M}^{\dag}, \:\:{\sf
P}^{\star}\stackrel{\rm{{def}}}{=}{\sf P}+{\sf M}^{\dag}{\sf
W}^{-1}{\sf M}  \label{29}
\end{equation}
are negative semidefinite and positive semidefinite respectively.}

Because ${\sf
M}{\sf P}^{-1}{\sf M}^{\dag}$ and ${\sf M}^{\dag}{\sf W}^{-1}{\sf M}$ are
positive semidefinite and negative semidefinite matrices respectively, the
new pair
${\sf W}^{\star}$ and ${\sf P}^{\star}$ represent sharpened versions of
${\sf W}$ and ${\sf P}$ respectively, in the same sense and for
essentially the same reasons as were explained in the single-output case
above.  Thus Eqs.\ (\ref{29}) provide the desired generalizations of
(\ref{26}) and (\ref{28}) to the multi-output case.  It might be mentioned
here that the new pair are actually optimal in the sense that they cannot
in general be sharpened any further.  It may also be added that the
occurrence of inverse matrices in Eqs.\ (\ref{29}) does not invalidate the
intended results in case $\sf W$ or $\sf P$ turns out to be singular,
although further analysis may be necessary to reestablish the appropriate
results.

It is appropriate at this point to pause and recall Theorem 2-(iv)
and our frequent references to the arbitrariness in the choice of
GCD's and the form of the resulting CSM's vis-\`{a}-vis the fact
that all such forms have an essentially unique source, namely the
local concavity of the objective function.  As an illustration of
this, consider the structure in (\ref{25}) which consists of a
family of CSM's for the (single output) profit maximization model.
Now this family is precisely of the form given in Theorem 2-(iv)
with ${\sf T}_{\mu \sigma}$ identified with ${\sf{\Delta}}_{\mu
\sigma}\stackrel{\rm{{def}}}{=}{\sf{\delta}}_{\mu
\sigma}-{l}_{\mu}{w}_{\sigma}/p$ and $\sf {\Omega}$ with ${\sf W}$,
where, it will be recalled, $\sf W$ stands for the negative
semidefinite matrix $\partial {x}_{\sigma}(\mathbf{w},p) / \partial
{w}_{\rho}$, the standard CSM of the model.  At this point it is
useful to know the spectrum of $\sf {\sf{\Delta}}$, which can be
found straightforwardly.  It consists of an $(M-1)$-fold degenerate
eigenvalue equal to unity (the corresponding eigenvectors being any
linearly independent set of $M-1$ vectors orthogonal to $\bf w$) and
a simple eigenvalue equal to $1-\mathbf{w}\cdot \mathbf{l}/p$ (the
corresponding eigenvector being $\bf l$).  Clearly, the determinant
of $\sf{\Delta}$ equals this last eigenvalue, so that as long as
$p-\mathbf{w}\cdot \mathbf{l} \neq 0$, the matrix $\sf {\Delta}$ is
nonsingular and the above family is congruent to $\sf W$, as in
Eq.~(\ref{10}).  However, for $p-\mathbf{w}\cdot \mathbf{l} =0$,
$\sf{\Delta}$ loses full rank and becomes singular, as in
Eq.~(\ref{11}), causing a breakdown of the congruence between $\sf
W$ and the above family of CSM's. Remarkably, this is exactly what
happens to $\sf Z$ at zero profits (although it has nothing to do
with profits as such; see the following paragraph).  To see this,
recall that $\sf Z$ is a member of the above family with
$\mathbf{l}=\mathbf{x}(\mathbf{w},p)/F(\mathbf{x}(\mathbf{w},p))$.
Thus the singular point for $\sf Z$ corresponds to
$pF(\mathbf{w},p)-\mathbf{w}\cdot \mathbf{x}(\mathbf{w},p)=0$, the
point of vanishing profits.  Of course this does not invalidate $\sf
Z$ as a CSM, as was pointed out in \S IIB, but it does reduce its
rank and the information it contains, relative to $\sf W$, at the
zero-profits point.

Why does the rank reduction just discussed occur at the point of
vanishing profits?  After all, the zero-profits point can be moved
by the addition of a fixed constant to the objective function and
therefore seems to be no different than any other profit level.  The
answer is revealed by an application of Theorem 2-(v).  Recall that
the latter allows a transformation of decision variables (as well as
the parameters), subject to the condition that the associated
Jacobian determinant not vanish.  Let us therefore consider a
transformation from $\mathbf{x}$ to $\tilde{\mathbf{x}}
\stackrel{\rm{def}}{=}\mathbf{x}/F(\mathbf{x})$.  Now consider the
matrix $\partial {\tilde{x}_{i}} / \partial {x}_{j}$, which is found
to equal
$[{\delta}_{ij}-{x}_{i}{w}_{j}/pF(\mathbf{x})]/F(\mathbf{x})$, where
first-order conditions have been used.  But this matrix is, save for
the scalar factor $1/F(\mathbf{x})$, precisely the one that was
found in the preceding paragraph to become singular at the point of
vanishing profits, or more explicitly, at the point where
$pF(\mathbf{x}(\mathbf{w},p))-\mathbf{x}(\mathbf{w},p) \cdot
\mathbf{w}$ vanishes. Therefore what goes wrong at this point,
causing the rank reduction noted above, is the vanishing of the the
Jacobian determinant for the transformation from ${x}_{i}$ to
${\zeta}_{i}={\tilde{x}}_{i}$. Therefore, as noted above, this
phenomenon has nothing to do with profits as such, a fact that can
be illustrated by considering a profit maximizing firm whose
operations entail a fixed cost, i.e., a firm with the objective
function $f(\mathbf{x}, \mathbf{a}) \stackrel{\rm{{def}}}{=}
pF(\mathbf{x})-\mathbf{x} \cdot \mathbf{w}-C$, $C>0$. Here, the
catastrophe would occur not at zero profits, but at the point where
$pF(\mathbf{x}(\mathbf{w},p))- \mathbf{x}(\mathbf{w},p) \cdot
\mathbf{w}$ vanishes, which corresponds to a loss equal to $C$.

The circumstance that a variety of CSM's is implied by the
underlying structure of a model--the local concavity of the
production function being the relevant property here--is a most
desirable one, especially from the standpoint of comparison with
empirical data.  The example given in Eq.~(\ref{26}) et seq. is a
case in point, since it demonstrates how the available variety of
CSM's can be explored for the purpose of finding the most incisive
implications of a given model for confrontation with measured data.

The last case to be considered in this section is a constrained variant of
the multi-output profit maximization model briefly considered above.  This
model, which has received some attention recently in connection with the
behavior of U.S. agricultural producers, is defined by
\begin{equation}
{\max \;}_\mathbf{x} {\;} s[\mathbf{p} \cdot
\mathbf{F}(\mathbf{x})-\mathbf{x} \cdot \mathbf{w}] {\:} {\;} s.t.
{\:} {\:} C- \mathbf{x} \cdot \mathbf{w}=0, \label{30}
\end{equation}
where, as before, the auxiliary parameter $s>0$ will eventually be
set equal to unity.  We shall refer to this problem as the
\textit{cost-constrained} profit maximization model; clearly, its
objective is to maximize profits subject to a prescribed total cost.
As in the unconstrained version of this model, there are $M$ input
factors and $G$ output products.  We shall take the parameter set to
be the $(M+G+2)$-dimensional vector
$\mathbf{a}\stackrel{\rm{{def}}}{=} (\mathbf{w},\mathbf{p},C,s)$.

Remarkably, Eq.~(\ref{30}) can be reinterpreted as a utility
maximization problem, a fact that appears to have escaped notice in
the literature.  To bring out this analogy, let us first note that
it is simpler to deal with the equivalent objective function
$s\mathbf{p} \cdot \mathbf{F}(\mathbf{x})$ by using the cost
constraint to replace the second term $s\mathbf{x} \cdot \mathbf{w}$
with $sC$, and then dropping it altogether since it is free of
decision variables. The formal analogy to the basic utility
maximization problem is now evident in the transformed optimization
problem ${\max \;}_\mathbf{x} {\;} s\mathbf{p} \cdot
\mathbf{F}(\mathbf{x}) {\:} {\;} s.t. {\:} {\:} C- \mathbf{x} \cdot
\mathbf{w}=0$, with $s\mathbf{p} \cdot \mathbf{F}(\mathbf{x})$
playing the role of the utility function and $C$ representing the
income term. In spite of this formal analogy, we will continue to
analyze this problem since it will lead to certain novel
observations in addition to illustrating the results of \S II.

It is convenient to start the analysis of this problem by exploiting
its symmetries in accordance with Theorem 2-(ii).  We do this by
choosing $J(\mathbf{x},\mathbf{a})\stackrel{\rm{{def}}}{=}
{b}_{1}{\sum}_{r=1}^{G} {p}_{r}(\partial / \partial
{p}_{r})+{b}_{2}( {\sum}_{r=1}^{M} {w}_{r}\partial / \partial
{w}_{r}+C \partial / \partial C)$, where ${b}_{1}$ and ${b}_{2}$ are
an arbitrary pair of real numbers.  This equation in effect defines
the pair of functions $\bf X$ and $\bf A$ introduced in connection
with the treatment of symmetries in \S IIB.  With this choice, the
conditions of Theorem 2-(ii) are satisfied (assuming a definite
Hessian for $\bf F$ at the maximum point), and we have the
conclusion that $J(\mathbf{x},\mathbf{a})
\mathbf{x}(\mathbf{w},\mathbf{p}, C)=0$. Choosing first ${b}_{1}
\neq 0, {b}_{2}=0$, then ${b}_{1}=0, {b}_{2} \neq 0$, we find from
the invariance result (via Euler's Theorem) two separate homogeneity
results as follows:

\textbf{Property (x)} \textit{The demand and supply functions are
homogeneous of degree zero in the \textit{output prices}, and,
independently, in the set $(\mathbf{w},C)$ as well.}

Note that these separate scale invariance properties constitute
a stronger condition than a joint one involving all the parameters.  In
particular, for $G=1$, the above homogeneity property in the (single)
output price implies a rather surprising result which we state as

\textbf{Property (xi)} \textit{In case of a single output, the
demand and supply functions are independent of the output price
$p$.}

In other words, the production decisions of a
cost-constrained, profit-maximizing firm producing a single output are
insensitive to changes in the output price.

Continuing with the analysis of the model, we turn to the
construction of the GCD's.  This is conveniently done by the method
of one term compensation using $s$ and $C$ for compensation on the
objective and constraint functions respectively.  In this manner we
find $D_{\alpha}(\mathbf{x},\mathbf{a})
\stackrel{\rm{{def}}}{=}\partial / \partial {w}_{\alpha} +
{x}_{\alpha} \partial / \partial C $ for $\alpha =1, \ldots, M$, and
$D_{\alpha}(\mathbf{x},\mathbf{a}) \stackrel{\rm{{def}}}{=}\partial
/
\partial {p}_{\alpha} - [{F}_{\alpha} / \mathbf{p}\cdot\mathbf{F}
(\mathbf{x})]s \partial / \partial s$ for $\alpha =M+1, \ldots,
M+G$. Using these derivatives, we can proceed to the construction of
the CSM according to Eq.~(\ref{7}). After some straightforward
algebra including the use of the constraint, and upon setting $s=1$,
we find that the matrix
\begin{equation}
\begin{bmatrix}
- \lambda [{\sf W}_{\mu \nu}+{x}_{\nu}(\mathbf{a})
\partial {x}_{\mu}(\mathbf{a}) / \partial C] &\,\, -\lambda {\sf M} \\
[{\sf Q}_{r \nu}+{x}_{\nu} (\mathbf{a})
\partial {F}_{r}(\mathbf{x}(\mathbf{a}))
/ \partial C] &\,\, {\sf P} \\
\end{bmatrix}
\end{equation} is positive
semidefinite. Here the matrices $\sf{W},\sf{M},\sf{P},\sf{Q}$ are
the same as those introduced in connection with the multi-output
profit maximization model, with their orders specified according to
$\mu, \nu =1,2,\ldots,M$ and $r,s=1,2,\ldots,G$.  Furthermore,
$\lambda \stackrel{\rm{{def}}}{=}{w}_{\mu}^{-1}
\partial [\mathbf{p}\cdot\mathbf{F}(\mathbf{x}(\mathbf{a}))] / \partial {x}_{\mu} > 0$
is a Lagrange multiplier which may be interpreted as the
\textit{marginal profitability of expenditure}, i.e., the increase
in maximized profits due to a marginal increase in allowable
expenditures $C$.  The full CSM given above is of order $M+G$, but
it has a smaller rank.  To see this, we appeal to Theorem 2-(iii)
which states that the rank of the above CSM is at most equal to the
smaller of $M-1$ and $M+G$.  Since $G \geq 1$, we conclude that the
CSM has a rank no larger than $M-1$ while its order is $M+G$,
confirming the above assertion.

The comparative statics properties of this model are thus summarized in

\textbf{Property (xii)} \textit{The two $M \times M$ and $G \times
G$ matrices
\[ -[\partial {x}_{\mu}(\mathbf{w},\mathbf{p},C) / \partial {w}_{\nu}+{x}_{\nu}
(\mathbf{w},\mathbf{p},C) \partial
{x}_{\mu}(\mathbf{w},\mathbf{p},C) /
\partial C], {\;\;} \partial {F}_{r}(\mathbf{w},\mathbf{p},C) / \partial
{p}_{s} \] are positive semidefinite.  Furthermore, we have the
equality}
\[ \lambda(\mathbf{w},\mathbf{p},C) \partial {x}_{\mu}(\mathbf{w},\mathbf{p},C) /
\partial {p}_{s} =- [\partial {F}_{s}(\mathbf{w},\mathbf{p},C) / \partial
{w}_{\mu} +{x}_{\mu}(\mathbf{w},\mathbf{p},C) \partial
{F}_{s}(\mathbf{w},\mathbf{p},C) / \partial C]. \] The similarity of
this model to the basic utility maximization problem is now evident
in the first of these matrices, which is in fact the analog of the
Slutsky-Hicks substitution matrix.

This concludes the application of our methods to certain models of profit
maximization.

\subsection{Generalized Utility Maximization}

In this section we shall be concerned with constrained optimization
problems of the general form ${\max \;}_\mathbf{x} {\;}
f(\mathbf{x},\mathbf{b}) {\:} {\;} s.t. {\:} {\:}
{g}^{k}(\mathbf{x},\mathbf{b}^{\prime})=0$, $k=1,2,\ldots,K$, where
$\mathbf{b}$ and $\mathbf{b}^{\prime}$, the parameters of interest
for comparative statics information, are separated into two sets,
one occurring in the objective function and the other in the
constraints.  We shall refer to this category as \textit{generalized
utility maximization problems}, reserving the label ``utility
maximization problem'' for those cases where the parameters of
interest for comparative statics information appear only in the
constraint functions.  Thus the latter assume the familiar structure
\[ {\max \;}_\mathbf{x} {\;} U(\mathbf{x}) {\:} {\;} s.t. {\:} {\:} {g}^{k}(\mathbf{x},\mathbf{a})=0, k=1,2,\ldots,K. \]
We shall start our discussion with the best-known example of this class,
the Slutsky-Hicks problem, and continue with a sequence of generalizations
toward more general forms, as well as applications to illustrate these.

In \S IIA we constructed a set of GCD's for the budget constraint that
appears in the Slutsky-Hicks problem, the basic utility maximization
model.  This model is defined by
\begin{equation}
{\max \;}_\mathbf{x} {\;} U(\mathbf{x}) {\:} {\;} s.t. {\:} {\:}
m-\mathbf{p} \cdot \mathbf{x}=0,  \label{31}
\end{equation}
where $U$ is a (quasi-concave, strongly monotonic) utility function.
The parameter set is $\mathbf{a} \stackrel{\rm{{def}}}{=}
(\mathbf{p},m)$, $N=M+1$ in number.  Using the construction
mentioned above, we find the well known result that (i) the Slutsky
matrix
\begin{equation}
{\sf{\Sigma}}_{\mu \nu}\stackrel{\rm{{def}}}{=}\partial
{x}_{\mu}(\mathbf{p},m) / \partial {p}_{\nu}
+{x}_{\nu}(\mathbf{p},m)
\partial {x}_{\mu}(\mathbf{p},m) / \partial m \label{32}
\end{equation}
is negative semidefinite.  Moreover, as an illustration of Theorem
2-(ii), we showed in \S IIB that (ii) the demand functions are
homogeneous of degree zero in the parameters $(\mathbf{p},m)$.

Recall that in our discussion of the various degrees of arbitrariness in
the construction of CSM's in \S IIB and IIIA, we deferred a more detailed
consideration of the selection and number of decision variables and
parameters to the present section.  We emphasized in our earlier
discussions, particularly in connection with Theorem 2-(v), that while one
ordinarily chooses certain decision variables and parameters as the
``natural'' or ``relevant'' ones to use for a given problem, such a choice
is, at least mathematically, only one among an infinite family of possible
complexions.  Here, using the Slutsky-Hicks problem as an example, we will
implement a contraction in parameter space followed by one in decision
space to illustrate and amplify a number of points mentioned earlier.

We have already stated that redundant comparative statics descriptions
(for which the CSM is less than full rank) are common, as well as useful,
in
economic applications.  In such cases, the CSM may be contracted to lower
dimensions without any loss of comparative statics information (in the
sense that the original CSM can be reconstructed from the reduced one).
For the Slutsky-Hicks problem under discussion, Theorem 2-(iii) implies
that the Slutsky matrix $\sf{\Sigma}$ cannot have a rank larger than $M-
1$.  Therefore, it should be possible to convey the comparative statics
information contained in $\sf{\Sigma}$ by means of a reduced, $(M-1)
\times (M-1)$ matrix.  In the following, we will implement this reduction,
in part to see whether there is any advantage to doing so.

As a preliminary step to reducing ${\sf{\Sigma}}$, let us eliminate
the income effect terms in favor of the price effect terms by using
the homogeneity condition of Property (ii) above.  The result is
${\sf{\Sigma}}_{\mu \nu}={\sum}_{\gamma=1}^{M}[\partial
{x}_{\mu}(\mathbf{p},m) / \partial {p}_{\gamma}]({\delta}_{\gamma
\nu}- {p}_{\gamma}{x}_{\nu}/m)$.  In effect, we have factored the
Slutsky matrix into the product form
$\sf{\Sigma}\stackrel{\rm{{def}}}{=}\sf{M} \sf{\Upsilon}$, where
$\sf{M}$ and $\sf{\Upsilon}$ can be identified by reference to the
subscripted form of the equation.  Now one can easily verify that
the price vector $\bf p$ is a null vector for $\sf{\Upsilon}$, i.e.,
$\sf{\Upsilon}\mathbf{p}=0$, as a direct consequence of the budget
constraint.  One can also verify that $\sf{\Upsilon}$ is idempotent,
i.e., that ${\sf{\Upsilon}}^{2}=\sf{\Upsilon}$, so that
$\sf{\Sigma}\sf{\Upsilon}=\sf{\Sigma}$.  Furthermore, since
$\sf{\Sigma}$ is symmetric, one can reexpress it as
${\sf{\Sigma}}^{\dag}={\sf{\Sigma}}={\sf{\Upsilon}}^{\dag}
{\sf{M}}^{\dag}$, and use idempotency to modify the latter to read
$\sf{\Sigma}={\sf{\Upsilon}}^{\dag} {\sf{M}}^{\dag}
{\sf{\Upsilon}}$. Finally, by a trivial change of parameters from
$\mathbf{p}$ to $\tilde{\mathbf{p}} \stackrel{\rm{{def}}}{=}
\mathbf{p}/m$, and another trivial change from $\sf{{\Sigma}}$ to
${\tilde{\sf{\Sigma}}}\stackrel{\rm{{def}}}{=}m{\sf{\Sigma}}$, we
can eliminate any reference to the parameter $m$ in the resulting
CSM;
\begin{equation}
{\tilde{\sf{\Sigma}}}_{\mu \nu}(\tilde{\mathbf{p}})
\stackrel{\rm{{def}}}{=} {\sum}_{\tau,\rho =1}^{M}[{\delta}_{\mu
\tau}-{x}_{\mu} \tilde{\mathbf{p}}) {\tilde{p}}_{\tau}]\partial
{x}_{\rho }(\tilde{\mathbf{p}}) /
\partial {\tilde{p}}_{\tau}[{\delta}_{\rho \nu }-{\tilde{p}}_{\rho}
{x}_{\nu}(\tilde{\mathbf{p}}) ].  \label{33}
\end{equation}

Clearly, this matrix is negative semidefinite, and it is fully
equivalent to the original Slutsky matrix in information content by
virtue of the budget constraint $1-\mathbf{x} \cdot
\tilde{\mathbf{p}}=0$.  We now have a description in terms of a
reduced parameter set $\tilde{\mathbf{p}}$, with the homogeneity
information implicit in the fact that the decision functions depend
on prices and income through the ratio $\mathbf{p}/m$. Note also
that the singular nature of the Slutsky matrix, $\tilde{\sf{\Sigma}}
\cdot \tilde{\mathbf{p}}=\tilde{\mathbf{p}} \cdot \tilde
{{\sf{\Sigma}}}=0$, is manifest in the new description. These
mathematical niceties notwithstanding, it is obvious that
$\tilde{\sf{\Sigma}}$ is not amenable to a clear economic
interpretation; the precious intuition afforded by the original
Slutsky matrix has been disguised in the new representation.

The second step in the reduction relies on the fact that
${\tilde{\sf{\Sigma}}}$ is singular, and as mentioned above,
redundant. To see how singularity implies redundancy, note that the
row and column vectors of a square, singular matrix must be linearly
dependent, so that one or more of these vectors may be expressed in
terms of the others. Therefore, as stated earlier, the dependent
rows and columns may be dropped with no loss of generality.  In our
example, any $(M-1) \times (M- 1)$ submatrix of
${\tilde{\sf{\Sigma}}}$ would still convey the information contained
in the full matrix.  We can therefore drop, e.g., the last row and
column of ${\tilde{\sf{\Sigma}}}$, thus arriving at a CSM in terms
of a reduced parameter set and with a rank which is consistent with
Theorem 2-(iii).  In general, this is the most economical
description of the comparative statics of the Slutsky-Hicks problem.
But it is certainly not the most cogent statement of the economics
of that model. Indeed the redundant Slutsky matrix expressed as a
function of a redundant set of parameters, as in Eq.~(\ref{32}),
provides a much better intuitive understanding of the economics than
the upper-left block of ${\tilde{\sf{\Sigma}}}$.  There is obviously
no expository advantage in adopting the reduced description in this
case.  On the other hand, it is entirely conceivable that the form
given in Eq.~(\ref{33}), or the fully reduced description just
discussed, will provide more suitable alternatives for comparison
with measured data and empirical tests of the underlying model.

Having discussed the basic utility model in some detail, we now
proceed to a number of generalizations whereby the single linear
constraint of Eq.~(\ref{31}) is replaced with multiple constraints,
linear and nonlinear. As a preliminary to considering an arbitrary
number of linear constraint equations, it is convenient to first
analyze the case of two such equations. The model in question then
reads
\begin{equation}
{\max \;}_\mathbf{x} {\;} U(\mathbf{x}) {\:} {\;} s.t. {\:} {\:}
{m}^{1}-{\mathbf{p}^{1}} \cdot \mathbf{x}= 0 {\:}{\:} and {\:}{\:}
{m}^{2} -\mathbf{p}^{2} \cdot \mathbf{x} = 0, \label{34}
\end{equation}
where, in order to avoid trivialities, we will take the two price
vectors to be nonparallel.  The construction of the CSM for this
model parallels the treatment of the basic model except for the
doubling of the parameter set here to
$\mathbf{a}\stackrel{\rm{{def}}}{=}(\mathbf{p}^{1},{m}^{1},\mathbf{p}^{2},{m}^{2})$,
for a total of $N=2(M+1)$ parameters.  Note that a given parameter
only appears in one of the constraint equations. Consequently, the
GCD's, also double as many as before, have the same basic structure
as before.  They are given by ${D}_{\alpha}(\mathbf{x},
\mathbf{a})\stackrel{\rm{{def}}}{=}{\partial \over
\partial {p}^{1}_{\alpha}} +{x}_{\alpha} {\partial \over \partial
{m}^{1}}$ for $\alpha=1,2,\ldots,M$, and by
${D}_{\alpha}(\mathbf{x},
\mathbf{a})\stackrel{\rm{{def}}}{=}{\partial \over
\partial {p}^{2}_{\alpha -M}} +{x}_{\alpha -M} {\partial \over \partial
{m}^{2}}$ for $\alpha=M+1,M+2,\ldots,2M$.  It is useful to note at this
point that according to the rank inequality formula of Theorem 2-(iii) the
$2M \times 2M$ CSM for this problem will have a rank no larger than
$\min(M-2,2M)=M-2$.  We therefore anticipate a highly redundant CSM for
this problem.

Returning to the construction of the CSM, we find by a
straightforward application of Eq.~(\ref{7}) that it takes the
``double-Slutsky'' form
$$
\begin{bmatrix} {\lambda}_{1}{\sf{\Sigma}}^{1} & {\lambda}_{1}{\sf{\Sigma}}^{2}
\\ {\lambda}_{2}  {\sf{\Sigma}}^{1} & {\lambda}_{2}
{\sf{\Sigma}}^{2}
\end{bmatrix},
$$
where ${\sf{\Sigma}}^{1}_{\mu
\nu}\stackrel{\rm{{def}}}{=}\partial {x}_{\mu}(\mathbf{a}) /
\partial {p}^{1}_{\nu} +{x}_{\nu}(\mathbf{a}) \partial
{x}_{\mu}(\mathbf{a}) /
\partial {m}^{1}$, ${\sf{\Sigma}}^{2}_{\mu
\nu}\stackrel{\rm{{def}}}{=}\partial {x}_{\mu}(\mathbf{a}) /
\partial {p}^{2}_{\nu} +{x}_{\nu}(\mathbf{a}) \partial
{x}_{\mu}(\mathbf{a}) /
\partial {m}^{2}$, $\mu,\nu=1,2,\ldots,M$, and $ {\lambda}_{1}$ and
${\lambda}_{2}$ are (positive) Lagrange multipliers associated with
the two constraint equations respectively.  Since the full CSM must
be negative semidefinite, the same must be true of its diagonal
blocks, ${\sf{\Sigma}}^{1}$ and ${\sf{\Sigma}}^{2}$.  Adding the
symmetry condition, which requires one off-diagonal matrix to be the
transpose of the other, we arrive at the following structure for the
CSM:
\begin{equation}
{\lambda}_{1}^{-1} \begin{bmatrix}
{{\lambda}_{1}}^{2}{\sf{\Sigma}}^{1} & {\lambda}_{1}{\lambda}_{2}
{\sf{\Sigma}}^{1} \\ {\lambda}_{2} {\lambda}_{1} {\sf{\Sigma}}^{1} &
{{\lambda}_{2}}^{2} {\sf{\Sigma}}^{1}
\end{bmatrix}.
\end{equation}
The high redundancy of this matrix is now fully manifest, since its
rank is at most equal to that of its building block,
${\sf{\Sigma}}^{1}$. Moreover, since this building block itself must
obey the two (independent) constraint conditions ${\sum}_{\mu
=1}^{M} {p}^{1}_{\mu}{\sf{\Sigma}}^{1}_{\mu \nu}={\sum}_{\mu =1}^{M}
{{p}}_{\mu}^{2}{{\sf{\Sigma}}^{1}}_{\mu \nu}=0$, its rank cannot
exceed $M-2$, implying the same for the full CSM.  This conclusion
confirms our earlier rank determination for the CSM on the basis of
Theorem 2-(iii).

The generalization of the two-constraint model of Eq.~(\ref{34}) to
the $K$-constraint case,
\begin{equation}
{\max \;}_\mathbf{x} {\;} U(\mathbf{x}) {\:} {\;} s.t. {\:} {\:}
{m}^{k}-{\mathbf{p}^{k}} \cdot \mathbf{x}= 0, {\:} k=1,2,\ldots,K
\leq M, \label{35}
\end{equation}
is straightforward and closely parallels the above development.
Accordingly, the resulting CSM is a simple extension of the block form
found above, as follows:

\textbf{Property (iii)} \textit{The CSM for the multiple-constraint
utility maximization model assumes a $K \times K$ block form with
the $(ij)$th block equal to $({\lambda}_{i}{\lambda}_{j} /
{\lambda}_{1}) {\sf{\Sigma}}^{1}$.}

The basic block ${\sf{\Sigma}}^{1}$, and consequently full CSM as well,
must now obey $K$ (independent) constraint equations ${\sum}_{\mu =1}^{M}
{p}^{k}_{\mu}{\sf{\Sigma}}^{1}_{\mu \nu}=0$, $k=1,2,\ldots,K$.
Consequently, the rank of the basic block ${\sf{\Sigma}}^{1}$ cannot
exceed $M-K$, with the same implied for the full CSM.  Clearly, each added
constraint lowers the rank of the CSM by adding a new zero to its spectrum
while in general reducing the optimized utility level.  As the number of
constraints approaches the dimension of the consumption bundle, i.e., as
$K \rightarrow M$, the optimization process becomes progressively less
relevant in determining the chosen bundle, while, correspondingly, the CSM
loses rank and information, until it finally vanishes altogether at $K=M$.

The next step in generalizing the basic utility maximization model is the
extension to nonlinear constraints.  This problem has received
considerable attention in the literature, as it represents a first step in
generalizing the prototype constrained optimization problem.  Indeed as
stated in the Introduction, Hatta's (1980) gain method was essentially
constructed to deal with (the multiple-constraint version of) this
generalization and represents the first successful treatment of this
problem.  In the following we shall deal with this generalization (and its
multiple-constraint version) according to our general procedure, and
subsequently illustrate it by a model of consumer market power.

The generalization in question consists in replacing the linear
constraint of Eq.~(\ref{31}) with ${B}-{E}(\mathbf{x},
\mathbf{b})=0$, where ${E}(\cdot)$ is a twice continuously
differentiable function, and $\mathbf{b}$ is an $L$-dimensional
vector of parameters.  With no loss in generality, we can assume
${B}$ to be nonnegative.  The parameter set is thus identified as
$\mathbf{a} \stackrel{\rm{{def}}}{=}(\mathbf{b},{B})$, of dimension
$N=L+1$. Carrying out a procedure parallel to that followed in \S
IIA and above for the construction of the GCD's, we are led to
$D_{\alpha}(\mathbf{x},\mathbf{a})\stackrel{\rm{{def}}}{=}{\partial
\over
\partial {b}_{\alpha}} +[
\partial {E}(\mathbf{x}, \mathbf{b}) / \partial {b}_{\alpha} ] {\partial \over
\partial {B}}$, $\alpha=1,2,\ldots,L$.  Substituting these in Eq.~(\ref{7}), we find the nonlinear generalization of the Slutsky matrix
(defined to be the negative of $\Omega$ divided by the Lagrange
multiplier) in the form
\begin{equation}
 {\sf{\Theta}}_{\alpha \beta} \stackrel{\rm{{def}}}{=} {\sum}_{i=1}^{M}
[ \partial {E}_{,i}(\mathbf{x}(\mathbf{b},{B}),\mathbf{b}) /
\partial {b}_{\alpha}]
 [\partial {x}_{i}(\mathbf{b},{ B}) / \partial {b}_{\beta}+{E}_{,\beta}(\mathbf{x}(\mathbf{b},{ B}),\mathbf{b}) \partial {x}_{i}(\mathbf{b},{ B}) / \partial {B}].
\label{36}
\end{equation}
This matrix is negative semidefinite, and has a rank not exceeding the
smaller of $M-1$ and $L$ according to Theorem 2-(iii).

As a first application, let us consider a consumer who, competing in a
small market where individual consumers have market power, acts to
maximize her utility on the basis of her best estimate of other consumers'
aggregate demand (which she takes as fixed).  Specifically, the model is
defined by
\begin{equation}
{\max \;}_\mathbf{x}{\;} U(\mathbf{x}) {\:} {\;} s.t. {\:} {\:} {m}-
{\sum}_{i=1}^{M}{x}_{i}{P}_{i}({x}_{i}+ {q}_{i}^{-})=0. \label{38}
\end{equation}
Here ${P}_{i}(\cdot)$ is the (twice continuously differentiable)
inverse supply function for the $i$th good and ${q}_{i}^{-}$ is the
aggregate demand of all the other consumers for that good.  The
consumer is thus confronted with a nonlinear version of the
Slutsky-Hicks problem treated above, with $B=m$, and
$\mathbf{b}=\mathbf{q}^{-}$.

Assuming as usual the existence of a solution
$\mathbf{x}(\mathbf{q}^{-}, {m})$, we can directly proceed to the
comparative statics of this model using Eq.~(\ref{36}).  Upon
multiplying the result by $1/{p}_{\alpha}^{\prime}$ on the left and
$1/{p}_{\beta}^{\prime}$ on the right, we find
\begin{equation}
{\sf{G}}_{\alpha
\beta}\stackrel{\rm{{def}}}{=}[1+{x}_{\alpha}(\mathbf{q}^{-},
{m}){p}_{\alpha}^{\prime \prime}/{p}_{\alpha}^{\prime}][{\partial
{x}_{\alpha}(\mathbf{q}^{-}, {m}) \over \partial {q}_{\beta}^{-}}{1
\over {p}_{\beta}^{\prime}} + {\partial {x}_{\alpha}(\mathbf{q}^{-},
{m}) \over
\partial m}{x}_{\beta}(\mathbf{q}^{-}, {m})], \label{39}
\end{equation}
where
${p}_{\alpha}\stackrel{\rm{{def}}}{=}{P}_{\alpha}({x}_{\alpha}(\mathbf{q}^{-},
{m})+{q}_{\alpha}^{-})$, and ${p}_{\alpha}^{\prime}$ and
${p}_{\alpha}^{\prime \prime}$ are respectively the first and second
derivatives of ${P}_{\alpha}({x}_{\alpha}(\mathbf{q}^{-},
{m})+{q}_{\alpha}^{-})$ with respect to its argument.  Thus $\sf{G}$
is negative semidefinite and has a rank not exceeding $M-1$.

At this point, one can proceed to a consideration of the Nash
equilibrium and associated comparative statics among the consumers
competing in the market described above.  We will not pursue
equilibrium considerations here, focusing our attention instead on
the implications of Eq.~(\ref{39}) for consumer behavior under
conditions of imperfect competition.  For that purpose, however,
Eq.~(\ref{39}) must first be rewritten in terms of $\bf p$, the
price vector.  The change of parameters from $(\mathbf{q}^{-},m)$ to
$(\mathbf{p},m)$ can be implemented by recourse to Theorem 2-(v) (or
by means of the chain rule).  The result is
\begin{equation}
{\sf{\tilde{G}}}_{\alpha \beta}{=}{\sum}_{\gamma
=1}^{M}{\sf{\Sigma}}_{\alpha \gamma}{[1+{x}_{\gamma}(\mathbf{p},
{m}){p}_{\gamma}^{\prime \prime}/{p}_{\gamma}^{\prime}]}^{-
1}{\sf{J}}_{\gamma \beta}, \label{40}
\end{equation}
where ${\sf{\Sigma}}_{\alpha \gamma}$ is the standard Slutsky form
${\partial {x}_{\alpha}(\mathbf{p}, {m}) \over \partial
{p}_{\gamma}}+ {\partial {x}_{\alpha}(\mathbf{p}, {m}) \over
\partial m}{x}_{\gamma}(\mathbf{p}, {m})$ and ${\sf{J}}_{\gamma
\beta}\stackrel{\rm{{def}}}{=}{\delta}_{\gamma \beta}-
{p}_{\gamma}^{\prime}{\partial {x}_{\beta}(\mathbf{p}, {m}) \over
\partial {p}_{\gamma}}$.  Furthermore, we have introduced in Eq.~(\ref{40}) a modified CSM ${\sf{\tilde{G}}}$ by pre- and
post-multiplying ${\sf{{G}}}$ according to Theorem 2-(iv) and
Eq.~(\ref{10}) using ${[1+{x}_{\gamma}(\mathbf{p},
{m}){p}_{\gamma}^{\prime \prime}/{p}_{\gamma}^{\prime}]}^{-1}$ and
${\sf{J}}^{\dag}$ as factors. We have in Eq.~(\ref{40}) a
generalization of the Slutsky matrix to the case where the consumer
has market power and full information (of the existing demand levels
and inverse supply functions, as stipulated above).

The first point to observe is that the limit $\mathbf{p}^{\prime}
\rightarrow 0$ represents the case of perfect competition, since it
removes all possible price variations and market power effects. Thus
we expect ${\sf{G}}$ to reduce to the standard Slutsky matrix
$\sf{\Sigma}$ in that limit, a fact which is readily verified by an
inspection of Eq.~(\ref{40}).  Next, we limit the remainder of this
discussion to the case of linear supply functions in order to
simplify the analysis.  In that case, $\mathbf{p}^{\prime}$ is a
constant vector while $\mathbf{p}^{\prime \prime}$ vanishes.  Thus
Eq.~(\ref{40}) reduces to
\begin{equation}
{\sf{\tilde{G}}}_{\alpha \beta}{=}{\sf{\Sigma}}_{\alpha \beta}-
{\sum}_{\gamma =1}^{M}{\sf{\Sigma}}_{\alpha
\gamma}{p}_{\gamma}^{\prime} {\partial {x}_{\beta}(\mathbf{p}, {m})
\over \partial {p}_{\gamma}} . \label{41}
\end{equation}
This form of the CSM for our model shows the effect of market power
as a modification to the Slutsky form.  The second contribution in
Eq.~(\ref{41}) is the market power term and has an intuitively
appealing interpretation.  To see this interpretation, let
${\epsilon}_{\beta
\gamma}^{DEM}\stackrel{\rm{{def}}}{=}({p}_{\gamma}/{x}_{\beta}){
\partial {x}_{\beta}\over \partial {p}_{\gamma}}$ represent the
price elasticity matrix of consumer's demand, and
${\epsilon}_{\gamma}^{SUP}\stackrel{\rm{{def}}}{=}{p}_{\gamma}/[({x}_{\gamma
}+{q}_{\gamma}^{-}){p}_{\gamma}^{\prime}]$ be the elasticity of
supply for the $\gamma$th good.  With these definitions,
Eq.~(\ref{41}) can be written in the form
\begin{equation}
{\sf{\tilde{G}}}_{\alpha \beta}{=}{\sum}_{\gamma =1}^{M}
{\sf{\Sigma}}_{\alpha \gamma}[{\delta}_{\gamma \beta}-
{{x}_{\gamma} \over {x}_{\gamma}+{q}_{\gamma}^{-}} {{x}_{\beta}\over
{x}_{\gamma}}
{{\epsilon}_{\beta \gamma}^{DEM}\over {\epsilon}_{\gamma}^{SUP}}].
\label{42}
\end{equation}
It is now clear that the correction term is scaled by the quantity of
consumer's demand as a fraction of the aggregate demand, and otherwise
involves the ratio of the consumer's demand elasticity to the market's
supply elasticity as well as a geometrical factor.  This is of course the
sort of result one would expect, and it clearly shows how the consumer's
market power is scaled by her share of the total demand.

Although we derived the above result assuming linear supply
functions, it is also valid in cases of ``weak'' market power, i.e.,
when ${x}_{i}{p}_{i}^{\prime}/{p}_{i} \ll 1$ (and the second-order
derivatives of the inverse supply functions are negligible).  For
such cases, Eq.~(\ref{42}) gives the small, leading correction to
the Slutsky matrix arising from the consumer's market power.  It
must be pointed out here that while $\sf{\Sigma}$ is identical in
\textit{form} to the Slutsky matrix associated with a linear budget
constraint, it is nevertheless different from the latter on account
of the nonlinear constraint in the present instance.  Indeed while
${\sf{\tilde{G}}}$ is symmetric and obeys
${\sf{\tilde{G}}}\mathbf{p} =\mathbf{p}^{\dag}{\sf{\tilde{G}}}=0$,
${\sf{\Sigma}}$ is not symmetric, and obeys $\mathbf{p}^{\dag}
\sf{\Sigma}=\sf{\Sigma} \tilde{\mathbf{p}}=0$, where
$\tilde{p}_{\alpha}={p}_{\alpha}+
{p}^{\prime}_{\alpha}{x}_{\alpha}$. Moreover, the simple homogeneity
condition stated as Property (ii) above is no longer valid here.
What replaces it is a modified Euler condition in the form
$\tilde{m}\partial {x}_{\mu} / \partial m +
{\sum}_{\nu=1}^{M}{\tilde{p}}_{\nu} \partial {x}_{\mu} / \partial
{{p}}_{\nu} =0$, where $\tilde{m}=
m+{\sum}_{\alpha=1}^{M}{{x}_{\alpha}}^{2}{p}_{\alpha}^{\prime}$.
Note that all terms representing modifications to the basic utility
maximization model involve a factor of ${p}_{\alpha}^{\prime}$ and
represent adjustments necessitated by the presence of market power
effects.   In summary, we have

\textbf{Property (iv)} \textit{The model of consumer demand with
market power defined in Eq.~(\ref{38}) leads to a CSM with a
modified Slutsky structure as given in Eqs.\ (\ref{39})-(\ref{42}).
The term arising from the consumer's market power is scaled by her
fraction of the total demand.}

Our final step in generalizing the utility maximization problem is
an extension of the single nonlinear case to that of $K$ nonlinear
constraints.  Thus we consider the constraint equations ${B}^{k}-
{E}^{k}(\mathbf{x}, \mathbf{b})=0$, $k=1,2,\ldots,K$, $K<M$, where
${E}^{k}(\cdot)$ are a set of twice continuously differentiable
functions.  With no loss in generality, we can assume the $K$
parameters ${B}^{k}$ to be nonnegative.  The parameter set is thus
identified as
$\mathbf{a}\stackrel{\rm{{def}}}{=}(\mathbf{b},\mathbf{B})$, of
dimension $N=L+K$.  A procedure closely parallel to that followed
above leads us to
$D_{\alpha}(\mathbf{x},\mathbf{a})\stackrel{\rm{{def}}}{=}{\partial
\over \partial {b}_{\alpha}} + {\sum}_{k=1}^{K}[\partial
{E}^{k}(\mathbf{x}, \mathbf{b}) / \partial {b}_{\alpha} ] {\partial
\over
\partial {B}^{k}}$, $\alpha=1,2,\ldots,L$. Substituting these in
Eq.~(\ref{7}), we find the nonlinear, multiple constraint
generalization of the Slutsky matrix (defined to be the negative of
$\Omega$) in the form
\begin{eqnarray}
{\sf{\Psi}}_{\alpha \beta}&& \stackrel{\rm{{def}}}{=}
{\sum}_{i=1}^{M} [{\sum}_{k=1}^{K}{\lambda}_{k} \partial
{E}^{k}_{,i}(\mathbf{x}(\mathbf{b},\mathbf{B}),\mathbf{b}) / \partial {b}_{\alpha}] \nonumber \\
&& \times [\partial {x}_{i}(\mathbf{b},\mathbf{B}) / \partial
{b}_{\beta}+{\sum}_{k=1}^{K}{E}^{k}_{,\beta}(\mathbf{x}(\mathbf{b},\mathbf{B}),\mathbf{b})
\partial {x}_{i}(\mathbf{b},\mathbf{B}) / \partial {B}^{k}], \label{37}
\end{eqnarray}
where ${\lambda}_{k}$ are Lagrange multipliers as in Eq.~(\ref{7}).
This matrix is negative semidefinite, and has a rank no larger than
the smaller of $M-K$ and $L$ according to Theorem 2-(iii).

This concludes our treatment of a series of progressively more general
constrained optimization problems as extensions of the basic utility
maximization model.  For the remainder of this section, we shall
illustrate the above results with two specific models, both of which
involve uncertainty.  Such models typically involve constraints
originating from the fact that the probability set has unit measure.  When
expressed in direct form, i.e., that the sum of all probabilities equals
unity, such conditions constrain the parameters but not the decision
variables, and therefore do not qualify as constraints in the usual sense
(e.g., they do not conform to the constraint qualification condition in
decision space).  Indeed they play no role in determining the solution to
the optimization problem.  However, if one is interested in comparative
statics information involving the probability set, i.e., if the
probability set is in fact included among the parameters of interest, then
there arises the question of how the constraint is to be implemented in
parameter space.  We will see in the following application that there is a
natural method of implementing such constraints which maintains the
intrinsic symmetries of the problem.

Our first model here is thus an illustration of Eq.~(\ref{37})
dealing with multiple nonlinear constraints, as well as of our
method of treating problems which involve uncertainty and
probability sets.  This model is the principal-agent problem with
hidden actions, already encountered in \S IIA, where a firm, the
principal, intends to hire an individual, the agent, to work on a
specific venture on a contractual basis (see, e.g., Mas-Collel et
al. 1995).  The principal's objective in this venture is to maximize
profits, while the agent is characterized as a utility maximizer
with a known utility function.  However, the eventual outcome of the
venture, including the realized profit and utility levels, are
uncertain for two reasons.  First, the \textit{effort level} of the
agent, which in our problem can take one of two possible values
\textit{high} and \textit{low}, is unknown and unobservable to the
principal, even after the venture is completed and profits are
realized.  On the other hand, this effort level is decided by the
agent through expected utility maximization, hence the asymmetry of
information between the principal and the agent.  Second, the
venture's profits are random and unpredictable. They are specified
stochastically according to a probability distribution which we take
to be discrete for simplicity.  As already stipulated, the realized
profit level, which is observable to both sides, does not reveal the
agent's effort level; in other words, any realizable profit level
can result from either effort level of the agent.  The principal's
problem is to design a contract that maximizes the firm's expected
profits. Since the agent's effort level is a choice variable in the
principal's profit maximization problem, the latter can be
formulated as a pair of maximization problems, one for each effort
level, and the optimum decided by comparing the results.

The basic problem then is to maximize the principal's expected profits,
assuming a given effort level for the agent.  Accordingly, it can be
formulated as follows:
\begin{eqnarray}
{\max \;}_\mathbf{x} {\;}{\sum}_{i=1}^{M}({\pi}_{i}-{x}_{i}){\cal
P}_{i}^{I} \;\; s.t. {\:} {\:} {\sum}_{i=1}^{M}v({x}_{i}){\cal
P}_{i}^{I}-{c}^{I} \geq
\bar{ u} && , \nonumber \\
{\sum}_{i=1}^{M}v({x}_{i}){\cal P}_{i}^{I}-{c}^{I} \geq
{\sum}_{i=1}^{M}v({x}_{i}){\cal P}_{i}^{II}-{c}^{II}, \;\;
{\sum}_{i=1}^{M}{\cal P}_{i}^{k}=1 && , \;k=I,II. \label{043}
\end{eqnarray}
Here ${\cal P}_{i}^{k}$, where $0 < {\cal P}_{i}^{k} < 1$,
$i=1,2,\ldots, M$, $k=I,II$, is the probability of profit level
${\pi}_{i}$ given the effort level $k$, with $I$ and $II$
corresponding to high and low effort, respectively.  The decision
variable ${x}_{i}$, on the other hand, is the agent's compensation
in case the $i$th profit level is realized, so that the vector
$\mathbf{x}$ specifies the compensation schedule and in effect
defines the contract.  The agent's utility function is of the
separable variety $v(x)-{c}^{k}$, where $x$ represents the
compensation and ${c}^{k}$ the disutility of working at effort level
$k$.  On the other hand, $\bar{ u}$ is the price of the agent's
services (determined exogenously through competitive markets) in
utility terms; it is often referred to as the agent's
\textit{reservation} utility level. Here we have assumed
${c}^{I}>{c}^{II}$, corresponding to the fact that the agent prefers
low effort to high effort, \textit{ceteris paribus}.  We have also
assumed a high effort level ($I$) for the agent in Eq.~(\ref{043}),
as is apparent from the objective function and the second constraint
in the above formulation.  Once the optimum contracts for the above
problem and its conjugate, which is arrived at by interchanging $I$
and $II$ in Eq.~(\ref{043}), are determined, the principal chooses
the more profitable compensation schedule and offers the
corresponding contract to the agent. We will assume in the following
that the principal finds $k=I$ to be the more profitable choice,
with no loss of generality.  We will also assume the utility
function $v(\cdot)$ to be twice continuously differentiable,
monotonically increasing and concave, and moreover that there is a
unique internal solution to the above maximization problem (and its
conjugate). Furthermore, we will assume that both inequality
constraints bind, since an inequality that does not bind (i.e., one
which is satisfied as a strict inequality at the assumed solution)
has no bearing on the comparative statics of the problem. It should
also be emphasized that the last pair of equations above only
constrain the parameters and are not constraints in the usual sense,
a point that was discussed following Eq.~(\ref{37}).

Under the above-stated assumptions, the principal's optimization problem
can be rewritten as follows:
\begin{equation}
{\min \;}_\mathbf{x} {\;}{\sum}_{i=1}^{M}{x}_{i}{\cal P}_{i}^{I} \;
\; s.t. {\:} {\:} {B}^{k}-{\sum}_{i=1}^{M}v({x}_{i}){\cal P}_{i}^{k}
=0, \; \; {s}^{k}-{\sum}_{i=1}^{M}{\cal P}_{i}^{k}=0, \; k=I,II,
\label{044}
\end{equation}
where ${B}^{k}\stackrel{\rm{{def}}}{=}{c}^{k}+\bar{ u}$ and $0<
{s}^{k}$ are a pair of auxiliary variables which will be set equal
to unity at a convenient point in the course of the analysis.

At this juncture we recognize that Eq.~(\ref{044}) is of the general
nonlinear, multiple-constraint variety treated in Eq.~(\ref{37}).
This problem was also considered in \S IIA, where we constructed a
set of GCD's by a two-term compensation procedure based upon the GCD
set for the basic utility maximization problem. The set thus
obtained is $D_{\alpha}(\mathbf{x},\mathbf{a})
\stackrel{\rm{{def}}}{=} \partial / \partial {\cal P}_{\alpha}^{I} +
v({x}_{\alpha}) \partial / \partial {B}^{I} + \partial /
\partial {s}^{I}$ for $\alpha =1,2,\ldots,M$ and $D_{\alpha}(\mathbf{x},\mathbf{a}) \stackrel{\rm{{def}}}{=} \partial / \partial {\cal P}_{\alpha - M}^{II}
+ v({x}_{\alpha - M}) \partial / \partial {B}^{II} + \partial /
\partial {s}^{II}$ for $\alpha =M+1,M+2,\ldots,2M$.  These
definitions can conveniently be combined according to
$D_{\alpha}(\mathbf{x},\mathbf{a}) \stackrel{\rm{{def}}}{=}
{d}^{I}_{\alpha}$ for $\alpha =1,2,\ldots,M$ and
$D_{\alpha}(\mathbf{x},\mathbf{a}) \stackrel{\rm{{def}}}{=}
{d}^{II}_{\alpha-M}$ for $\alpha =M+1,M+2,\ldots,2M$ to obtain
$${d}^{k}_{i}={\partial \over \partial {\cal P}_{i}^{k}}+ v({x}_{i})
{\partial \over \partial {B}^{k}} + {\partial \over \partial {s}^{k}},
{\;}i=1,2,\ldots,M, {\;}k=I,II.$$

To proceed with the construction of the CSM, let us recall our
earlier stipulation that a comparison of the maximized profits for
the two effort levels $I$ and $II$ by the principal has revealed the
former to be optimal.  To derive the comparative statics information
corresponding to this case, we substitute the GCD set derived above
in Eq.~(\ref{7}) to arrive at the desired CSM.  The result is a
negative semidefinite matrix which may be written in the
block-matrix form
\begin{equation}
\begin{bmatrix}
{\sf{\Phi}}^{11}&{\sf{\Phi}}^{12} \\ {\sf{\Phi}}^{21}&
{\sf{\Phi}}^{22} \end{bmatrix},
\end{equation}

 where ${\sf{\Phi}}^{kk'}_{ij}
\stackrel{\rm{{def}}}{=} [2-k-
{\lambda}_{k}{v}^{\prime}({x}_{i})]{d}^{k'}_{j}
{x}_{i}(\mathbf{a})$, $i,j=1,2,\ldots,M$ and $k,k'=I,II$, with the
numerical assignments of $I \stackrel{\rm{{def}}}{=}1$ and $II
\stackrel{\rm{{def}}}{=}2$ understood here. Furthermore,
${v}^{\prime}(\cdot)$ is the derivative of $v(\cdot)$ and
${\lambda}_{k}$ is a Lagrange multiplier in the foregoing
expressions.  It is also useful to record the first order conditions
at this point:
$${\cal P}^{I}_{i}-{\lambda}_{I}{\cal P}^{I}_{i}{v}^{\prime}({x}_{i})
-{\lambda}_{II}{\cal P}^{II}_{i}{v}^{\prime}({x}_{i})=0.$$
We will use these conditions in the course of the analysis.

It is convenient at this point to eliminate the auxiliary parameters
${s}^{I,II}$ from our results.  Since partial derivatives with
respect to these parameters occur only in terms of the form
${d}^{k}_{j} {x}_{i}$ in our expressions, they can be eliminated by
recourse to a certain scale symmetry of the underlying problem. This
symmetry can be readily recognized by an inspection of
Eq.~(\ref{044}).  To that end, observe that a rescaling of the
(augmented) parameter set $(\mathbf{a}) \stackrel{\rm{{def}}}{=}
(\mathbf{a}^{I},\mathbf{a}^{II}) \stackrel{\rm{{def}}}{=}({\cal
P}^{I}, {B}^{I}, {s}^{I},{\cal P}^{II}, {B}^{II}, {s}^{II})$
according to $\mathbf{a}^{k} \rightarrow {\mu}^{k}\mathbf{a}^{k}$,
where the scale constants ${\mu}^{I,II}$ are a pair of positive
numbers, leaves the problem unchanged.  This invariance condition
then implies, via Theorem 2-(ii), the following homogeneity
condition:
\[
\{ {\sum}_{j=1}^{M}{\cal P}_{j}^{k} {\partial \over \partial {\cal
P}_{j}^{k}} + {B}^{k}{\partial \over \partial {B}^{k}} + {s}^{k}
{\partial \over \partial {s}^{k}} \} \mathbf{x}(\mathbf{a})=0, \; \;
k=I,II.
\]
Using this homogeneity condition, one can eliminate all derivative terms
with respect to the auxiliary parameters ${s}^{k}$, and then proceed to
set them equal to unity.  The result of this elimination is simply the
replacement of the GCD set ${d}^{k}_{j}$ with
\[ {\cal D}_{i}^{k} \stackrel{\rm{{def}}}{=} {\partial \over \partial {\cal
P}_{i}^{k}} - {\sum}_{j=1}^{M}{\cal P}_{j}^{k} {\partial \over \partial
{\cal P}_{j}^{k}} + [v({x}_{i}) - {B}^{k}]{\partial \over \partial
{B}^{k}}, {\;}i=1,2,\ldots, M, {\;}k=I,II. \]

We pause at this point to underline the role played by the auxiliary
parameters in dealing with the pair of constraints ${\sum}_{i=1}^{M}
{\cal P}_{i}^{I,II}=1$ in Eq.~(\ref{044}) which constrain its
parameter space but not its decision space.  In view of the fact
that these constraints do not allow a change in one of the
probabilities while the others are kept fixed, a partial derivative
of the sort $\partial \mathbf{x}(\mathbf{a}) / \partial {\cal
P}_{i}^{k}$ does not correspond to a realizable scenario in the real
world.  On the other hand, one can envision a change in a given
probability with the compensating change required by the constraint
symmetrically allocated to all the probabilities.  It is precisely
this objective of enforcing the constraint in a symmetrical manner
that is accomplished by the introduction of the auxiliary parameters
${s}^{I,II}$.  The end result, which can be gleaned by an inspection
of the expression for ${\cal D}_{i}^{k}$ given above, is the
replacement of the naive derivative ${\partial / \partial {\cal
P}_{i}^{k}}$ by the constraint-conforming combination ${\partial /
\partial {\cal P}_{i}^{k}} - {\sum}_{j=1}^{M}{\cal P}_{j}^{k} {\partial /
\partial {\cal P}_{j}^{k}}$.  Observe that the role played by the
auxiliary variables here is entirely analogous to that of Lagrange
multipliers in a constrained optimization problem, the major distinction
being that the latter are used to enforce decision space constraints.
Needless to say, the parameter constraints present here are typical of any
optimization problem which involves uncertainty and is defined in terms of
random variables drawn from probability sets.  The method of auxiliary
variables used here is thus a natural and effective way of dealing with
problems involving uncertainty.

Returning to the CSM derived above, we now proceed to exploit its
symmetry properties.  Thus considering the $2 \times 2$ block form
of the CSM, we infer on the basis of its semidefiniteness that the
diagonal blocks ${\sf{\Phi}}^{11}$ and ${\sf{\Phi}}^{22}$ are
symmetric, $M \times M$ matrices, while
${\sf{\Phi}}^{12}={{\sf{\Phi}}^{21}}^{\dag}$.  The first and last of
these conditions together with the first order conditions given
above may be used to infer the equality ${\sf{\Phi}}^{22}_{ij}=
{R}_{i}{\sf{\Phi}}^{11}_{ij}{R}_{j}$, while the remaining condition
then implies that ${\sf{\Phi}}^{12}_{ij}=
{\sf{\Phi}}^{11}_{ij}{R}_{j}$, where
$0>{R}_{i}\stackrel{\rm{{def}}}{=}-{\cal P}^{I}_{i}/{\cal
P}^{II}_{i}$. These relations lead to the ``double-Slutsky'' form
already found for the model defined in Eq.~(\ref{34}), namely
\[
\begin{bmatrix}{\sf{\Phi}}&{\sf{\Phi}}{\sf{R}} \\
{\sf{R}}{\sf{\Phi}}& {\sf{R}}{\sf{\Phi}}{\sf{R}} \end{bmatrix},
\]
where $\sf{\Phi} \stackrel{\rm{{def}}}{=} \sf{\Phi}^{11}$ and
$\sf{R}$ is a diagonal matrix whose $i$th diagonal element is
${R}_{i}$. This structure clearly implies that the CSM is highly
redundant, its rank being the same as that of each of its blocks,
and since ${\sf{R}}$ is a negative definite matrix, the properties
of the full CSM are fully conveyed by the statement that
${\sf{\Phi}}$ is a negative semidefinite matrix.  Note that this
matrix is at most of rank $M-2$, since the pair of (decision space)
constraints in Eq.~(\ref{044}) imply the existence of a pair of zero
eigenvalues for it.  Thus the rank of the full CSM is no larger than
$M- 2$, which is precisely what is predicted by Theorem 2-(iii).

To interpret the CSM just derived, let us first denote the minimized
value of the objective function in Eq.~(\ref{044}) by
$C(\mathbf{a})$.  This is of course the principal's minimum expected
cost for inducing the agent to perform at effort level $I$.  On the
other hand, observe that since
${B}^{k}\stackrel{\rm{{def}}}{=}{c}^{k}+\bar{ u}$, we have the
envelope result that $\partial C(\mathbf{a}) / \partial {c}^{k}
={\lambda}_{k} $, where we recall that ${c}^{k}$ is the agent's
disutility for performing at effort level $k$.  Therefore we find,
not surprisingly in the light of the basic utility maximization
problem, that ${\lambda}_{k}$ represents the expected marginal cost
to the principal of the agent's disutility for working at effort
level $k$.  Recalling that ${c}_{I} > {c}_{II}$, and that the
contract offered to the agent is intended to induce level $I$
performance, we can conclude that ${\lambda}_{I} \geq 0$ while
${\lambda}_{II} \leq 0$.  A more detailed analysis involving the
first-order conditions and the constraints together with the
properties of the agent's value function confirms these intuitive
conclusions. Furthermore, since by assumption ${\cal P}^{k}_{i}>0$
as well as ${v}^{\prime}({x}_{i})>0$ for $k=I,II$ and
$i=1,2,\ldots,M$, we may conclude, using the first order conditions,
that $1- {\lambda}_{I}{v}^{\prime}({x}_{i}) \leq 0$.

We are now in a position to state the comparative statics results
for the principal-agent problem in a more useful form.  First, let
us consider the negative semidefinite matrix
${\sf{\Phi}}^{22}_{ij}=- {\lambda}_{II}{v}^{\prime}({x}_{i}){\cal
D}^{II}_{j} {x}_{i}(\mathbf{a})$. Since ${\lambda}_{II} \leq 0$, we
have the result that the $M \times M$ matrix
\begin{equation}
{\sf{H}}_{ij} \stackrel{\rm{{def}}}{=} {\partial
v({x}_{i}(\mathbf{a})) \over
\partial {\cal P}_{j}^{II}} - {\sum}_{l=1}^{M}{\cal P}_{l}^{II} {\partial
v({x}_{i}(\mathbf{a})) \over \partial {\cal P}_{l}^{II}} +
[v({x}_{j}(\mathbf{a})) - {c}^{II}- \bar{ u}]{\partial
v({x}_{i}(\mathbf{a})) \over \partial {c}^{II}} \label{045}
\end{equation}
is also negative semidefinite and has a rank no larger than $M-2$.
Note that we have set $\partial {x}_{j}(\mathbf{a}) / \partial
{B}^{II}=\partial {x}_{j}(\mathbf{a}) / \partial {c}^{II}$, an
equality which can be ascertained by reference to
${B}^{k}\stackrel{\rm{{def}}}{=}{c}^{k}+\bar{ u}$.  Note also that
we have stated the CSM in terms of
${v}^{\prime}({x}_{j}(\mathbf{a}))$ rather than
${x}_{j}(\mathbf{a})$, a choice that simplifies the resulting
expressions.  The structural similarity of the result in
Eq.~(\ref{045}) to the Slutsky-Hicks result is quite evident.  The
extra terms here, on the other hand, result from the appearance of
the parameter constraints and the (generally nonlinear) function $v(
\cdot )$ in the present problem.  For example, Eq.~(\ref{045})
implies that the appropriately compensated change in
${x}_{i}(\mathbf{a})$, the wages offered by the principal in case
profit level $i$ is realized, as a result of an increase in ${\cal
P}_{j}^{II}$, the probability of the $j$th profit level conditional
on effort level $II$, is negative.  This is of course the expected
response inasmuch as the principal's objective is to induce the
agent to effort level $I$, hence \textit{away} from effort level
$II$.

The result in Eq.~(\ref{045}) can be equivalently stated in terms of
level $I$ parameters.  To wit, the negative semidefinite matrix
${\sf{\Phi}}^{11}$ can be restated as
\begin{equation}
{\sf{\Phi}}_{ij}^{11}=[1-{\lambda}_{l}{v}^{\prime}({x}_{i}(\mathbf{a}))]{\cal
D}^{I}_{j} {x}_{i}(\mathbf{a})={\cal D}^{I}_{j} {x}_{i}(\mathbf{a})
- {\lambda}_{I} {\cal D}^{I}_{j} v({x}_{i}(\mathbf{a})), \label{046}
\end{equation}
where ${\cal D}^{I}_{j} $ is the level $I$ counterpart of the
compensated derivative in Eq.~(\ref{045}).  The full
semidefiniteness results in Eqs.\ (\ref{045}) and \ (\ref{046})
imply, among other things, the semidefiniteness of the respective
diagonal elements, as summarized in

\textbf{Property (v)} \textit{The comparative statics properties of
the principal-agent problem defined above are conveyed by the
negative semidefinite matrices $\sf{\Phi}$ and $\sf H$.  In
particular, we have the inequalities}
\begin{eqnarray}
{\partial {x}_{i}(\mathbf{a}) \over \partial {\cal P}_{i}^{I}} -
{\sum}_{j=1}^{M}{\cal P}_{j}^{I} {\partial {x}_{i}(\mathbf{a}) \over
\partial {\cal P}_{j}^{I}} + [v({x}_{i}(\mathbf{a})) - {c}^{I}- \bar{
u}]{\partial
{x}_{i}(\mathbf{a}) \over \partial {c}^{I}} \geq 0 , \nonumber \\
{\partial {x}_{i}(\mathbf{a}) \over \partial {\cal P}_{i}^{II}} -
{\sum}_{j=1}^{M}{\cal P}_{j}^{II} {\partial {x}_{i}(\mathbf{a})
\over
\partial {\cal P}_{j}^{II}} + [v({x}_{i}(\mathbf{a})) - {c}^{II}- \bar{
u}]{\partial {x}_{i}(\mathbf{a}) \over \partial {c}^{II}} \leq 0.
\nonumber
\end{eqnarray}
Note the sign reversal in going from level $I$ derivatives to those of
level $II$.  Recalling that the contract was designed to induce the agent
to level $I$ effort, one can comprehend the above inequalities as
compensated adjustments in the contract's wage schedule, $\bf x$, in
response to changing expectations, $\cal P$, which counteract any expected
depreciation in the level $I$ prospects or appreciation in those of level
$II$, while keeping costs as low as possible.  As already emphasized, the
above analysis and its results pertain to the case where the contract is
designed to induce effort level $I$.  While the analysis for the
complementary case is parallel to the foregoing, the results are not
expected to be symmetrical with respect to an interchange of $I$ and $II$,
since the condition ${c}^{I} > {c}^{II}$ breaks the symmetry between the
two cases.

We conclude our discussion of the principal-agent problem by
remarking that the formulation of this problem given in
Eq.~(\ref{044}) is formally identical to that of cost minimization
subject to constraints, a fact that explains the observed
similarities to the Slutsky-Hicks problem.

Our final application deals with the problem of selecting an efficient
portfolio from a set of $M$ financial assets ${s}_{i}$ under idealized
market conditions, including the possibility of a riskless asset and short
sales.  Treating returns as random variables, the model seeks to find the
portfolio which achieves a prescribed level of expected return with a
minimum variance.  This problem is of the generalized utility maximization
category defined earlier.  It is defined by (Fama and Miller 1972, Elton
and Gruber 1991)
\begin{equation}
{\min \;}_\mathbf{x} {\;}
\mathbf{x}^{\dag}{\sf{\sigma}}\mathbf{x}{\:} {\;} s.t. {\:} {\:}
{\cal W}=\mathbf{w}\cdot \mathbf{x}{\:} {\:} \rm{ and} {\:} {\:}
{\cal R} =\mathbf{r}\cdot\mathbf{x}, \label{43}
\end{equation}
where ${x}_{i}$ represents the fraction of the total investment
allocated to asset ${s}_{i}$, ${\sf{\sigma}}$ stands for the
covariance matrix of the asset set, ${\cal W}\stackrel{\rm{
def}}{=}1$ and ${w}_{i}\stackrel{\rm{ def}}{=}1$, $i=1,2,\dots,M$,
are auxiliary variables introduced for notational convenience,
${r}_{i}$ represents the expected rate of return for asset
${s}_{i}$, and $\cal R$ is the prescribed rate of return for the
portfolio.  The possibility of short sales implies that the vector
$\bf x$ is unrestricted in the sign and magnitude of its components.

The problem posed in Eq.~(\ref{43}) and variants of it have been
extensively analyzed by various methods.  Our objective here is to
derive and present the comparative statics information associated
with it in a compact and useful manner.  As it stands, the
formulation given in Eq.~(\ref{43}) leads to a highly redundant CSM.
For example, the block containing partial derivatives with respect
to ${\sf{\sigma}}_{ij}$ will be of order $M(M+1)/2$ (the number of
independent elements in the covariance matrix) while its rank cannot
exceed $M-2$, which for a typical case with $M \gg 1$ is far lower
than its order.  The scenario underlying this difficulty, as well as
its converse where the number of decision variables far exceeds the
number of parameters, is not uncommon and naturally arises in a
number of problems.  We will deal with this difficulty by
transforming Eq.~(\ref{43}) to a more suitable form which, it turns
out, is quite useful in other respects as well.  The basic idea here
is the simple observation that the entire analysis would be greatly
simplified if the available assets were in fact uncorrelated, i.e.,
if the covariance matrix were diagonal.  Actual assets are of course
significantly correlated, but one can always construct certain mixes
of them, to be called \textit{principal portfolios}, that are
entirely uncorrelated.  These special mixes are directly determined
by the eigenvectors of the covariance matrix.  In effect, the
problem of stock selection from existing, \textit{correlated assets}
is traded for the simpler problem of choosing from a set of
\textit{uncorrelated principal portfolios}, which for all intents
and purposes act as new assets themselves. Mathematically, this
transformation amounts to a simple change of basis from the initial
assets to the principal portfolio basis where the covariance matrix
has a diagonal representation.

Let us now implement the above ideas.  First, note that by virtue of
its symmetry, $\sf{\sigma}$ admits of a set of $M$ orthogonal
eigenvectors $\mathbf{e}^{\mu}$, $\mu=1,2,\ldots,M$.  With no loss
in generality, we will take these vectors to be of unit length, so
that $\mathbf{e}^{\mu}\cdot\mathbf{e}^{\nu}={\delta}_{\mu \nu}$. The
covariance matrix itself can then be represented in terms of its
eigenvalues and eigenvectors in the form
${\sf{\sigma}}_{ij}={\sum}_{\mu=1}^{M}
{\sigma}_{\mu}^{2}{e}^{\mu}_{i} {e}^{\mu}_{j}$, where
${\sigma}_{\mu}^{2} \geq 0$ are the eigenvalues of the covariance
matrix, or the \textit{principal variances}.  The principal
portfolios are now defined by ${S}_{\mu}\stackrel{\rm{
def}}{=}{\sum}_{i=1}^{M}{e}^{\mu}_{i}{s}_{i}$, i.e., the principal
portfolio ${S}_{\mu}$ is a standard mix which contains an amount
${e}^{\mu}_{i}$ of asset ${s}_{i}$.  Each principal portfolio is
characterized by a \textit{weight} ${W}_{\mu}\stackrel{\rm{
def}}{=}{\sum}_{i=1}^{M}{e}^{\mu}_{i}{w}_{i}$, an expected rate of
return ${R}_{\mu}\stackrel{\rm{ def}}{=}
{\sum}_{i=1}^{M}{e}^{\mu}_{i}{r}_{i}/{W}_{\mu}$, and a variance
${\sum}_{j=1}^{M}{\sum}_{i=1}^{M} {e}_{i}^{\mu}{\sf{\sigma}}_{ij}
{e}_{j}^{\mu}={\sigma}_{\mu}^{2}$, to be referred to as
\textit{principal} weight, return, and variance, respectively.  A
general portfolio, initially specified by ${x}_{i}$, is now
described by ${X}_{\mu}$, which determines what share of the total
investment is to be allocated to principal portfolio ${S}_{\mu}$.
These equivalent specifications are related by the transformation
equations ${x}_{i}={\sum}_{\mu=1}^{M}{e}_{i}^{\mu}{X}_{\mu}$ and
${X}_{\mu}={\sum}_{i=1}^{M}{e}^{\mu}_{i}{x}_{i}$.  The pair of
constraints ${\cal W}-\mathbf{w}\cdot\mathbf{x}=0$ and ${\cal
R}-\mathbf{r}\cdot\mathbf{x}=0$, on the other hand, are transformed
to ${\cal W}-\mathbf{W}\cdot\mathbf{X}=0$ and ${\cal
R}-\mathbf{R}\cdot\mathbf{X}=0$.  A riskless asset, if present,
would constitute a principal portfolio by itself, with vanishing
variance, unit weight, and return equal to the riskless rate of
return.  As is well known, the efficient portfolio in the presence
of a riskless asset is a simple mix of the latter and an efficient
portfolio composed of the risky ones. Accordingly, we will treat the
risky assets separately, excluding riskless assets from
consideration in the following.  We will thus assume the principal
variances to be positive definite, since a principal portfolio with
a vanishing principal variance is essentially equivalent to a
riskless asset and can be treated separately as stipulated.

The efficient portfolio is now defined by the optimization problem
\begin{equation}
{\min \;}_\mathbf{X} {\;}
{\sum}_{\mu=1}^{M}{\sigma}_{\mu}^{2}{{X}_{\mu}}^{2} {\:} {\;} s.t.
{\:} {\:} {\cal W}-\mathbf{W}\cdot\mathbf{X}=0 {\:}{\:}\rm{ and}{\:}
{\:} {\cal R}-\mathbf{R}\cdot \mathbf{X}=0, \label{44}
\end{equation}
where ${\sigma}_{\mu}^{2} >0$.  This is the principal-portfolio
version of the problem posed in Eq.~(\ref{43}).  The problem posed
in Eq.~(\ref{44}), on the other hand, can be solved in a
straightforward manner. The solution is most conveniently expressed
in terms of the rescaled vectors $\bar{{W}}_{\mu}\stackrel{\rm{
def}}{=}{W}_{\mu}/{\sigma}_{\mu}$ and $\bar{{R}}_{\mu}\stackrel{\rm{
def}}{=}{R}_{\mu}/{\sigma}_{\mu}$, where ${\sigma}_{\mu}$ is the
positive square root of ${\sigma}_{\mu}^{2}$.  We find, in terms of
these rescaled vectors, the solution
\begin{equation}
{X}_{\mu}(\vec{{\sigma}^2},\mathbf{W},{\cal W},\mathbf{R},{\cal R})=
[{\lambda}_{1}(\bar{\mathbf{W}},{\cal W},\bar{\mathbf{R}},{\cal
R})\bar{{W}}_{\mu}+{\lambda}_{2}(\bar{\mathbf{W}},{\cal
W},\bar{\mathbf{R}},{\cal R})\bar{{R}}_{\mu}] /2{\sigma}_{\mu},
\label{45}
\end{equation}
where ${\lambda}_{1}\stackrel{\rm{ def}}{=}2 [
\bar{\mathbf{R}}\cdot({\cal W}\bar{\mathbf{R}}-{\cal
R}\bar{\mathbf{W}})]/{\cal D}$ and ${\lambda}_{2}\stackrel{\rm{
def}}{=}-2[\bar{\mathbf{W}}\cdot({\cal W}\bar{\mathbf{R}}-{\cal
R}\bar{\mathbf{W}})]/{\cal D}$ with ${\cal D}\stackrel{\rm{
def}}{=}(\bar{\mathbf{W}}\cdot\bar{\mathbf{W}})(\bar{\mathbf{R}}\cdot\bar{\mathbf{R}})-
{(\bar{\mathbf{W}} \cdot \bar{\mathbf{R}})}^{2}$, are the pair of
Lagrange multipliers associated with the two constraints in
Eq.~(\ref{44}), and $\vec{{\sigma}^2}\stackrel{\rm{ def}}{=}
({\sigma}^2_{1},\ldots,{\sigma}^2_{M})$.  On the other hand, the
minimized value of the objective function, to be called the
\textit{portfolio variance} and denoted by ${\sigma}^{2}_{P}$, is
given by the simple formula
\begin{equation}
{\sigma}^{2}_{P}(\bar{\mathbf{W}},{\cal W},\bar{\mathbf{R}},{\cal
R})={ ({\cal W}\bar{\mathbf{R}}-{\cal R}\bar{\mathbf{W}})\cdot
({\cal W}\bar{\mathbf{R}}-{\cal R}\bar{\mathbf{W}}) \over
(\bar{\mathbf{W}}\cdot\bar{\mathbf{W}})(\bar{\mathbf{R}}\cdot\bar{\mathbf{R}})-{(\bar{\mathbf{W}}
\cdot \bar{\mathbf{R}})}^{2} } \label{46}
\end{equation}
The emergence of a quadratic dependence of ${\sigma}^{2}_{P}$ on the
portfolio return $\cal R$, a well-known result of portfolio
analysis, should be noted.  The minimum of ${\sigma}^{2}_{P}$ with
respect to $\cal R$ occurs at ${\cal R}^{\star}\stackrel{\rm{
def}}{=}{\cal
W}\bar{\mathbf{W}}\cdot\bar{\mathbf{R}}/\bar{\mathbf{W}}\cdot\bar{\mathbf{W}}$,
and it simply equals ${\cal
W}^{2}{(\bar{\mathbf{W}}\cdot\bar{\mathbf{W}})}^{-1}={\cal
W}^{2}{[{\sum}_{\mu=1}^{M} {({W}_{\mu}/{\sigma}_{\mu})}^{2}]}^{-1}$,
again a remarkably simple result.  In summary, then, we have

\textbf{Property (vi)} \textit{Eqs.\ (\ref{45}) and (\ref{46}) et
seq. give a complete solution of the efficient portfolio problem by
the method of principal portfolios.}

Returning to the comparative statics of the principal portfolio
problem posed in Eq.~(\ref{44}), we note that this is a simple
extension of the problem posed in Eq.~(\ref{34}) thanks to the
parameter separation between the objective and constraint functions.
This separation guarantees that the partial derivatives with respect
to ${\sigma}_{\mu}^{2}$ are already in compensated form with respect
to the constraints.  Hence we can take the GCD's to be defined by
${D}_{\alpha}(\mathbf{X},
\mathbf{a})\stackrel{\rm{{def}}}{=}\partial /
\partial {\sigma}_{\alpha}^{2}$ for $\alpha=1,2,\ldots,M$, by
${D}_{\alpha}(\mathbf{X},
\mathbf{a})\stackrel{\rm{{def}}}{=}\partial /
\partial {W}_{\alpha}+{X}_{\alpha} \partial / \partial {\cal W}$ for
$\alpha=M+1,M+2,\ldots,2M$, and by ${D}_{\alpha}(\mathbf{X},
\mathbf{a})\stackrel{\rm{{def}}}{=}\partial / \partial {R}_{\alpha}
+{X}_{\alpha}
\partial / \partial {\cal R}$ for $\alpha=2M+1,2M+2,\ldots,3M$.  Using
these in Eq.~(\ref{7}), one can derive the CSM in a straightforward
manner.  The result is the $3 \times 3$ block matrix
\begin{equation}
\begin{bmatrix}
\partial {{X}_{\mu}}^{2} / \partial {\sigma}_{\nu}^{2}
& 2{X}_{\mu}{\sf{\Sigma}}^{W}_{\mu \nu} &
2{X}_{\mu}{\sf{\Sigma}}^{R}_{\mu \nu}\cr -{\lambda}_{1} \partial
{{X}_{\mu}} / \partial {\sigma}_{\nu}^{2} &-
{\lambda}_{1}{\sf{\Sigma}}^{W}_{\mu \nu}&-
{\lambda}_{1}{\sf{\Sigma}}^{R}_{\mu \nu} \\ -{\lambda}_{2}
\partial {{X}_{\mu}} / \partial {\sigma}_{\nu}^{2} &-
{\lambda}_{2}{\sf{\Sigma}}^{W}_{\mu \nu}&-
{\lambda}_{2}{\sf{\Sigma}}^{R}_{\mu \nu} \end{bmatrix}, \label{47}
\end{equation}
where ${\sf{\Sigma}}_{\mu \nu}^{W} \stackrel{\rm{ def}}{=}\partial
{X}_{\mu} / \partial {{W}}_{\nu}+{X}_{\nu} \partial {X}_{\mu}/
\partial {\cal W}$, ${\sf{\Sigma}}_{\mu \nu}^{R} \stackrel{\rm{
def}}{=}\partial {X}_{\mu} /
\partial {{R}}_{\nu}+{X}_{\nu} \partial {X}_{\mu}/ \partial {\cal R}$,
with $\mu,\nu=1,2,\dots,M$.  The full matrix in Eq.~(\ref{47}) is
negative semidefinite, which implies that its diagonal blocks,
$\partial {{X}_{\mu}}^{2} / \partial {\sigma}_{\nu}^{2}$,
$-{\lambda}_{1}{\sf{\Sigma}}^{W}_{\mu \nu}$, and
$-{\lambda}_{2}{\sf{\Sigma}}^{R}_{\mu \nu}$ are also negative
semidefinite. The symmetry property of the full matrix, on the other
hand, implies that $\partial {{X}_{\mu}}^{2} / \partial
{\sigma}_{\nu}^{2}=-4 {X}_{\mu} {\sf{\Sigma}}^{W}_{\mu \nu}
{X}_{\nu}/{\lambda}_{1}=-4 {X}_{\mu} {\sf{\Sigma}}^{R}_{\mu \nu}
{X}_{\nu}/{\lambda}_{2}$.  Finally, the constraints in
Eq.~(\ref{44}) together with the symmetry of the CSM imply that $\bf
W$ and $\bf R$ are null vectors for the three matrices
${{X}_{\mu}}^{-1}(\partial {{X}_{\mu}}^{2} / \partial
{\sigma}_{\nu}^{2}){{X}_{\nu}}^{-1}$, ${\sf{\Sigma}}^{W}_{\mu \nu}$,
and ${\sf{\Sigma}}^{R}_{\mu \nu}$.  Therefore the above $3M\times3M$
matrix is in fact redundant, with a rank no larger than $M-2$.  It
is also evident by an inspection of Eq.~(\ref{44}) that
$\mathbf{X}(\vec{{\sigma}^2},\mathbf{W},{\cal W},\mathbf{R},{\cal
R})$ is homogeneous of order zero in the parameters
$\vec{{\sigma}^2}$, $(\bar{\mathbf{W}},{\cal W})$, and
$(\bar{\mathbf{R}},{\cal R})$ separately.  Thus the intuition
associated with the Slutsky-Hicks problem can be brought to bear on
the portfolio problem as formulated in Eq.~(\ref{44}) in its
entirety (albeit with a trivial reversal of signs). Note that the
negative semidefiniteness of $\partial {{X}_{\mu}}^{2} / \partial
{\sigma}_{\nu}^{2}$ implies that the portfolio variance is a concave
function of the principal variances, a result that is evident in
Eq.~(\ref{44}) by inspection and may also be deduced from
Eq.~(\ref{44}) by an application of the envelope relations.  In
summary, then, we have

\textbf{Property (vii)}\textit{The comparative statics properties of
the principal portfolio analysis are summarized in three
semidefinite matrices of order $M$ and rank no larger than $M-2$ as
given in Eq.~(\ref{47}) et seq.}

At this juncture the results of the principal portfolio analysis can be
related to the original form.  First, let us write
\begin{equation}
\mathbf{x}(\vec{{\sigma}^2},\vec{\mathbf{e}},\mathbf{w},{\cal W},
\mathbf{r}, {\cal R})={\sum}_{\mu=1}^{M}\mathbf{e}^{\mu}{X}_{\mu}
(\vec{{\sigma}^2},\mathbf{W},{\cal W},\mathbf{R},{\cal R}),
\label{48}
\end{equation}
where $\vec{\mathbf{e}} \stackrel{\rm{ def}}{=}(\mathbf{e}^{1},
\mathbf{e}^{2},\ldots,\mathbf{e}^{M})$ is a vector of unit vectors.
Note that we have replaced the covariance matrix by its principal
variance ``vector'' $\vec{{\sigma}^2}$ and the vector of its
eigenvectors, $\vec{\mathbf{e}}$, among the arguments of
$\mathbf{x}$. One measure of the progress achieved in going from
$\bf x$ to $\bf X$ is the reduction in the number of arguments
needed to describe these functions.  This number is of the order of
${M}^{2}/2$, as can be estimated from Eq.~(\ref{48}).  From a
numerical point of view, the price for this progress is the cost of
computing the eigenquantities of the covariance matrix (a standard,
albeit nontrivial, task).

If desired, the comparative statics of the problem deduced above can
be translated into the variables of the original formulation using
Eq.~(\ref{48}).  It is important to realize, however, that except
for the pair of Slutsky matrices, the clarity and compactness of the
above results will be obscured when stated in terms of the original
variables.  For this reason, we will forgo a detailed calculation
here and limit the discussion to a statement of the comparative
statics information with respect to $(\bf w,{\cal W})$ and $(\bf
r,{\cal R})$.  That information is already discernible from the
original form in Eq.~(\ref{43}) by analogy to the linear,
multiple-constraint utility maximization model, and it can also be
deduced from the matrices ${\sf{\Sigma}}^{W}$ and
${\sf{\Sigma}}^{R}$.  In summary, we have

\textbf{Property (vii)} \textit{The Slutsky matrices ${\lambda}_{1}
{\sf{\Sigma}}^{w}$ and ${\lambda}_{2} {\sf{\Sigma}}^{r}$, where
${\sf{\Sigma}}^{w}_{ij} \stackrel{\rm{def}}{=}\partial {x}_{i}/
\partial {w}_{j} +{x}_{j}\partial {x}_{i}/ \partial {\cal W}$ and
${\sf{\Sigma}}^{r}_{ij}\stackrel{\rm{def}}{=}\partial {x}_{i}/
\partial {r}_{j} +{x}_{j}\partial {x}_{i} / \partial {\cal R} $, are
positive semidefinite and related by the equation ${\lambda}_{2}
{\sf{\Sigma}}^{w}={\lambda}_{1} {\sf{\Sigma}}^{r}$.  Each has a rank
no larger than $M-2$, and possesses two null vectors, $\mathbf{w}$
and $\bf r$. Furthermore, $\mathbf{x}(\vec{{\sigma}^2},
\vec{\mathbf{e}},\mathbf{w},{\cal W}, \mathbf{r}, {\cal R})$ is
homogeneous of order zero in the parameters $\vec{{\sigma}^2}$,
$(\bf w,{\cal W})$, and $(\mathbf{r}, {\cal R})$ separately.}

Of course there is no direct interest attached to ${\sf{\Sigma}}^{w}$
since the parameters $(\bf w,{\cal W})$ are fixed by definition.  We
conclude our discussion of the portfolio problem (with short sales
allowed) by noting that the foregoing analysis has amply confirmed our
earlier statements that the principal
portfolios are the natural variables for analyzing the efficient portfolio
problem.  This is of course another instance of a golden rule in applied
analysis that, where quadratic forms are involved, a reformulation of the
problem in terms of the principal axes is often quite advantageous.

\section{Concluding Remarks}

The main objective of this work, namely the derivation of unrestricted
comparative statics matrices for a general, differentiable optimization
problem, has been realized.  The result, which is stated in Theorem 1, in
effect completes the program initiated and developed by Samuelson (1947),
generalized and streamlined by Silberberg (1974), and further advanced by
Hatta (1980) and other authors, most of whom were mentioned in the
introduction.  We have primarily established our results constructively,
devoting considerable effort to explaining the details of the
construction, as well as to developing a clear, geometric picture of its
workings.  We have also developed a number of new results and extensions,
mainly summarized in Theorem 2 and the corollaries to Theorem 1, which
further characterize the properties of the CSM's and thereby serve to
broaden the power and reach of the analysis.  In particular, the
realization that comparative statics results for a given problem can
assume a wide range of shapes and forms significantly strengthens their
role in empirical comparisons and hypothesis testing, the primary raison
d'\^{e}tre for a CSM.

In \S IID we constructed a maximal, universal CSM for an arbitrary
optimization problem, summarized in Theorem 3.  For theoretical purposes,
this theorem is an important basic result as it defines a framework in
which any other method must be subsumed.  Not surprisingly, this universal
method is not as practically convenient, or intuitively appealing, as the
method of generalized compensated derivatives used throughout this work.
On the other hand, recalling that the universal construction is applicable
to any differentiable system governed by an extremum principle, a category
which includes numerous physical and mathematical problems in quite
diverse fields of inquiry, one realizes the vast scope and reach of this
theorem.

Our applications cover a wide range from the very familiar to the novel,
and are intended to illustrate the theorems and demonstrate the
effectiveness of the method.  Remarkably, a surprising number of new and
significant results have emerged from these applications (as summarized in
italicized paragraphs in \S III), even in the case of very familiar
models.  In choosing the applications, we have strived to illustrate our
method by means of meaningful and interesting economic models, avoiding
excessively complex problems which, although more effective in
demonstrating the power of the formalism, would be lacking in clear
economic meaning or intuitive sense.

Throughout, we have dealt with the comparative statics of a given,
interior solution to an optimization problem.  As such, there is no need
to deal with inequality constraints, since those that bind can be treated
as equality constraints, and those that don't can be ignored altogether.
Nor have we concerned ourselves with issues of integrability, since these
are primarily relevant to utility maximization problems of a particular
structure.  Similarly, although we have not dealt with problems involving
discrete-time, finite-horizon, intertemporal
optimization, this and other categories can be treated
straightforwardly by our method.

\textbf{Allen, Roy G.D.} \textit{Mathematical Analysis for
Economists.} London: Macmillan and Co., Limited, 1938.

\textbf{Antonelli, Giovanni B.} ``Sulla Teoria Matematica della
Economia Politica,'' monograph published in \textit{Pisa, Nella
Tipogrefia del Folchetto}, 1886, translated as ``On the Mathematical
Theory of Political Economy," in  J. S. Chipman, L. Hurwicz, M. K.
Richter, and H. F. Sonnenschein, eds., \textit{Preferences, Utility,
and Demand}.  New York: Harcourt Brace Jovanovich, Inc. 1971.

\textbf{Chichilnisky, Graciela and Kalman, Peter J.} ``Properties of
Critical Points and Operators in Economics.'' \textit{Journal of
Mathematical Analysis and Applications}, February 1977, 57(2), pp.
340-49.

\textbf{Chichilnisky, Graciela and Kalman, Peter J.} ``Comparative
Statics of Less Neoclassical Agents.'' \textit{International
Economic Review}, February 1978, 19(1), pp. 141-48.

\textbf{Cournot, Augustin A.} \textit{Researches into the
Mathematical Principles of the Theory of Wealth,} 1838, translated
by Nathaniel T. Bacon, New York: The Macmillan Company, 1897.

\textbf{Elton, Edwin J. and Gruber, Martin J.} \textit{Modern
Portfolio Theory and Investment Analysis}, 4th Ed. New York: Wiley,
1991.

\textbf{Fama, Eugene F. and Miller, Merton H.} \textit{The Theory of
Finance}. Hinsdale: Dreyden Press, 1972, p. 282.

\textbf{Hatta, Tatsuo.} ``Structure of the Correspondence Principle
at an Extremum Point.'' \textit{Review of Economic Studies}, October
1980, 47(5), pp. 987-97.

\textbf{Hicks, John R.} \textit{Value and Capital}. London: Oxford
University Press, 1939.

\textbf{Houthakker, Hendrick S.} ``Compensated Changes in Quantities
and Qualities Consumed.'' \textit{Review of Economic Studies},
1951-52, 19, pp. 155-64.

\textbf{Kalman, Peter and Intriligator, Michael D.} ``Generalized
Comparative Statics with Applications to Consumer and Producer
Theory.'' \textit{International Economic Review}, June 1973, 14(2),
pp. 473-86.

\textbf{Mas-Colell, Andreu, Whinston, Michael D., and Green, Jerry
R.} \textit{Microeconomic Theory.} New York: Oxford University
Press, 1995, pp. 956-58.

\textbf{Milgrom, Paul and Shannon, Chris.} ``Monotone Comparative
Statics.'' \textit{Econometrica}, January 1994, 62(1), pp. 157-80.

\textbf{Milgrom, Paul and Roberts, John.} ``Comparing Equilibria.''
\textit{American Economic Review}, June 1994, 84(3), pp. 441-59.

\textbf{Paris, Quirino, Caputo, Michael R. and Holloway, Garth J.}
``Keeping the Dream of Rigorous Hypothesis Testing Alive.''
\textit{American Journal of Agricultural Economics}, October 1993,
75(6), pp. 25-40.

\textbf{Samuelson, Paul A.} \textit{Foundations of Economic
Analysis}. Cambridge: Harvard University Press, 1947.

\textbf{Samuelson, Paul A.} ``Using Full Duality to Show That
Simultaneously Additive Direct and Indirect Utilities Implies
Unitary Price Elasticity of Demand.'' \textit{Econometrica}, October
1965, 33(4), pp. 781-96.

\textbf{Silberberg, Eugene.} ``A Revision of Comparative Statics
Methodology in Economics, or, How to Do Comparative Statics on the
Back of an Envelope.'' \textit{Journal of Economic Theory}, February
1974, 7(2), pp. 159-72.

\textbf{Silberberg, Eugene.} \textit{The Structure of Economics: A
Mathematical Analysis}, 2nd Ed. New York: McGraw-Hill Publishing
Company, 1990, pp. 213 and 216.

\textbf{Slutsky, E. E.} ``Sulla Teoria del Bilancio del
Consumatore.'' \textit{Giornale degli Economisti e Rivista di
Statistica}, 1915, 51, 1-26, translated as ``On the Theory of the
Budget of the Consumer,'' in G. L. Stigler and K. Boulding, eds.,
\textit{Readings in Price Theory}, Homewood: Richard Irwin, 1952.

\textbf{Takayama, Akira.} \textit{Mathematical Economics}, 2nd Ed.
Cambridge: Cambridge University Press, 1985, p. 159.

\end{document}